\documentclass[11pt,twoside, leqno]{article}

\usepackage{amssymb}
\usepackage{amsmath}
\usepackage{mathrsfs}
\usepackage{amsthm}
\usepackage{amsfonts}
\usepackage{txfonts}
\usepackage{enumerate}
\usepackage{hyperref}
\usepackage{esint}
\usepackage{latexsym,amssymb}
\usepackage{color}
\usepackage{marginnote}
\usepackage{multirow}
\usepackage{makecell}
\usepackage{float}
\usepackage{microtype}

\usepackage{yhmath}

\usepackage{indentfirst}
\usepackage{graphicx}
\usepackage{pgf,tikz}
\usepackage{mathrsfs}
\usetikzlibrary{arrows}

\definecolor{ttttff}{rgb}{0.2,0.2,1.}
\definecolor{ttffcc}{rgb}{0.2,1.,0.8}
\definecolor{qqqqff}{rgb}{0.,0.,1.}

\definecolor{zzttqq}{rgb}{0.6,0.2,0.}
\definecolor{qqqqff}{rgb}{0.,0.,1.}

\allowdisplaybreaks

\hypersetup{
 colorlinks = true,
 citecolor = blue,
 linkcolor = blue,
 urlcolor = blue
}

\def\rr{{\mathbb R}}
\def\rn{{\mathbb{R}^n}}

\def\zz{{\mathbb Z}}

\def\nn{{\mathbb N}}

\def\cm {{\mathcal{M}}}

\def\fz{\infty }

\def\lf{\left}
\def\r{\right}

\def\gfz{\genfrac{}{}{0pt}{}}

\def\msn{{\mathbb{S}^{n-1}}}

\def\dis{\displaystyle}

\def\loc{{\mathop\mathrm{loc\,}}}
\def\supp{\mathop\mathrm{supp\,}}

\def\uc{\varepsilon}

\newtheorem{theorem}{Theorem}[section]
\newtheorem{lemma}[theorem]{Lemma}

\newtheorem{proposition}[theorem]{Proposition}
\newtheorem{example}[theorem]{Example}
\theoremstyle{definition}
\newtheorem{remark}[theorem]{Remark}
\newtheorem{definition}[theorem]{Definition}
\renewcommand{\appendix}{\par
\setcounter{section}{0}%
\setcounter{subsection}{0}%
\setcounter{subsubsection}{0}%
\gdef\thesection{\@Alph\c@section}%
\gdef\thesubsection{\@Alph\c@section.\@arabic\c@subsection}%
\gdef\theHsection{\@Alph\c@section.}%
\gdef\theHsubsection{\@Alph\c@section.\@arabic\c@subsection}%
\csname appendixmore\endcsname
}

\numberwithin{equation}{section}

\definecolor{qqqqff}{rgb}{0,0,1}
\definecolor{ffqqqq}{rgb}{1,0,0}

\newcounter{rea}
\def\MSq{\hyperref[def-MScond]{{\rm (MS)}_{q}}}
\def\NMSq{\hyperref[def-MScond]{{\rm (NMS)}_{q}}}
\def\VMSq{\hyperref[def-MScond]{{\rm (VMS)}_{q}}}

\def\BBMmq{\hyperref[def-BBMcond]{{\rm (BBM)}_{m,\,q}}}
\def\NBBMmq{\hyperref[def-BBMcond]{{\rm (NBBM)}_{m,\,q}}}
\def\VBBMmq{\hyperref[def-BBMcond]{{\rm (VBBM)}_{m,\,q}}}

\def\mathringBVm{\smash{\mathring{\mathrm{BV}}}^m}
\def\BVm{\mathrm{BV}^m}

\begin{document}

\arraycolsep=1pt

\title{\vspace{-1.8cm}
Pointwise and $\Gamma$-Convergence of BBM and MS formulae for Nonlocal Besov Type Energies
\footnotetext{\hspace{-0.35cm}
2020 {\it Mathematics Subject Classification}. Primary: 46E35;
Secondary: 49J45, 42B35.
\endgraf {\it Key words and phrases}. $\Gamma$-convergence, BBM formula, MS formula, Besov space, logarithmic Besov space, Sobolev space.
\endgraf
L. Liu is supported by the National Natural Science Foundation of China (\# 12371102);
Q. Zhang is supported by the Fundamental Research Funds for the Central Universities, and the Research Funds of Renmin University of China (\# RUC25QSDL132).}}
\author{Liguang Liu and Qishun Zhang}
\date{}
\maketitle


\begin{abstract}
For $p\in[1,\infty)$, $q\in(0,\infty)$ and $m\in\mathbb N$, we study
nonlocal Besov type energies $\mathscr{E}_{p,\,q}^{m,\,\varphi_j}$
defined by $m$th-order differences and nonnegative kernels
$\varphi_j$. We establish necessary and
sufficient kernel conditions for the Maz'ya--Shaposhnikova (MS) limit when
$m=1$ and for the existence of subsequential Bourgain--Brezis--Mironescu (BBM)
limits for arbitrary $m$. In each limiting regime, we prove that the same kernel condition
characterizes pointwise convergence, $\Gamma$-convergence in the strong
$L^p$ topology, and strong $\Gamma$-convergence with constant recovery
sequences.

For $m=1$, the MS limit of $\{\mathcal E^{1,\,\varphi_j}_{p,\,q}(f)\}_{j\in\nn}$ for $f\in W^{1,\,p}(\mathbb R^n)$
is $2^{q/p}\|f\|_{L^p(\mathbb R^n)}^q$.
Each mode of convergence holds along the full
sequence if and only if the weighted measures
$\min\{|h|^q,\,1\}\varphi_j(h)\,dh$ have mass tending to one outside
$B_1$ and to zero in $B_\delta$ for every $\delta\in[1,\infty)$, where $B_\delta$ denotes the ball of radius $\delta$ centered at the origin.

For arbitrary $m$, the existence of a subsequential
BBM limit in any one of the three senses is
equivalent to the existence of a subsequence along which $\min\{|h|^{mq},\,1\}\varphi_j(h)\,dh$ have mass tending to
one in $B_1$ and to zero outside $B_\delta$ for every $\delta\in(0,1]$.
The limit on $W^{m,\,p}(\mathbb R^n)$
averages the $q$th powers of the $L^p$ norms of the $m$th-order
directional derivatives against a Borel probability measure $\mu$ on
the unit sphere. At $p=1$, the characterization extends to functions
of bounded variation of order $m$, with distributional directional
derivatives measured by total variation. For radial kernels, the
characterization holds along the full sequence with $\mu$ equal to
normalized surface measure. A counterexample shows that passage to a
subsequence can be necessary for nonradial kernels.
\end{abstract}

\section{Introduction and main results}\label{sec1}

The classical Bourgain--Brezis--Mironescu (BBM) and Maz'ya--Shaposhnikova (MS) formulae \cite{Bourgain-Brezis-Mironescu-2001,
Mazya-Shaposhnikova-JFA-2002} describe the pointwise asymptotic limits of fractional Sobolev energies as the differentiation
order $s\in(0,1)$ approaches the two endpoints of the scale. The BBM limit ($s\to1_-$) describes the localization of nonlocal
energies to a local gradient energy, whereas the MS limit ($s\to0_+$) describes the opposite regime where the limiting energy
reduces to a pure mass term.

The $\Gamma$-convergence framework provides a complementary variational perspective: through liminf inequalities and recovery
sequences, it identifies variational limits and, in the presence of suitable compactness, guarantees the convergence of minima
and the stability of minimizers. In general, pointwise convergence and $\Gamma$-convergence are distinct notions, and neither
implies the other \cite[Theorems 1 and 2]{Nguyen-DUKE-2011}.

This paper proves, however, that for a broad class of nonlocal energies -- allowing arbitrary difference orders in the BBM regime
and first difference order in the MS regime, general exponents, and nonradial kernels -- sharp kernel conditions simultaneously characterize
pointwise convergence, $\Gamma$-convergence, and strong $\Gamma$-convergence in both the BBM and MS regimes.


\subsection{Pointwise convergence of BBM and MS formulae}\label{ss1.1}

We begin by recalling the classical BBM and MS formulae.
Let $p\in[1,\fz)$ and $s\in(0,1)$. Define the \emph{homogeneous fractional Sobolev space} $\mathring{W}^{s,\,p}(\rn)$
to be the collection of all locally integrable functions $f$ on $\rn$ such that
\begin{align*}
\|f\|_{\mathring{W}^{s,\,p}(\rn)}
:= \lf(\int_\rn \int_\rn \frac{|f(x)-f(y)|^p}{|x-y|^{n+sp}}\,dy\,dx \r)^\frac{1}{p}<\fz.
\end{align*}
Let $W^{1,\,p}(\rn)$ be the \emph{classical inhomogeneous Sobolev space}
of order $1$.
The Bourgain--Brezis--Mironescu formula \cite{Bourgain-Brezis-Mironescu-2001} asserts that,
for every $f\in W^{1,\,p}(\rn)$,
\begin{align}\label{eq-BBM-Wsp}
\lim_{s\to1_-} (1-s)p\,\|f\|_{\mathring{W}^{s,\,p}(\rn)}^p
= \int_\msn \|\xi\cdot\nabla f\|_{L^p(\rn)}^p\,d\sigma(\xi),
\end{align}
where $d\sigma$ denotes the surface measure on the unit sphere $\msn$ in $\rn$.
For the endpoint $p=1$, by the argument of D\'avila \cite{Davila-CVPDE-2002},
identity \eqref{eq-BBM-Wsp} holds for all functions $f$ with bounded variation on $\rn$, with the norm
$\|\xi\cdot\nabla f\|_{L^1(\rn)}$ replaced by the total variation $|\xi\cdot\nabla f|(\rn)$. In contrast,
the Maz'ya--Shaposhnikova formula \cite[Theorem 3]{Mazya-Shaposhnikova-JFA-2002} states
that if $f\in L^p(\rn)\cap\mathring{W}^{s_0,\,p}(\rn)$ for some $s_0\in(0,1)$, then
\begin{align}\label{eq-MS-Wsp}
\lim_{s\to0_+} sp\,\|f\|_{\mathring{W}^{s,\,p}(\rn)}^p
= 2\omega_{n-1}\,\|f\|_{L^p(\rn)}^p,
\end{align}
where $\omega_{n-1}=\sigma(\msn)$.
Identifying the two endpoint limits of the fractional Sobolev scale, the BBM and MS formulae
have since served as prototypes for numerous extensions.

Nowadays, the BBM and MS formulae have attracted considerable attention and have undergone intensive development.
For example, Milman \cite{Milman-TransAMS-2005} gave a new proof of \eqref{eq-BBM-Wsp} and \eqref{eq-MS-Wsp} via
real interpolation and, more importantly, established the continuity of interpolation scales at the endpoints, which
paves the way for extending the BBM and MS formulae to general scales of function spaces. This is evidenced by Karadzhov,
Milman and Xiao \cite{Karadzhov-Milman-Xiao-JFA-2005}, who established the BBM and MS formulae for various higher-order
Besov spaces $\mathring B^{sm}_{p,\,q}(\rn)$ as $s\to1_-$ or $s\to0_+$, with endpoint limits in $\mathring{W}^{m,\,p}(\rn)$
or $L^p(\rn)$, respectively. For further developments along this line, we refer the reader to \cite{Astashkin-Lykov-Milman-JAT-2019,
Dai-Grafakos-Pan-Yang-Yuan-Zhang-MathAnn-2024,Dominguez-Milman-CVPDE-2023,Leoni-Spector-JFA-2011,Liu-Xiao-Zhang-2026-submit,
Liu-Zhang-2026-submit} and the references therein.

Motivated by these developments, we develop a unified framework for studying the BBM and MS formulae beyond the
classical fractional Sobolev setting. We consider a general class of nonlocal energies based on
differences of order $m$, with general exponents $p,\,q$ and nonnegative locally integrable kernels $\varphi$.
The corresponding local energies are described by directional Sobolev-type functionals, which may be
anisotropic and are parameterized by Borel probability measures on the unit sphere. We formulate this
framework in the following definition.

\begin{definition}\label{def-energy}
Let $p\in[1,\fz)$, $q\in(0,\fz)$ and $m\in\nn$.
\begin{enumerate}[{\rm (i)}]
\item (\emph{nonlocal energy})
Let $\varphi$ be a nonnegative locally integrable function on $\rn \setminus\{0\}$, with $\varphi(0)$ possibly equal to $+\infty$.
For any $f\in L_{\rm loc}^1(\rn)$,
define the \emph{Besov energy} by
\begin{align}\label{eq-energy}
\mathscr{E}_{p,\,q}^{m,\,\varphi}(f)
:=\int_\rn \|\Delta_h^mf\|_{L^p(\rn)}^q\,\varphi(h)\,dh,
\end{align}
where $\Delta_h^m$ denotes the \emph{$m$-order difference}, given by
\begin{align*}
\Delta_hf(x):=\Delta_h^1f(x):=f(x+h)-f(x)
\qquad\text{for}\ \, x\in\rn,
\end{align*}
and
\begin{align*}
\Delta_h^mf:=\Delta_h\bigl(\Delta_h^{m-1}f\bigr)
\qquad\text{for}\ \, m\ge2.
\end{align*}
When $m=1$, we simply write $\mathscr{E}_{p,\,q}^{\varphi}(f):= \mathscr{E}_{p,\,q}^{1,\,\varphi}(f)$.

\item (\emph{local energy}) Let $\mu$ be a Borel probability measure on $\msn$. For any $f\in L_{\rm loc}^1(\rn)$,
define the \emph{anisotropic $\mu$-Sobolev energy} by
\begin{align}\label{eq-Dpqmf}
\mathscr{D}_{p,\,q}^{m,\,\mu}(f)
:=\int_\msn \|(\xi\cdot\nabla)^mf\|_{L^p(\rn)}^q\,d\mu(\xi)
\end{align}
and the \emph{anisotropic $\mu$-bounded variation energy} by
\begin{align}\label{eq-Dqmf-BV}
\mathscr{D}_q^{m,\,\mu}(f)
:= \int_\msn \|(\xi\cdot\nabla)^mf\|_{\mathcal{M}(\rn)}^q\,d\mu(\xi),
\end{align}
where $\mathcal{M}(\rn)$ is the space of finite signed Radon measures with the total variation norm. In
both local energies a directional norm is assigned the value $\fz$ if the distributional derivative does
not belong to the indicated space. In particular, if $d\mu=(\omega_{n-1})^{-1}\,d\sigma$, then we omit
the superscript $\mu$, and write $\mathscr{D}_{p,\,q}^{m,\,\mu}$ and $\mathscr{D}_q^{m,\,\mu}$ simply as
$\mathscr{D}_{p,\,q}^m$ and $\mathscr{D}_q^m$, respectively.
\end{enumerate}
\end{definition}

For any $s\in(0,\infty)$, taking $m=\lfloor s\rfloor+1$ and $\varphi(h)=|h|^{-(n+sq)}$ in \eqref{eq-energy}
yields the classical \textit{homogeneous Besov space $\mathring{\Lambda}_{p,\,q}^s(\rn)$}, whose
(quasi-)seminorm is given by
\begin{align*}
\|f\|_{\mathring{\Lambda}_{p,\,q}^s(\rn)}
:= \lf(\mathscr{E}^{m,\,\varphi}_{p,\,q}(f) \r)^\frac{1}{q}
= \lf(\int_\rn \|\Delta_h^{\lfloor s\rfloor+1}f\|_{L^p(\rn)}^q\,|h|^{-(n+sq)}\,dh \r)^\frac{1}{q}.
\end{align*}
We refer the reader to Triebel's monographs \cite{Triebel-function1-1983,Triebel-functionspace2-1992} for a systematic
treatment of Besov spaces. Another important example of \eqref{eq-energy} is obtained by considering the logarithmic kernel
\begin{align*}
\varphi(h):=|h|^{-(n+mq)} \lf(\log\frac{e}{|h|\wedge1} \r)^{bq},
\end{align*}
where $bq<-1$. This choice leads to the \textit{logarithmic Besov spaces} studied in \cite{Haroske-ZAnalA-2000,
Dominguez-Haroske-Tikhonov-JMPA-2020,Dominguez-Tikhonov-AMS-2023} and the references therein. Further examples
and related discussions will be given in
Section \ref{sec-example}.

The local energies defined in \eqref{eq-Dpqmf} and \eqref{eq-Dqmf-BV} are closely related to the homogeneous
Sobolev space $\mathring{W}^{m,\,p}(\rn)$ and the higher-order bounded variation space $\mathringBVm(\rn)$, respectively;
see Section \ref{sec-Perliminary1} for details. Indeed, for every Borel probability measure $\mu$ on $\msn$,
the energy $\mathscr{D}_{p,\,q}^{m,\,\mu}$ is finite on $\mathring{W}^{m,\,p}(\rn)$, while $\mathscr{D}_{q}^{m,\,\mu}$
is finite on $\mathringBVm(\rn)$. Moreover, when $\mu$ is the normalized surface measure on $\msn$, we have
\begin{align*}
\lf(\mathscr{D}_{p,\,q}^m(f) \r)^\frac{1}{q}
\approx \|f\|_{\mathring{W}^{m,\,p}(\rn)}
\qquad\text{and}\qquad
\lf(\mathscr{D}^m_{q}(f) \r)^\frac{1}{q}
\approx \|f\|_{\mathringBVm(\rn)},
\end{align*}
where the implicit constants depend only on $m$, $n$ and $q$; see Proposition \ref{Prop-equivalent-Sobolevnorm}.

We remark that the kernel $\varphi$ in \eqref{eq-energy} need not be radial,
which corresponds to direction-dependent interactions and gives rise to
anisotropic local limits of the form \eqref{eq-Dpqmf} or \eqref{eq-Dqmf-BV}.
Nonradial kernels arise naturally in the variational theory of micromagnetism modeling  on antisymmetric exchange interactions (see \cite{Davoli-DiFratta-Giorgio2024JAM, Davoli-DiFratta-Pagliari-ProcRSEC-2026}), and in nonlocal elliptic equations (see
\cite{Caffarelli-Silvestre-CPAM-2009,Caffarelli-Silvestre-Annals-2011}); the
resulting directional structure is encoded in the spherical measure $\mu$ of
the limiting energy (see \cite{Gennaioli-Stefani-ARMA-2026} for more details).
We also note that second-order differences occur naturally in fully nonlinear nonlocal elliptic equations; see, for example,
\cite{Caffarelli-Silvestre-Annals-2011}.

The classical BBM formula \eqref{eq-BBM-Wsp} and MS formula \eqref{eq-MS-Wsp} arise as the two endpoint limits, as
$s\to1_-$ and $s\to0_+$, respectively, of the same fractional Sobolev energy $\mathscr{E}_{p,\,p}^{1,\,\varphi_s}$, where
$\varphi_s(h)=|h|^{-(n+sp)}$. In our previous work \cite{Liu-Zhang-2026-submit}, we established BBM and MS formulae
for the much broader family of energies $\mathscr{E}_{p,\,q}^{m,\,\varphi}$, allowing arbitrary difference orders $m$
and general exponents $p$ and $q$, under suitable \emph{sufficient} conditions on \emph{radial} kernels. The results
most relevant to the present paper are summarized below.

\begin{theorem}[\cite{Liu-Zhang-2026-submit}]\label{thm-LiuZhang}
Let $p\in[1,\fz)$, $q\in(0,\fz)$, $m\in\nn$ and $\{\varphi_j\}_{j\in\nn}$ be a sequence of nonnegative radial
functions in $L_{\rm loc}^1(\rn\setminus\{0\})$. Then the following implications hold:
\begin{enumerate}[{\rm (i)}]
\item \textbf{(MS formula)} Assume that the sequence $\{\varphi_j\}_{j\in\nn}$ satisfies the following $\MSq$ condition:
\begin{align}\label{def-MScond}
\begin{cases}
 \dis \lim_{j\to\fz} \int_{|h|\ge 1} \varphi_j(h)\, dh=1; \hfill \NMSq \\[4mm]
 \dis \lim_{j\to\fz} \int_{|h|<\delta} \varphi_j(h) \min\lf\{|h|^q,\,1 \r\} dh=0
\ \ \text{for any}\ \,\delta\in[1,\fz). \quad \hfill \VMSq
\end{cases}\tag*{$\MSq$}
\end{align}
Then, for all $f\in W^{1,\,p}(\rn)$,
\begin{align*}
\lim_{j\to\fz} \mathscr{E}_{p,\,q}^{\varphi_j}(f)
= 2^\frac{q}{p}\,\|f\|_{L^p(\rn)}^q.
\end{align*}

\item \textbf{(BBM formula)} Assume that the sequence $\{\varphi_j\}_{j\in\nn}$ satisfies the following $\BBMmq$ condition:
\begin{align}\label{def-BBMcond}
\begin{cases}
\dis \lim_{j\to\fz} \int_{|h|<1} |h|^{mq}\,\varphi_j(h)\,dh=1;\hfill \NBBMmq \\[4mm]
\dis \lim_{j\to\fz} \int_{|h|\ge\delta} \varphi_j(h) \min\lf\{|h|^{mq},\,1 \r\} dh=0
\ \ \text{for any}\ \,\delta\in(0,1]. \quad \hfill \VBBMmq
\end{cases}\tag*{$\BBMmq$}
\end{align}
Then, for all $f\in W^{m,\,p}(\rn)$,
\begin{align}\label{eq-BBM}
\lim_{j\to\fz} \mathscr{E}_{p,\,q}^{m,\,\varphi_j}(f)
= \mathscr{D}_{p,\,q}^m(f).
\end{align}
\end{enumerate}
\end{theorem}

\begin{remark}\label{eq-N-V-conds}
The conditions $\MSq$ and $\BBMmq$ were introduced in \cite{Liu-Zhang-2026-submit} for radial kernel sequences.

\begin{enumerate}[\rm (i)]
\item For later use, we refer to the first and second conditions in $\MSq$ as the \emph{normalization condition} $\NMSq$ and the
\emph{vanishing condition} $\VMSq$, respectively. Note that $\NMSq$ normalizes the total mass outside the unit ball, while
$\VMSq$ requires the weighted mass on every fixed ball to vanish.

\item Similarly, we call the first and second conditions in $\BBMmq$ the \emph{normalization condition} $\NBBMmq$
and the \emph{vanishing condition} $\VBBMmq$, respectively. In contrast to $\MSq$, the normalization and vanishing
conditions in $\BBMmq$ interchange the roles of the origin and infinity: $\NBBMmq$ normalizes the weighted mass near
the origin, while $\VBBMmq$ requires the weighted mass away from the origin to vanish.
\end{enumerate}
\end{remark}

Regarding Theorem \ref{thm-LiuZhang}(ii), Foghem \cite[Section 9.2]{Foghem-NoDEA-2025} proved that, in the
special case $m=1$, $p=q\in(1,\fz)$, and for radial kernel sequences, the condition $\BBMmq$ is not only sufficient
but also necessary for the validity of the BBM formula \eqref{eq-BBM}. These results naturally raise the following question:
\begin{center}
\begin{minipage}{13cm}
\emph{For general exponents $p$ and $q$, and nonradial kernels $\{\varphi_j\}_{j\in\nn}$, are the conditions
$\MSq$ and $\BBMmq$ necessary and sufficient for the corresponding pointwise MS and BBM formulae?}
\end{minipage}
\end{center}
We answer this question in the affirmative; see Theorems \ref{thm-MS} and \ref{thm-BBM} below.
For the BBM case, the characterization is established up to passing to a suitable subsequence.

\subsection{$\Gamma$-Convergence of BBM and MS formulae}\label{ss1.2}

The classical notion of $\Gamma$-convergence, introduced by De Giorgi and Franzoni \cite{Giorgi-Franzoni-AttiAccadNazLincei-1975},
plays an important role in establishing the local limits of sequences of nonlocal energies (see, for example, \cite{Kumagai-Sturm-JMKU-2005,
Sturm-PotentialAnal-1998,AlonsoRuiz-Baudoin-Pisa-2025,AlonsoRuiz-Baudoin-NonlinearAnal-2025,Crismale-DeLuca-Kubin-Ninno-Ponsiglione-JFA-2023}).
Nonlocal functionals related to the total variation have also been investigated from the $\Gamma$-convergence viewpoint, with
connections to image processing \cite{Brezis-Nguyen-AnnPDE-2018}.
In nonlinear elasticity, Friesecke, James and M\"uller derived a hierarchy of plate models by $\Gamma$-convergence
\cite{Friesecke-James-Muller-ARMA-2006}.




Following the classical liminf and recovery-sequence formulation of $\Gamma$-convergence in \cite{Giorgi-Franzoni-AttiAccadNazLincei-1975},
we introduce a notion of $\Gamma$-convergence for the energy sequence $\{\mathscr{E}_{p,\,q}^{m,\,\varphi_j}\}_{j\in\nn}$
on the function space $X\subset L^{p}(\rn)$ and with respect to the strong $L^p(\rn)$ topology.

\begin{definition}\label{def-Gamma}
Let $p\in[1,\fz)$, $q\in(0,\fz)$, $m\in\nn$ and $\{\varphi_j\}_{j\in\nn}$ be a sequence of nonnegative
functions in $L_{\mathrm{loc}}^1(\rn\setminus\{0\})$. Given a function space $X\subset L^p(\rn)$, let
$ \mathscr{D}:\ X\to[0,\fz]$ be a functional. Suppose that the following two conditions hold:
\begin{itemize}
\item[$\mathbf{(\Gamma)_1}$]
for any $f\in X$, there exists a sequence $\{f_j\}_{j\in\nn}\subset L^p(\rn)$ such that
$f_j\to f$ in $L^p(\rn)$ as $j\to\fz$ and
\begin{align*}
\limsup_{j\to\fz}\mathscr{E}_{p,\,q}^{m,\,\varphi_j}(f_j)\le \mathscr{D}(f);
\end{align*}

\item[$\mathbf{(\Gamma)_2}$]
for any $f\in X$ and any sequence $\{f_j\}_{j\in\nn}\subset L^p(\rn)$ converging to
$f$ in $L^p(\rn)$,
\begin{align*}
\liminf_{j\to\fz}\mathscr{E}_{p,\,q}^{m,\,\varphi_j}(f_j)\ge \mathscr{D}(f).
\end{align*}
\end{itemize}
Then we say that the sequence $\{\mathscr{E}_{p,\,q}^{m,\,\varphi_j}\}_{j\in\nn}$ \emph{$\Gamma$-converges
to $ \mathscr{D}$ on $X$ with respect to the $L^p(\rn)$ topology}. In this case, we write
\begin{align*}
\Gamma\text{-}\lim_{j\to\fz}\mathscr{E}_{p,\,q}^{m,\,\varphi_j}
= \mathscr{D}
\quad\text{on}\,\ X.
\end{align*}
\end{definition}

In accordance with \cite{Davoli-DiFratta-Pagliari-ProcRSEC-2026,Gennaioli-Stefani-ARMA-2026}, we introduce the following
notion of \emph{strong $\Gamma$-convergence} for the energy sequence $\{\mathscr{E}_{p,\,q}^{m,\,\varphi_j}\}_{j\in\nn}$.

\begin{definition}\label{def-strong-Gamma}
Let $p$, $q$, $m$, $\{\varphi_j\}_{j\in\nn}$, $X$ and $ \mathscr{D}$ be as in Definition \ref{def-Gamma}.
If condition $\mathbf{(\Gamma)_2}$ of Definition \ref{def-Gamma} holds, and condition $\mathbf{(\Gamma)_1}$
is strengthened to
\begin{itemize}
\item[$\mathbf{(s\Gamma)_1}$]
for any $f\in X$,
\begin{align*}
\limsup_{j\to\fz}\mathscr{E}_{p,\,q}^{m,\,\varphi_j}(f)\le \mathscr{D}(f),
\end{align*}
\end{itemize}
then we say that the sequence $\{\mathscr{E}_{p,\,q}^{m,\,\varphi_j}\}_{j\in\nn}$ \emph{strongly $\Gamma$-converges
to $ \mathscr{D}$ on $X$ with respect to the $L^p(\rn)$ topology}. In this case, we write
\begin{align*}
s\Gamma\text{-}\lim_{j\to\fz}\mathscr{E}_{p,\,q}^{m,\,\varphi_j}
= \mathscr{D}
\quad\text{on}\,\ X.
\end{align*}
\end{definition}

\begin{remark}\label{rem-strong-Gamma}
It is clear that condition $\mathbf{(s\Gamma)_1}$ is precisely condition $\mathbf{(\Gamma)_1}$ with the constant
recovery sequence $f_j\equiv f$, which implies
$
\mathbf{(s\Gamma)_1}\ \Rightarrow\ \mathbf{(\Gamma)_1}
$
and, hence,
\begin{align*}
\text{strong $\Gamma$-convergence}\ \Rightarrow\ \text{$\Gamma$-convergence}.
\end{align*}
Further, taking $f_j\equiv f$ in $\mathbf{(\Gamma)_2}$, together with $\mathbf{(s\Gamma)_1}$, gives the pointwise convergence
\begin{align}\label{eq-pointwise}
\lim_{j\to\fz}\mathscr{E}_{p,\,q}^{m,\,\varphi_j}(f)= \mathscr{D}(f)
\quad\text{for every}\,\ f\in X,
\end{align}
so that
\begin{align*}
\text{strong $\Gamma$-convergence}\ \Rightarrow\ \text{pointwise convergence}.
\end{align*}
Conversely, \eqref{eq-pointwise} yields $\mathbf{(s\Gamma)_1}$, so that $\Gamma$-convergence together with
\eqref{eq-pointwise} implies strong $\Gamma$-convergence. Combining these implications, we conclude that
\begin{align*}
\text{strong $\Gamma$-convergence}\ \Leftrightarrow\
\text{$\Gamma$-convergence} \ +\ \text{pointwise convergence}.
\end{align*}
In general, $\Gamma$-convergence neither implies nor is implied by pointwise convergence. Interestingly, for
MS and BBM formulae, pointwise convergence turns out to be equivalent to $\Gamma$-convergence; see
Theorems \ref{thm-MS} and \ref{thm-BBM} below.
\end{remark}

\begin{remark}\label{rem-Gamma-cases}
We present the instances of Definition \ref{def-Gamma} to be used throughout this paper.
\begin{enumerate}[{\rm (i)}]
\item If $m=1$, $X=W^{1,\,p}(\rn)$ and $ \mathscr{D}= \mathscr{D}^0_{p,\,q}$, where
\begin{align*}
\mathscr{D}^0_{p,\,q}(f)
:= 2^\frac{q}{p}\,\|f\|_{L^p(\rn)}^q,
\end{align*}
then we say that \emph{the sequence $\{\mathscr{E}_{p,\,q}^{\varphi_j}\}_{j\in\nn}$ $\Gamma$-converges to
$\mathscr{D}^0_{p,\,q}$ on $W^{1,\,p}(\rn)$} provided that
\begin{align*}
\Gamma\text{-}\lim_{j\to\fz}\mathscr{E}_{p,\,q}^{\varphi_j}
= \mathscr{D}^0_{p,\,q}
\quad\text{on}\,\ W^{1,\,p}(\rn).
\end{align*}
Alternatively, we simply say that the \emph{MS formula holds in the sense of $\Gamma$-convergence}.

\item If $m\in\nn$, $p\in(1,\infty)$, $X=W^{m,\,p}(\rn)$ and $ \mathscr{D}= \mathscr{D}^{m,\,\mu}_{p,\,q}$,
then we say that \emph{the sequence $\{\mathscr{E}_{p,\,q}^{m,\,\varphi_j}\}_{j\in\nn}$ $\Gamma$-converges
to $\mathscr{D}_{p,\,q}^{m,\,\mu}$ on $W^{m,\,p}(\rn)$} provided that
\begin{align*}
\Gamma\text{-}\lim_{j\to\fz} \mathscr{E}_{p,\,q}^{m,\,\varphi_j}
= \mathscr{D}^{m,\,\mu}_{p,\,q}
\quad\text{on}\,\ W^{m,\,p}(\rn).
\end{align*}
When $p=1$, we replace $W^{m,\,p}(\rn)$ and $\mathscr{D}_{p,\,q}^{m,\,\mu}$ by ${\rm BV}^m(\rn)$ and
$\mathscr{D}_q^{m,\,\mu}$, respectively. In either case, for simplicity, we say that the \emph{BBM
formula holds in the sense of $\Gamma$-convergence}.

\end{enumerate}
The same terminology applies verbatim to strong $\Gamma$-convergence, with $\Gamma\text{-}\lim$ replaced by $s\Gamma\text{-}\lim$.
\end{remark}


Several recent works have investigated the variational convergence of nonlocal energies in the BBM setting.
For $m=1$, $p=q\in[1,\fz)$ and nonradial kernel sequences $\{\varphi_j\}_{j\in\nn}$, Gennaioli and Stefani
\cite{Gennaioli-Stefani-ARMA-2026} considered weights $\rho_j(h)=|h|^p\,\varphi_j(h)$. In our notation, their
unit-mass conditions take the form
\begin{align}\label{eq-ARMA1}
\sup_{r\in(0,\fz)} \limsup_{j\to\fz} \int_\rn \varphi_{j}(h) \min\lf\{|h|^p,\,r^p \r\} dh<\fz
\end{align}
and, for any $\psi\in C_{\rm c}(\rn)$,
\begin{align}\label{eq-ARMA2}
\lim_{j\to\fz} \int_\rn \psi(h)\,|h|^p\,\varphi_{j}(h)\,dh
= \psi(0).
\end{align}
Indeed, the combination of \eqref{eq-ARMA1} and \eqref{eq-ARMA2} is equivalent to the $\BBMmq$ condition
(with $p=q$ and $m=1$); see Proposition \ref{Prop-ARMA=BBM-MS} below. Gennaioli
and Stefani \cite[Theorem 1.1]{Gennaioli-Stefani-ARMA-2026} allow the vague limit $\alpha\delta_0$ with
$\alpha\ge0$ and a finite positive measure of mass $\alpha$. In the normalized case $\alpha=1$, conditions
\eqref{eq-ARMA1} and \eqref{eq-ARMA2} along some subsequence are necessary and sufficient for the existence
of a Borel probability measure $\mu$ on $\msn$ and a subsequence $\{\varphi_{j_k}\}_{k\in\nn}$ along with
the strong $\Gamma$-convergence
\begin{align*}
s\Gamma\text{-}\lim_{k\to\fz} \mathscr{E}_{p,\,p}^{1,\,\varphi_{j_k}}
= \mathscr{D}_{p,\,p}^{1,\,\mu}
\quad\text{on } W^{1,\,p}(\rn)
\end{align*}
holds. When $p=1$, the space $W^{1,\,p}(\rn)$ and the limiting energy $\mathscr{D}_{p,\,p}^{1,\,\mu}$ are replaced by
$\mathrm{BV}(\rn)$ and $\mathscr{D}_{1}^{1,\,\mu}$, respectively. In particular, the $\Gamma$-limsup inequality in
\cite{Gennaioli-Stefani-ARMA-2026} can be realized by the constant recovery sequence $f_k\equiv f$, which explains
why their $\Gamma$-convergence result is strong in the terminology adopted in this paper. Moreover, the argument of \cite{Gennaioli-Stefani-ARMA-2026}
gives that \eqref{eq-ARMA1}--\eqref{eq-ARMA2} is equivalent to the subsequential pointwise convergence
\begin{align}\label{eq-ARMA-pw}
\lim_{k\to\fz} \mathscr{E}_{p,\,p}^{1,\,\varphi_{j_k}}(f)
= \mathscr{D}_{p,\,p}^{1,\,\mu}(f)
\quad\text{for every } f\in W^{1,\,p}(\rn).
\end{align}
In the case $p=q=2$, the pointwise convergence \eqref{eq-ARMA-pw} had previously been established by
Davoli, Di Fratta and Pagliari \cite[Theorems 1.1 and 1.2]{Davoli-DiFratta-Pagliari-ProcRSEC-2026} for nonradial
kernel sequences, using Fourier transform techniques. In \cite[Section 5.3]{Davoli-DiFratta-Pagliari-ProcRSEC-2026},
they further observed that their sufficient conditions, which correspond to the case $p=q=2$ of \eqref{eq-ARMA1}--\eqref{eq-ARMA2},
yield standard $\Gamma$-convergence; however, the \emph{necessity} of these conditions for standard
$\Gamma$-convergence was left unresolved.

For further results concerning $\Gamma$-convergence in the BBM setting, we refer the reader to \cite{Brezis-Nguyen-AnnPDE-2018,
Crismale-DeLuca-Kubin-Ninno-Ponsiglione-JFA-2023,Ponce-CVPDE-2004}. A common feature of these works is that they address only
the case $m=1$. These developments lead naturally to the following question:
\begin{center}
\begin{minipage}{13cm}
\emph{For general exponents $p$ and $q$ and nonradial kernels $\{\varphi_j\}_{j\in\nn}$, do the same
kernel conditions $\MSq$ and $\BBMmq$ also characterize the $\Gamma$-convergence (in either the standard or strong
sense) of the corresponding energy sequences?}
\end{minipage}
\end{center}
Building on and extending \cite{Davoli-DiFratta-Pagliari-ProcRSEC-2026,Gennaioli-Stefani-ARMA-2026}, we answer this question
in the affirmative; see Theorems \ref{thm-MS} and \ref{thm-BBM} below. For the BBM formula, the characterization is established
after passing to a suitable subsequence of $\{\varphi_j\}_{j\in\nn}$, if necessary.

\subsection{Main results}\label{ss1.3}

In this paper, we establish the equivalence between the $\MSq$ condition and the validity of the MS formula
in each of the following senses: pointwise convergence, $\Gamma$-convergence and strong $\Gamma$-convergence.

\begin{theorem}\label{thm-MS}
Let $p\in[1,\fz)$ and $q\in(0,\fz)$. For any $j\in\nn$, assume that $0\le \varphi_j\in L_{\rm loc}^1(\rn \setminus\{0\})$,
with $\varphi_j(0)$ possibly equal to $+\infty$. Then, the following statements are equivalent.
\begin{enumerate}[{\rm (i)}]
\item The sequence $\{\varphi_j\}_{j\in\nn}$ satisfies the $\MSq$ condition.

\item The MS formula holds in the sense of pointwise convergence, that is,
\begin{align*}
\lim_{j\to\fz}\mathscr{E}_{p,\,q}^{\varphi_j}(f)
= \mathscr{D}^0_{p,\,q}(f)
\quad\text{for every}\,\ f\in W^{1,\,p}(\rn).
\end{align*}

\item The MS formula holds in the sense of $\Gamma$-convergence, that is,
\begin{align*}
\Gamma\text{-}\lim_{j\to\fz} \mathscr{E}_{p,\,q}^{\varphi_j}
= \mathscr{D}^0_{p,\,q}
\quad\text{on}\,\ W^{1,\,p}(\rn).
\end{align*}

\item The MS formula holds in the sense of strong $\Gamma$-convergence, that is,
\begin{align*}
s\Gamma\text{-}\lim_{j\to\fz} \mathscr{E}_{p,\,q}^{\varphi_j}
= \mathscr{D}^0_{p,\,q}
\quad\text{on}\,\ W^{1,\,p}(\rn).
\end{align*}
\end{enumerate}
\end{theorem}

Next, we establish that, after passing to a suitable subsequence of $\{\varphi_j\}_{j\in\nn}$ if necessary, the
$\BBMmq$ condition is equivalent to the validity of the BBM formula in each of the following senses: pointwise
convergence, $\Gamma$-convergence and strong $\Gamma$-convergence.

\begin{theorem}\label{thm-BBM}
Let $p\in[1,\fz)$, $q\in(0,\fz)$ and $m\in\nn$. For any $j\in\nn$, assume that $0\le \varphi_j\in L_{\rm loc}^1(\rn\setminus\{0\})$,
with $\varphi_j(0)$ possibly equal to $+\infty$. Then, the following statements are equivalent.
\begin{enumerate}[{\rm (i)}]
\item There exists a subsequence $\{\varphi_{j_k}\}_{k\in\nn}$ satisfying the $\BBMmq$ condition.

\item The BBM formula holds in the sense of pointwise convergence, that is, there exist a Borel probability
measure $\mu$ on $\msn$ and a subsequence $\{\mathscr{E}_{p,\,q}^{m,\,\varphi_{j_k}}\}_{k\in\nn}$ such that
\begin{align*}
\lim_{k\to\fz} \mathscr{E}_{p,\,q}^{m,\,\varphi_{j_k}}(f)
= \mathscr{D}_{p,\,q}^{m,\,\mu}(f)
\quad\text{for every}\,\ f\in W^{m,\,p}(\rn).
\end{align*}

\item The BBM formula holds in the sense of $\Gamma$-convergence, that is, there exist a Borel probability
measure $\mu$ on $\msn$ and a subsequence $\{\mathscr{E}_{p,\,q}^{m,\,\varphi_{j_k}}\}_{k\in\nn}$ such that
\begin{align*}
\Gamma\text{-}\lim_{k\to\fz} \mathscr{E}_{p,\,q}^{m,\,\varphi_{j_k}}
= \mathscr{D}_{p,\,q}^{m,\,\mu}
\quad\text{on}\,\ W^{m,\,p}(\rn).
\end{align*}

\item The BBM formula holds in the sense of strong $\Gamma$-convergence, that is, there exist a Borel probability
measure $\mu$ on $\msn$ and a subsequence $\{\mathscr{E}_{p,\,q}^{m,\,\varphi_{j_k}}\}_{k\in\nn}$ such that
\begin{align*}
s\Gamma\text{-}\lim_{k\to\fz} \mathscr{E}_{p,\,q}^{m,\,\varphi_{j_k}}
= \mathscr{D}_{p,\,q}^{m,\,\mu}
\quad\text{on}\,\ W^{m,\,p}(\rn).
\end{align*}
\end{enumerate}
In statements {\rm (ii)}, {\rm (iii)} and {\rm (iv)}, when $p=1$, the space $W^{m,\,p}(\rn)$
and the local energy $\mathscr{D}_{p,\,q}^{m,\,\mu}$ can also be replaced by $\BVm(\rn)$ and $\mathscr{D}_{q}^{m,\,\mu}$, respectively.
\end{theorem}

\begin{remark}\label{rem-MS-BBM}
In both Theorems \ref{thm-MS} and \ref{thm-BBM}, the implications (iv) $\Rightarrow$ (iii) are immediate by
Remark \ref{rem-strong-Gamma}; moreover, the implications (iii) $\Rightarrow$ (ii) are established in Section
\ref{sec-Gamma-to-point}, based on the admissibility discussion of the sequence $\{\varphi_j\}_{j\in\nn}$.
The remaining implications (ii) $\Rightarrow$ (i) $\Rightarrow$ (iv) of Theorem \ref{thm-MS} are established
in Section \ref{sec-MS}, while those of Theorem \ref{thm-BBM} are established in Section \ref{sec-BBM}.

Taken together, Theorems \ref{thm-MS} and \ref{thm-BBM}
establish a unified characterization of the MS and BBM formulae in a substantially broader framework, showing that
the corresponding kernel conditions simultaneously characterize pointwise convergence, the standard $\Gamma$-convergence
and the strong $\Gamma$-convergence. In view of Remark \ref{rem-strong-Gamma} and Theorems \ref{thm-MS} and \ref{thm-BBM},
we further obtain
\begin{align*}
\text{strong $\Gamma$-convergence}
\quad\Leftrightarrow\quad
\text{standard $\Gamma$-convergence}
\quad\Leftrightarrow\quad
\text{pointwise convergence}
\end{align*}
for the MS formula along the full sequence and for the BBM formula with the subsequential quantifiers in Theorem \ref{thm-BBM}.
These results serve as a definitive reference for subsequent research on the asymptotic analysis of nonlocal energies across
various areas, including partial differential equations, the calculus of variations, and physical modeling.
\end{remark}

\begin{remark}\label{rem-BBM}
According to Remark \ref{rem-MS-BBM}, Theorem \ref{thm-BBM} is proved by checking the implications
\begin{align*}
\rm (iv) \ \Rightarrow\ (iii) \ \Rightarrow\ (ii) \ \Rightarrow\ (i)\
\Rightarrow\ (iv).
\end{align*}
Passing to a subsequence in Theorem \ref{thm-BBM} is necessary only for the implication (i) $\Rightarrow$ (iv);
all the other implications hold along the same subsequence. Moreover, this passage to a subsequence is, in general,
unavoidable. Indeed, in Example \ref{eg-BBM-counterexam} below, we construct a sequence $\{\varphi_j\}_{j\in\nn}$
satisfying $\BBMmq$ along the full sequence for which the full energy sequence
 has no pointwise limit at
a particular $f\in W^{m,\,p}(\rn)$.

\end{remark}

For radial kernels, the full-sequence versions of (i)--(iv) in Theorem \ref{thm-BBM} are equivalent and, moreover, the angular measure $d\mu$ reduces to $(\omega_{n-1})^{-1}\,d\sigma$. More precisely, we have the following theorem.

\begin{theorem}\label{thm-BBM-radial}
Let $p\in[1,\fz)$, $q\in(0,\fz)$ and $m\in\nn$. For any $j\in\nn$, assume that $0\le\varphi_j\in L_{\rm loc}^1(\rn \setminus\{0\})$
is radial, with $\varphi_j(0)$ possibly equal to $+\infty$. Then the following statements are equivalent.
\begin{enumerate}[{\rm (i)}]
\item The sequence $\{\varphi_j\}_{j\in\nn}$ satisfies the $\BBMmq$ condition.

\item The BBM formula holds in the sense of pointwise convergence, that is,
\begin{align*}
\lim_{j\to\fz} \mathscr{E}_{p,\,q}^{m,\,\varphi_j}(f)
= \mathscr{D}_{p,\,q}^{m}(f)
\quad\text{for every}\,\ f\in W^{m,\,p}(\rn).
\end{align*}

\item The BBM formula holds in the sense of $\Gamma$-convergence, that is,
\begin{align*}
\Gamma\text{-}\lim_{j\to\fz} \mathscr{E}_{p,\,q}^{m,\,\varphi_j}
= \mathscr{D}_{p,\,q}^m
\quad\text{on}\,\ W^{m,\,p}(\rn).
\end{align*}

\item The BBM formula holds in the sense of strong $\Gamma$-convergence, that is,
\begin{align*}
s\Gamma\text{-}\lim_{j\to\fz} \mathscr{E}_{p,\,q}^{m,\,\varphi_j}
= \mathscr{D}_{p,\,q}^m
\quad\text{on}\,\ W^{m,\,p}(\rn).
\end{align*}
\end{enumerate}
In statements {\rm (ii)}, {\rm (iii)} and {\rm (iv)}, when $p=1$, the space $W^{m,\,p}(\rn)$
and the local energy $\mathscr{D}_{p,\,q}^m$ can also be replaced by $\BVm$ and $\mathscr{D}_q^m$, respectively.
\end{theorem}

\subsection{Proof strategies and organization of the paper}\label{ss1.4}

We briefly outline the main ideas behind the proofs of Theorems \ref{thm-MS} and \ref{thm-BBM}. In both cases,
by Remark \ref{rem-MS-BBM}, it suffices to establish
\begin{align*}
\rm (iii)\Rightarrow(ii)\Rightarrow(i)\Rightarrow(iv).
\end{align*}
We first explain the common argument for $\rm (iii)\Rightarrow(ii)$ and then discuss the remaining implications in the MS and BBM settings.

\medskip

\textit{From standard $\Gamma$-convergence to pointwise convergence.}
The main difficulty is that standard $\Gamma$-convergence provides a recovery sequence $f_j\to f$, rather than
convergence of the energies evaluated at the fixed function $f$. To overcome this, we establish in Proposition
\ref{prop-BBM-Gamma-to-cond} a stability estimate showing that, whenever $f_j\to f$ in $L^p(\rn)$ with
$f\in W^{m,\,p}(\rn)$ (or $f\in\BVm(\rn)$ when $p=1$) and the energies of $f_j$ remain bounded,
\begin{align}\label{eq-intro-Gamma-stability}
\limsup_{j\to\infty} \mathscr{E}_{p,\,q}^{m,\,\varphi_j}(f)
\le \limsup_{j\to\infty} \mathscr{E}_{p,\,q}^{m,\,\varphi_j}(f_j).
\end{align}
The proof mainly relies on uniform admissibility of the kernel sequence, established in Section \ref{sec-Perliminary3}. Applying \eqref{eq-intro-Gamma-stability} to a recovery
sequence and combining it with the $\Gamma$-liminf inequality for the constant sequence yields pointwise convergence.

\medskip

\textit{Proof strategies for the MS formula.}
To prove $\rm (ii)\Rightarrow(i)$ in Theorem \ref{thm-MS}, we use the identity
\begin{align}\label{eq-limit-Delta}
\lim_{|h|\to\infty} \|\Delta_hf\|_{L^p(\rn)}
=2^\frac{1}{p}\,\|f\|_{L^p(\rn)}
\end{align}
for every $f\in L^p(\rn)$; see Lemma \ref{Lem-LXlem2829}(v). By choosing suitable compactly supported test
functions and applying a scaling argument, we show that, for every $\delta\in(0,\infty)$,
\begin{align}\label{eq-limit-delta-MS}
\lim_{j\to\infty}
\int_{|h|\ge\delta}\varphi_j(h)\,dh=1;
\end{align}
see Lemmas \ref{Lem-MS-upper} and \ref{Lem-MS-lower}. In particular, taking $\delta=1$ gives the normalization
condition $\NMSq$. The vanishing condition $\VMSq$ then follows by comparing the contribution from large
increments with the limiting energy and using a lower bound for first-order differences near the origin.

For $\rm (i)\Rightarrow(iv)$, the key observation is that $\MSq$ forces the kernel mass to concentrate at infinity,
as expressed by \eqref{eq-limit-delta-MS}. The pointwise convergence follows from \eqref{eq-limit-Delta}, while the
strong $\Gamma$-liminf inequality is obtained using
\begin{align*}
\|\Delta_hf_j-\Delta_hf\|_{L^p(\rn)}
\le 2\,\|f_j-f\|_{L^p(\rn)}
\end{align*}
together with the normalization of the kernel mass at infinity.

\medskip

\textit{Proof strategies for the BBM formula.}
To prove $\rm (ii)\Rightarrow(i)$ in Theorem \ref{thm-BBM}, we use radial test functions and compare higher-order
differences with directional derivatives. The radiality assumption is particularly useful because, for such functions,
the quantity $\|(\xi\cdot\nabla)^mf\|_{L^p(\rn)}$ is independent of $\xi\in\msn$. This allows us to separate the
weighted kernel mass from the directional derivative norm. Combining a higher-order difference approximation with suitable
dilations, we obtain the upper and lower estimates for the weighted mass near the origin; see Lemmas \ref{Lem-DeltaNabla},
\ref{Lem-upper}, and \ref{Lem-lower}. These estimates yield the normalization condition $\NBBMmq$. The vanishing condition
$\VBBMmq$ is then established using compactly supported test functions to control the contribution from increments bounded
away from the origin.

The proof of $\rm (i)\Rightarrow(iv)$ presents an additional difficulty in identifying the limiting measure $\mu$,
since the angular distributions of nonradial kernels may vary with $j$. To overcome this, we associate with each kernel
$\varphi_j$ a finite Borel measure $\mu_j$ on $\msn$, defined by
\begin{align}\label{eq-intro-angular-measure}
d\mu_j(\xi):= \left(\int_0^1 t^{n+mq}\,\varphi_j(t\xi)\,\frac{dt}{t} \right) d\sigma(\xi).
\end{align}
By $\NBBMmq$ and the compactness of $\msn$, we can extract a subsequence $\{\mu_{j_k}\}_{k\in\nn}$ converging weakly
to a Borel probability measure $\mu$ on $\msn$. This measure determines the limiting directional energy; see Lemma \ref{Lem-lim-varphi=mu}.

The strong $\Gamma$-limsup inequality follows from the comparison between higher-order differences and directional derivatives,
together with the vanishing condition $\VBBMmq$. For the strong $\Gamma$-liminf inequality, we mollify a sequence $f_k\to f$
in $L^p(\rn)$ and prove in Lemma \ref{Lem2-Gamma-converge} that
\begin{align*}
\lim_{\delta\to0_+} \lim_{k\to\fz} \int_{|h|<\delta} \|\Delta_h^mf_{k,\uc}\|_{L^p(\rn)}^q\,\varphi_{j_k}(h)\,dh
= \mathscr{D}_{p,\,q}^{m,\,\mu}(f_\uc).
\end{align*}
Finally, letting $\uc\to0_+$ and using Lemma \ref{Lem-uc-to-0} yields the desired liminf inequality.
This approach is completely different from that of
\cite{Gennaioli-Stefani-ARMA-2026}, since the latter treats only first-order differences
and the case $p=q$.

When the kernels are radial, the measures $\mu_j$ in \eqref{eq-intro-angular-measure} converge to the normalized surface measure
$(\omega_{n-1})^{-1}\sigma$ without passing to a subsequence. Consequently, the BBM convergence holds along the full sequence,
as stated in Theorem \ref{thm-BBM-radial}. The counterexample in Section \ref{ss-eg} shows that this conclusion generally fails
for nonradial kernel sequences.

\medskip

\textit{Organization of the paper.}
Section \ref{sec-Perliminary} contains the preliminary results. In Section \ref{sec-Perliminary1}, we establish
the connections between the directional energies $\mathscr{D}_{p,\,q}^{m,\,\mu}$ and $\mathscr{D}_q^{m,\,\mu}$
and the classical Sobolev and bounded variation spaces $\mathring{W}^{m,\,p}(\rn)$ and $\mathringBVm(\rn)$,
respectively. Section \ref{sec-Perliminary2} collects basic properties of higher-order differences used throughout
the paper. In Section \ref{ss-BBM-MS-cond}, we characterize the conditions $\MSq$ and $\BBMmq$ in terms of measure
convergence and uniform kernel bounds. In particular, we establish the equivalence between $\BBMmq$ and the
conditions \eqref{eq-ARMA1}--\eqref{eq-ARMA2} introduced in \cite{Gennaioli-Stefani-ARMA-2026}.

 Section~\ref{sec-Gamma-to-point} is devoted to showing that standard
$\Gamma$-convergence implies pointwise convergence for the nonlocal energies
considered here. In particular, Section~\ref{sec-Perliminary3} is devoted to uniform admissibility: we show
that each of the assertions in Theorems~\ref{thm-MS} and~\ref{thm-BBM}
independently guarantees the uniform admissibility of the full MS sequence
or of the indicated BBM subsequence. These results are used in
Sections~\ref{sec3.2} and~\ref{ss-Gamma=>pw}, where we prove the implication
(iii) $\Rightarrow$ (ii) in both Theorems~\ref{thm-MS} and~\ref{thm-BBM}.

Section \ref{sec-MS} completes the proof of Theorem \ref{thm-MS}. The implications $\rm (ii)\Rightarrow(i)$ and
$\rm (i)\Rightarrow(iv)$ are established in Sections \ref{sec-MS1} and \ref{sec-MS2}, respectively.

Section \ref{sec-BBM} is devoted to the proofs of Theorems \ref{thm-BBM} and \ref{thm-BBM-radial}. The implications
$\rm (ii)\Rightarrow(i)$ and $\rm (i)\Rightarrow(iv)$ in Theorem \ref{thm-BBM} are established in Sections
\ref{sec-BBM1} and \ref{sec-BBM2}, respectively. We then prove Theorem \ref{thm-BBM-radial} in Section \ref{sec-BBM4}
and construct a counterexample in Section \ref{ss-eg}, showing that subsequence extraction is, in general, necessary for nonradial kernels.

Finally, in Section \ref{sec-example}, we apply Theorems \ref{thm-MS} and \ref{thm-BBM-radial} to establish
pointwise convergence, $\Gamma$-convergence, and strong $\Gamma$-convergence of classical and logarithmic
Besov energies.

\bigskip

\noindent{\bf Notation.\,}

\begin{itemize}
\item Let $\nn:=\{1,2,\dots\}$ and $\zz_+:=\{0\}\cup\nn$.

\item Denote by $C^\fz(\rn)$ the space of infinitely differentiable functions on $\rn$,
and by $C_{\rm c}^\fz(\rn)$ the space of functions in $C^\fz(\rn)$ with compact support.


\item For any Lebesgue measurable set $E\subset\rn$, denote by $|E|$ and $\mathbf{1}_E$, respectively, the {\it Lebesgue measure} of $E$
and the {\it characteristic function} of $E$.

\item Let $L_\loc^1(\rn)$ be the space of locally integrable functions on $\rn$. For any
$E\subset\rn$ with $|E|\in (0,\fz)$, the {\it integral average of $f$
over $E$} is defined by
\begin{align*}
\fint_E f(x)\,dx:=\frac{1}{|E|} \int_E f(x)\,dx.
\end{align*}

\item Let $d\sigma$ be the surface measure on the unit sphere $\msn$ in $\rn$ and $\omega_{n-1}:=\sigma(\msn)$.
For any $\sigma$-integrable function $g$ (or nonnegative measurable $g$) on $\msn$, set
\begin{align*}
\fint_{\msn} g(\xi)\,d\sigma(\xi):=\frac1{\omega_{n-1}}\int_{\msn} g(\xi)\,d\sigma(\xi).
\end{align*}

\item The notation $U\lesssim V$ or $V\gtrsim U$ means that $U\le C V$ for some positive constant $C$
independent of the variables under consideration. The symbol $U\approx V$ means that $U\lesssim V\lesssim U$.
\end{itemize}

\section{Preliminaries}\label{sec-Perliminary}

\subsection{Higher-order Sobolev and {\rm BV} spaces}\label{sec-Perliminary1}

Let $\mathcal{M}(\rn)$ denote the space of finite signed Radon measures on $\rn$. For any $\mu\in\mathcal{M}(\rn)$, set
\begin{align*}
\|\mu\|_{\mathcal{M}(\rn)}:=|\mu|(\rn)
:= \sup\lf\{\int_\rn \phi(x)\,d\mu(x):\,\phi\in C_{\rm c}(\rn),\,\ |\phi|\le1 \r\},
\end{align*}
where $C_{\rm c}(\rn)$ denotes the space of continuous functions on $\rn$ with compact support.

\begin{definition}\label{def-Sobolev-BV}
Let $m\in\nn$, $p\in[1,\fz)$ and $q\in(0,\fz)$.
\begin{enumerate}[{\rm (i)}]
\item The \emph{homogeneous Sobolev space $\mathring{W}^{m,\,p}(\rn)$ of order $m$}
consists of all $f\in L^1_{\rm loc}(\rn)$ such that
\begin{align}\label{eq-Sobolev-norm}
\|f\|_{\mathring{W}^{m,\,p}(\rn)}
:= \sum_{\{\alpha\in\zz_+^n:\,|\alpha|=m\}} \frac{m!}{\alpha!}\,\|\partial^\alpha f\|_{L^p(\rn)}<\fz.
\end{align}
The corresponding \emph{inhomogeneous Sobolev space of order $m$} is defined by
\begin{align*}
W^{m,\,p}(\rn)
:= L^p(\rn)\cap\mathring{W}^{m,\,p}(\rn).
\end{align*}

\item The \emph{homogeneous bounded variation space $\mathringBVm(\rn)$ of order $m$}
consists of all $f\in L^1_{\rm loc}(\rn)$ such that
\begin{align}\label{eq-BV-norm}
\|f\|_{\mathringBVm(\rn)}
:= \sum_{\{\alpha\in\zz_+^n:\,|\alpha|=m\}} \frac{m!}{\alpha!}\,\|\partial^\alpha f\|_{\mathcal{M}(\rn)}<\infty.
\end{align}
The corresponding \emph{inhomogeneous bounded variation space of order $m$}
is defined by
\begin{align*}
\BVm(\rn):=L^1(\rn)\cap\mathringBVm(\rn).
\end{align*}
When $m=1$, we simply write $\mathringBVm(\rn)$ and $\BVm(\rn)$ as $\mathring{\mathrm{BV}}(\rn)$ and $\mathrm{BV}(\rn)$, respectively.
\end{enumerate}
\end{definition}

\begin{remark}\label{rem-Wmp-BV}
In view of the expression
\begin{align}\label{eq-def-xi-cdot-nabla}
(h\cdot\nabla)^mf
= \sum_{\{\alpha\in\zz_+^n:\,|\alpha|=m\}} \frac{m!}{\alpha!}\,h^\alpha\,\partial^\alpha f,
\end{align}
it is convenient for us to using the seminorms in \eqref{eq-Sobolev-norm} and \eqref{eq-BV-norm}.
Indeed, by \eqref{eq-def-xi-cdot-nabla}, we immediately obtain
\begin{align}\label{eq-hWmp}
\|(h\cdot\nabla)^mf\|_{L^p(\rn)}
\le |h|^m\,\|f\|_{\mathring{W}^{m,\,p}(\rn)}
\end{align}
and
\begin{align}\label{eq-hBVm}
\|(h\cdot\nabla)^mf\|_{{\mathcal{M}(\rn)}}
\le |h|^m\,\|f\|_{\mathringBVm(\rn)}.
\end{align}
Consequently, if $\mu$ is a probability measure on the sphere $\msn$, then
\begin{align}\label{eq-DmuWBV}
\lf(\mathscr{D}_{p,\,q}^{m,\,\mu}(f) \r)^\frac{1}{q}
\le \|f\|_{\mathring{W}^{m,\,p}(\rn)}
\qquad\text{and}\qquad
\lf(\mathscr{D}_q^{m,\,\mu}(f) \r)^\frac{1}{q}
\le \|f\|_{\mathringBVm(\rn)}.
\end{align}
\end{remark}

The following proposition shows that, up to suitable constants, the converses of the inequalities in
\eqref{eq-DmuWBV} are valid provided that $d\mu=(\omega_{n-1})^{-1}\,d\sigma$.

\begin{proposition}\label{Prop-equivalent-Sobolevnorm}
Let $p\in[1,\fz),\,q\in(0,\infty)$ and $m\in\nn$. Then, there exists a positive constant $C_{m,n}$,
depending only on $m$ and $n$, such that for any $f\in W^{m,\,p}(\rn)$,
\begin{align}\label{eq-Sobolevnorm}
(C_{m,\,n})^{-\max\{1,\,\frac{1}{q}\}}\,\|f\|_{\mathring{W}^{m,\,p}(\rn)}
\le \lf(\mathscr{D}_{p,\,q}^m(f) \r)^\frac{1}{q}
\le \|f\|_{\mathring{W}^{m,\,p}(\rn)}.
\end{align}
At the endpoint $p=1$, instead of \eqref{eq-Sobolevnorm}, it also holds that for all $f\in \mathrm{BV}^{m}(\rn)$,
\begin{align}\label{eq-BVnorm}
(C_{m,\,n})^{-\max\{1,\,\frac{1}{q}\}}\,\|f\|_{\mathringBVm(\rn)}
\le \lf(\mathscr{D}_q^m(f) \r)^\frac{1}{q}
\le \|f\|_{\mathringBVm(\rn)}.
\end{align}
\end{proposition}

The proof of Proposition \ref{Prop-equivalent-Sobolevnorm} relies on the following lemma.

\begin{lemma}\label{Lem-existence}
For $m\in\nn$ and $\mathcal{I}_m:= \{\gamma\in\zz_+^n:\,|\gamma|=m\}$, there exists a family $\{F_\gamma:\,\gamma\in \mathcal{I}_m\},$ with each
$F_\gamma$ being a polynomial, such that for all $\alpha,\gamma\in \mathcal{I}_m$,
\begin{align}\label{eq-duality}
\frac1{\omega_{n-1}}\int_\msn \xi^\alpha\,F_\gamma(\xi)\,d\sigma(\xi)=
\begin{cases}
1&\quad\text{if}\,\ \gamma=\alpha; \\
0&\quad\text{if}\,\ \gamma\neq\alpha.
\end{cases}
\end{align}
\end{lemma}

\begin{proof}
Denote by $\dot{\mathscr{P}}_m(\rn)$ the collection of all homogeneous polynomials on $\rn$ of degree $m$. Let
\begin{align*}
\dot{\mathscr{P}}_m(\msn)
:= \lf\{P|_\msn:\,P\in \dot{\mathscr{P}}_m(\rn) \r\}.
\end{align*}
Then, the family of monomials
\begin{align*}
\mathscr{B}:=\lf\{Q_\beta(\xi):=\xi^\beta:\, \beta\in\mathcal{I}_m\r\}
\end{align*}
forms a basis of $\dot{\mathscr{P}}_m(\msn)$. Applying the Gram--Schmidt procedure to the monomial basis $\mathscr{B}$
produces an orthonormal basis of $\dot{\mathscr{P}}_m(\msn)$, denoted by $\{P_\beta\}_{\beta\in\mathcal{I}_m}$, such that
\begin{align*}
\frac1{\omega_{n-1}}\int_\msn P_\beta(\xi)\,P_{\beta'}(\xi)\,d\sigma(\xi)=
\begin{cases}
 1&\quad\text{if}\,\ \beta=\beta'; \\
0&\quad\text{if}\,\ \beta\neq\beta'.
\end{cases}
\end{align*}

Define the matrix $M:= (M_{\eta,\,\beta})_{\eta,\,\beta\in\mathcal{I}_m}$, with each element $M_{\eta,\,\beta}$ defined by
\begin{align*}
M_{\eta,\,\beta}:= \frac{1}{\omega_{n-1}}\int_\msn \xi^\eta\,P_\beta(\xi)\,d\sigma(\xi).
\end{align*}
Since $\{P_\beta\}_{\beta\in\mathcal{I}_m}$ is an orthonormal basis, for each
$\eta\in\mathcal{I}_m$ the polynomial $Q_\eta(\xi)=\xi^\eta$ has the expansion
\begin{align*}
Q_\eta(\xi)
= \sum_{\beta\in \mathcal{I}_m} M_{\eta,\,\beta}\,P_\beta(\xi)
\quad\, \text{for all}\ \, \xi\in\msn.
\end{align*}
Thus $M$ is the transition matrix from the orthonormal basis
$\{P_\beta\}_{\beta\in\mathcal{I}_m}$ to the monomial basis
$\{Q_\eta\}_{\eta\in\mathcal{I}_m}$, which implies that $M$ is invertible.
Let
$N:= (N_{\eta,\,\beta})_{\eta,\,\beta\in \mathcal{I}_m}$
be the inverse matrix of $M$. Then
\begin{align}\label{eq1-orthogonal}
\sum_{\beta\in\mathcal{I}_m} M_{\alpha,\,\beta}\,N_{\beta,\,\gamma}=
\begin{cases}
 1&\quad\text{if}\,\ \alpha=\gamma; \\
 0&\quad\text{if}\,\ \alpha\neq\gamma.
\end{cases}
\end{align}

Now, for each $\gamma\in \mathcal{I}_m$, define
\begin{align*}
F_\gamma(\xi):= \sum_{\beta\in\mathcal{I}_m} N_{\beta,\,\gamma}\,P_\beta(\xi)
\quad\, \text{for all}\ \, \xi\in\msn.
\end{align*}
Since each $P_\beta$ is a polynomial, so is every $F_\gamma$. Moreover,
\begin{align*}
\frac1{\omega_{n-1}}\int_\msn \xi^\alpha\,F_\gamma(\xi)\,d\sigma(\xi)
&= \sum_{\beta\in\mathcal{I}_m} N_{\beta,\,\gamma}
\lf(\frac{1}{\omega_{n-1}} \int_\msn \xi^\alpha\,P_\beta(\xi)\,d\sigma(\xi)\r)\\
&
= \sum_{\beta\in\mathcal{I}_m} M_{\alpha,\,\beta}\,N_{\beta,\,\gamma}.
\end{align*}
This, together with \eqref{eq1-orthogonal}, yields \eqref{eq-duality}.
\end{proof}

\begin{proof}[Proof of Proposition \ref{Prop-equivalent-Sobolevnorm}]
By \eqref{eq-DmuWBV}, it suffices to prove the first inequalities in \eqref{eq-Sobolevnorm}
and
\eqref{eq-BVnorm}. To this end, we let $\{F_\gamma:\,\gamma\in\zz_+^n,\, |\gamma|=m\}$ be the polynomials determined in Lemma \ref{Lem-existence}.
Since each $F_\gamma$ is continuous, $\msn$ is compact and there are only finitely many $\gamma$ with $|\gamma|=m$, it follows that
\begin{align}\label{eq-F-gamma-bounded}
C_{m,\,n}:= \sum_{|\gamma|=m}\,\sup_{\xi\in\msn} |F_\gamma(\xi)|
\end{align}
is a positive finite constant.

Now, we let $f\in W^{m,\,p}(\rn)$ and prove the first inequality in \eqref{eq-Sobolevnorm}. For almost every $x\in\rn$, since
the family $\{F_\gamma:\gamma\in\zz_+^n,\ |\gamma|=m\}$ satisfies \eqref{eq-duality} and
\begin{align*}
(\xi\cdot\nabla)^mf(x)
= \sum_{|\gamma|=m} \frac{m!}{\gamma!}\,\xi^\gamma\,\partial^\gamma f(x),
\end{align*}
we deduce that, for any $\alpha\in\zz_+^n$ with $|\alpha|=m$,
\begin{align*}
\frac{m!}{\alpha!}\,\partial^\alpha f(x)
&= \sum_{|\gamma|=m} \frac{m!}{\gamma!}\,\partial^\gamma f(x)
\lf(\fint_\msn \xi^\gamma\,F_\alpha(\xi)\,d\sigma(\xi) \r)\\
&
= \fint_\msn (\xi\cdot\nabla)^mf(x)\,F_\alpha(\xi)\,d\sigma(\xi).
\end{align*}
From $p\in[1,\infty)$ and the Minkowski inequality, it follows that
\begin{align}\label{eq-norm-Sob}
\frac{m!}{\alpha!}\,\|\partial^\alpha f\|_{L^p(\rn)}
&\le \lf(\sup_{\xi\in\msn} |F_\alpha(\xi)|\r)\,\lf(\fint_\msn \|(\xi\cdot\nabla)^mf\|_{L^p(\rn)}\,d\sigma(\xi) \r).
\end{align}
If $q\in[1,\infty)$, then applying the H\"older inequality to the integral average in \eqref{eq-norm-Sob} yields
\begin{align*}
\frac{m!}{\alpha!}\,\|\partial^\alpha f\|_{L^p(\rn)}
&\le \lf(\sup_{\xi\in\msn} |F_\alpha(\xi)|\r)\,\lf(\fint_\msn \|(\xi\cdot\nabla)^mf\|_{L^p(\rn)}^q\, d\sigma(\xi) \r)^\frac{1}{q}\\
&
= \lf(\sup_{\xi\in\msn} |F_\alpha(\xi)|\r)\, \lf(\mathscr{D}_{p,\,q}^m(f) \r)^\frac{1}{q}.
\end{align*}
Summing over all $\alpha\in\zz_+^n$ with $|\alpha|=m$ and using the definition of $C_{m,n}$ in \eqref{eq-F-gamma-bounded}, we obtain
\begin{align}\label{eq-q>1}
\|f\|_{\mathring{W}^{m,\,p}(\rn)}
&= \sum_{|\alpha|=m}\frac{m!}{\alpha!}\,\|\partial^\alpha f\|_{L^p(\rn)}\\
&\le \sum_{|\alpha|=m}\lf(\sup_{\xi\in\msn} |F_\alpha(\xi)|\r)\, \lf(\mathscr{D}_{p,\,q}^m(f) \r)^\frac{1}{q} \notag
\\
&= C_{m,\,n} \lf(\mathscr{D}_{p,\,q}^m(f) \r)^\frac{1}{q}. \notag
\end{align}
If $q\in(0,1)$, then applying \eqref{eq-hWmp} to \eqref{eq-norm-Sob} yields
\begin{align*}
\frac{m!}{\alpha!}\,\|\partial^\alpha f\|_{L^p(\rn)}
&\le \lf(\sup_{\xi\in\msn} |F_\alpha(\xi)|\r)\,\fint_\msn \|(\xi\cdot\nabla)^mf\|_{L^p(\rn)}^q \,\|f\|_{\mathring{W}^{m,\,p}(\rn)}^{1-q}\, d\sigma(\xi)\\
&
= \lf(\sup_{\xi\in\msn} |F_\alpha(\xi)|\r)\,\|f\|_{\mathring{W}^{m,\,p}(\rn)}^{1-q}\,\mathscr{D}_{p,\,q}^m(f).
\end{align*}
In this case, summing over all $\alpha\in\zz_+^n$ with $|\alpha|=m$, we have
\begin{align*}
\|f\|_{\mathring{W}^{m,\,p}(\rn)}
\le C_{m,\,n}\,\|f\|_{\mathring{W}^{m,\,p}(\rn)}^{1-q}\,\mathscr{D}_{p,\,q}^m(f),
\end{align*}
that is,
\begin{align}\label{eq-q<1}
\|f\|_{\mathring{W}^{m,\,p}(\rn)}
\le \lf(C_{m,\,n}\,\mathscr{D}_{p,\,q}^m(f) \r)^\frac{1}{q}.
\end{align}
Combining \eqref{eq-q>1} and \eqref{eq-q<1} proves the first inequality in \eqref{eq-Sobolevnorm}.

The validity of the first inequality in \eqref{eq-BVnorm} for $f\in \mathrm{BV}^{m}(\rn)$
follows similarly. In this case, the norm $\|\cdot\|_{L^p(\rn)}$ in \eqref{eq-norm-Sob}
should be replaced by $\|\cdot\|_{\cm(\rn)}$. Then, arguing as in the proof of \eqref{eq-q>1} for
$q\in[1,\infty)$, and as in the proof of \eqref{eq-q<1} for $q\in(0,1)$ (but now using \eqref{eq-hBVm}
instead of \eqref{eq-hWmp}), respectively, we obtain
\begin{align*}
\|f\|_{\mathringBVm(\rn)}
\le(C_{m,\,n})^{\max\{1,\,\frac{1}{q}\}} \lf(\mathscr{D}_q^m(f) \r)^\frac{1}{q}.
\end{align*}
This proves the first inequality in \eqref{eq-BVnorm}.
\end{proof}

\subsection{Basic lemmas on higher-order differences}\label{sec-Perliminary2}

The following basic facts are taken from \cite[Lemmas 2.8 and 2.9]{Liu-Zhang-2026-submit}, except for \eqref{eq-Delta-BV-upper}.

\begin{lemma}[\cite{Liu-Zhang-2026-submit}]\label{Lem-LXlem2829}
Let $p\in [1,\fz)$ and $m\in\nn$.
\begin{enumerate}[{\rm (i)}]
\item For any $f\in L^p(\rn)$ and $h\in\rn$,
$\|\Delta_h^mf\|_{L^p(\rn)}\le 2^m\,\|f\|_{L^p(\rn)}.$

\item For any $f\in W^{m,\,p}(\rn)$ and $h\in \rn$,
\begin{align}\label{eq-Delta-Lp-upper}
\|\Delta_{h}^mf\|_{L^p(\rn)}
\le \|(h\cdot\nabla)^mf\|_{L^p(\rn)}
\le |h|^m\,\|f\|_{\mathring{W}^{m,\,p}(\rn)}.
\end{align}
Additionally, for any $f\in \mathrm{BV}^m(\rn)$ and $h\in\rn$,
\begin{align}\label{eq-Delta-BV-upper}
\|\Delta_{h}^mf\|_{L^1(\rn)}
\le \|(h\cdot\nabla)^mf\|_{\mathcal{M}(\rn)}
\le |h|^m\,\|f\|_{\mathringBVm(\rn)}.
\end{align}

\item For any $0\nequiv f\in L^p(\rn)$, there exist $t_f,\,C_f\in (0,\fz)$, such that for
all $t\in (0,t_f)$ and $\xi\in \msn$,
\begin{align*}
\|\Delta_{t\xi}^mf\|_{L^p(\rn)}\ge C_f\,t^m.
\end{align*}

\item For any $f\in L_{\rm loc}^1(\rn)$ supported in a ball $B_R$ with $R\in(0,\fz)$, and any $h\in\rn$ with $|h|>2R$,
\begin{align*}
\|\Delta_{h}^mf\|_{L^p(\rn)}^p
= \sum_{j=0}^m \binom{m}{j}^p\,\|f\|_{L^p(\rn)}^p.
\end{align*}

\item For a general $f\in L^p(\rn)$,
\begin{align*}
\lim_{|h|\to\fz} \|\Delta_h^mf\|_{L^p(\rn)}^p
= \sum_{j=0}^m \binom{m}{j}^p\, \|f\|_{L^p(\rn)}^p.
\end{align*}
\end{enumerate}
\end{lemma}

\begin{proof}
It suffices to show \eqref{eq-Delta-BV-upper}; the others have been established in \cite[Lemmas 2.8 and 2.9]{Liu-Zhang-2026-submit}.
Fix $f\in \mathrm{BV}^m(\rn)$ and $h\in\rn$. Set
\begin{align*}
\mu_h:= (h\cdot\nabla)^mf.
\end{align*}
Then $\mu_h\in\mathcal{M}(\rn)$.
Choose $0\le\eta\in C_{\rm c}^\fz(\rn)$ such that $\int_\rn \eta(x)\,dx=1$. For any $\uc\in(0,1)$ and $x\in\rn$, define
\begin{align*}
\eta_\uc(x):= \uc^{-n}\,\eta(\uc^{-1}x)
\qquad\text{and}\qquad
 f_\uc:= f\ast\eta_\uc.
\end{align*}
Since $f_\uc\in W^{m,\,1}(\rn)$, it follows from \eqref{eq-Delta-Lp-upper} and the fact $\int_\rn \eta(x)\,dx=1$ that
\begin{align*}
\|\Delta_{h}^mf_\uc\|_{L^1(\rn)}
&\le \|(h\cdot\nabla)^mf_\uc\|_{L^1(\rn)} \\
&= \int_\rn \lf|\lf\langle(h\cdot\nabla)^mf,\,\eta_\uc(x-\cdot)\r\rangle \r| dx \\
&= \int_\rn \lf|\int_\rn \eta_\uc(x-y)\,d\mu_h(y) \r| dx \\
&\le \int_\rn \int_\rn \eta_\uc(x-y)\,dx\,d|\mu_h|(y) \\
&= |\mu_h|(\rn).
\end{align*}
Further, by \eqref{eq-hBVm}, we have
\begin{align*}
|\mu_h|(\rn)=\|(h\cdot\nabla)^mf\|_{\mathcal{M}(\rn)}
\le |h|^m\,\|f\|_{\mathringBVm(\rn)}.
\end{align*}
Combining the last two estimates yields
\begin{align*}
\|\Delta_{h}^mf\|_{L^1(\rn)}
= \lim_{\uc\to0_+} \|\Delta_{h}^mf_\uc\|_{L^1(\rn)}
\le |\mu_h|(\rn)
\le |h|^m\,\|f\|_{\mathringBVm(\rn)}.
\end{align*}
This yields \eqref{eq-Delta-BV-upper}.
\end{proof}

\begin{lemma}\label{lem-diffLpNorm}
Let $p\in[1,\fz]$ and $m\in\nn$. Then, for all $f\in C^\infty(\rn)$ and $h\in\rn$,
\begin{align}\label{eq-diff-LpNorm}
\|\Delta_h^mf-(h\cdot\nabla)^m f\|_{L^p(\rn)}
\le m \,|h|^{m+1}\, \sup_{\xi\in\msn} \|(\xi\cdot\nabla)^{m+1}f\|_{L^p(\rn)}.
\end{align}
\end{lemma}

\begin{proof}
Let $f\in C^\fz(\rn)$. For any $x,h\in\rn$, by Grafakos \cite[Proposition 1.4.5]{Grafakos-GTM250-ModernFourierAnalysis}, we know that
\begin{align*}
\Delta_h^mf(x)
= \int_{[0,1]^m} (h\cdot\nabla)^mf\Big(x+(t_1+\ldots+t_m)h \Big)\, dt_1\,\ldots\,dt_m.
\end{align*}
Combining this with the identity
\begin{align*}
&(h\cdot\nabla)^mf\Big(x+(t_1+\ldots+t_m)h \Big)-(h\cdot\nabla)^m f(x)\\
&\quad= \int_0^1 \partial_\tau \bigg((h\cdot\nabla)^mf\Big(x+\tau(t_1+\ldots+t_m)h \Big)\bigg)\,d\tau\\
&\quad= \int_0^1 (t_1+\ldots+t_m)\, \bigg((h\cdot\nabla)^{m+1}f\Big(x+\tau(t_1+\ldots+t_m)h \Big)\bigg)\,d\tau,
\end{align*}
we obtain
\begin{align*}
&|\Delta_h^mf(x)-(h\cdot\nabla)^m f(x)|\\
&\quad= \lf|\int_{[0,1]^m} \lf((h\cdot\nabla)^mf\Big(x+(t_1+\ldots+t_m)h \Big)-(h\cdot\nabla)^m f(x) \r) dt_1\,\ldots\,dt_m\r|\\
&\quad\le m\int_{[0,1]^{m+1}} \bigg|(h\cdot\nabla)^{m+1}f\Big(x+\tau(t_1+\ldots+t_m)h \Big)\bigg|\,d\tau\, dt_1\,\ldots\,dt_m.
\end{align*}
Then, applying the Minkowski inequality yields
\begin{align*}
\|\Delta_h^mf-(h\cdot\nabla)^m f\|_{L^p(\rn)}
&\le m\int_{[0,1]^{m+1}} \|(h\cdot\nabla)^{m+1}f\|_{L^p(\rn)} \,d\tau\, dt_1\,\ldots\,dt_m\\
&=m\,\|(h\cdot\nabla)^{m+1}f\|_{L^p(\rn)}\notag\\
&\le m\,|h|^{m+1} \sup_{\xi\in\msn} \|(\xi\cdot\nabla)^{m+1}f\|_{L^p(\rn)}. \notag
\end{align*}
This proves the desired estimate in \eqref{eq-diff-LpNorm}.
\end{proof}


\begin{lemma}\label{lem-norm-minus}
Let $p\in[1,\fz)$, $q\in(0,\fz)$ and $r:=\max\{q,\,1\}$. Suppose that $A\subset\rn$ is Lebesgue measurable and
$f,\,g:\,\rn\times\rn\to\rr$ satisfy
\begin{align*}
\int_A \|f(\cdot, h)\|_{L^p(\rn)}^q\,dh+\int_A \|g(\cdot, h)\|_{L^p(\rn)}^q\,dh<\infty.
\end{align*}
Then,
\begin{align}\label{eq-norm-minus}
&\lf|\lf(\int_A \|f(\cdot, h)\|_{L^p(\rn)}^q\,dh\r)^\frac{1}{r}
-\lf( \int_A \|g(\cdot, h)\|_{L^p(\rn)}^q\,dh\r)^\frac{1}{r} \r|\\
&\quad\le \lf(\int_A \|f(\cdot,h)-g(\cdot,h)\|_{L^p(\rn)}^q\,dh \r)^\frac{1}{r}. \notag
\end{align}
In particular, for any $m\in\nn$ and $0\le\varphi\in L_{\mathrm{loc}}^1(\rn\setminus\{0\})$, and any measurable functions $u,v$ on $\rn$,
\begin{align}\label{eq-norm-minus-Euv}
&\lf( \mathscr{E}_{p,\,q}^{m,\,\varphi}(u+v)\r)^\frac 1{r}\le \lf( \mathscr{E}_{p,\,q}^{m,\,\varphi}(u)\r)^\frac{1}{r}
+\lf( \mathscr{E}_{p,\,q}^{m,\,\varphi}(v)\r)^\frac 1{r}.
\end{align}
\end{lemma}

\begin{proof}
Once \eqref{eq-norm-minus} has been obtained,
we then apply it with $A=\rn$ and
\begin{align*}
f(x,h)=\Delta_h^m(u+v)(x)\left[\varphi(h)\right]^{\frac 1 q},\qquad
g(x,h)=\Delta_h^m u(x)\left[\varphi(h)\right]^{\frac 1 q}.
\end{align*}
This, together with the definition of $\mathscr{E}_{p,\,q}^{m,\,\varphi}$, yields \eqref{eq-norm-minus-Euv}.
Hence, it remains to prove \eqref{eq-norm-minus}.

By $p\in[1,\infty)$ and the Minkowski inequality, we have for almost every $h\in\rn$ that
\begin{align}\label{eq-f-g}
\lf|\|f(\cdot, h)\|_{L^p(\rn)}-\|g(\cdot, h)\|_{L^p(\rn)}\r|
\le \|f(\cdot, h)-g(\cdot, h)\|_{L^p(\rn)}
\end{align}
We now prove \eqref{eq-norm-minus} using \eqref{eq-f-g}, by considering the cases $q\in(0,1)$ and $q\in[1,\infty)$ separately.

 If $q\in(0,1)$, then $r=1$. In this case, using \eqref{eq-f-g} together with the elementary inequality
\begin{align*}
|a^q-b^q|\le |a-b|^q \qquad\text{for all} \ \, a,\,b\in(0,\infty),
\end{align*}
we obtain
\begin{align*}
 \lf| \int_A \|f(\cdot, h)\|_{L^p(\rn)}^q\,dh
-\int_A \|g(\cdot, h)\|_{L^p(\rn)}^q\,dh\r|
 &\le \int_A \lf|\|f(\cdot, h)\|_{L^p(\rn)}
 -\|g(\cdot, h)\|_{L^p(\rn)}\r|^q \,dh\\
 &\le \int_A \|f(\cdot, h)-g(\cdot, h)\|_{L^p(\rn)}^q \,dh,
\end{align*}
which proves \eqref{eq-norm-minus} in the case $q\in(0,1)$.

If $q\in[1,\infty)$, then $r=q$. In this case, applying \eqref{eq-f-g} and the Minkowski inequality yields
\begin{align*}
&\lf|\lf(\int_A \|f(\cdot,h)\|_{L^p(\rn)}^q\,dh \r)^\frac{1}{q}
- \lf(\int_A \|g(\cdot,h)\|_{L^p(\rn)}^q\,dh \r)^\frac{1}{q} \r| \\
&\quad\le \lf(\int_A \lf|\|f(\cdot,h)\|_{L^p(\rn)}
- \|g(\cdot,h)\|_{L^p(\rn)}\r|^q\,dh \r)^\frac{1}{q} \\
&\quad\le \lf(\int_A \|f(\cdot,h)-g(\cdot, h)\|_{L^p(\rn)}^q\,dh \r)^\frac{1}{q}.
\end{align*}
Thus, \eqref{eq-norm-minus} holds in the case $q\in[1,\infty)$.
\end{proof}

\subsection{Equivalent characterizations of $\MSq$ and $\BBMmq$}\label{ss-BBM-MS-cond}

The following proposition shows that the \(\BBMmq\) condition is equivalent to the higher-order analogue of the unit-mass conditions in \cite{Gennaioli-Stefani-ARMA-2026}, and also gives an analogous characterization of the \(\MSq\) condition.

\begin{proposition}\label{Prop-ARMA=BBM-MS}
Let $q\in(0,\fz)$ and $m\in\nn$. For any $j\in\nn$, suppose $0\le \varphi_j\in L_{\rm loc}^1(\rn\setminus\{0\})$,
with $\varphi_j(0)$ possibly equal to $+\infty$.
\begin{enumerate}[\rm (i)]
\item The sequence $\{\varphi_j\}_{j\in\nn}$ satisfies the $\BBMmq$ condition if and only if
\begin{align}\label{eq-ARMA1-mq}
\mathcal{N}_{m,\,q}:=\sup_{r\in(0,\fz)} \limsup_{j\to\fz} \int_\rn \varphi_j(h) \min\left\{|h|^{mq},\,r^{mq} \right\} dh<\fz
\end{align}
and, moreover, for any $\psi\in C_{\mathrm{c}}(\rn)$,
\begin{align}\label{eq-ARMA2-mq}
\lim_{j\to\fz} \int_\rn \psi(h)\,|h|^{mq}\,\varphi_j(h)\,dh
= \psi(0).
\end{align}

\item The sequence $\{\varphi_j\}_{j\in\nn}$ satisfies the $\MSq$ condition if and only if
\begin{align}\label{eq-ARMA1-MS}
\mathcal{N}_q:=\sup_{r\in(0,\fz)} \limsup_{j\to\fz} \int_\rn \varphi_j(h) \min\lf\{\frac{|h|^q}{r^q},\,1 \r\} dh<\fz
\end{align}
and, moreover, for any $\psi\in C_{\mathrm{c}}(\rn)$,
\begin{align}\label{eq-ARMA2-MS}
\lim_{j\to\fz} \int_\rn \psi\lf(\frac{h}{|h|^2} \r) \varphi_j(h)\,dh=\psi(0).
\end{align}

\end{enumerate}
\end{proposition}

\begin{remark}\label{rem-vague-conv-BBM-MS}
Set $d\nu_j(h):=|h|^{mq}\,\varphi_j(h)\,dh$. Condition \eqref{eq-ARMA2-mq} means that $\nu_j$ converges vaguely to the Dirac measure
$\delta_0$ as $j\to\infty$. Condition \eqref{eq-ARMA2-MS} admits an analogous interpretation. Indeed, define
\begin{align*}
I(h):=\frac{h}{|h|^2}
\qquad \text{for all}\ \, h\in\rn\setminus\{0\},
\end{align*}
and let $\widetilde\nu_j$ be the pushforward of the measure $\varphi_j(h)\,dh$ under $I$, namely,
\begin{align*}
\widetilde\nu_j:=I_\#\left(\varphi_j(h)\,dh\right).
\end{align*}
 Then \eqref{eq-ARMA2-MS} is equivalent to the vague convergence of $\widetilde\nu_j$ to $\delta_0$ as $j\to\infty$.
\end{remark}

\begin{proof}[Proof of Proposition \ref{Prop-ARMA=BBM-MS}]
We divide the proof into four parts.

\medskip

\textbf{Part 1: Proof that $\BBMmq$ implies
\eqref{eq-ARMA1-mq} and \eqref{eq-ARMA2-mq}.}
For any $r\in(0,\infty)$, we have
\begin{align*}
&\int_\rn \varphi_j(h) \min\lf\{|h|^{mq},\,r^{mq} \r\} dh\\
&\quad \le \int_{|h|<1} |h|^{mq}\,\varphi_j(h)\,dh
+ r^{mq} \int_{|h|\ge 1} \varphi_j(h)\,dh \\
&\quad = \int_{|h|<1} |h|^{mq}\,\varphi_j(h)\,dh
+ r^{mq}\int_{|h|\ge 1} \varphi_j(h) \min\lf\{|h|^{mq},\,1 \r\} dh.
\end{align*}
By $\NBBMmq$ and $\VBBMmq$, the right-hand side converges to $1$
as $j\to\infty$. Consequently,
\begin{align*}
\mathcal{N}_{m,\,q}
= \sup_{r\in(0,\infty)} \limsup_{j\to\infty} \int_\rn \varphi_j(h) \min\lf\{|h|^{mq},\,r^{mq} \r\} dh\le1<\infty.
\end{align*}
This proves \eqref{eq-ARMA1-mq}.

Let $\psi\in C_{\rm c}(\rn)$ such that $\supp\psi\subset B(0,R)$
for some $R\in(1,\fz)$. For any $\uc\in(0,1)$, there exists $\delta\in(0,1)$ such that
$|\psi(h)-\psi(0)|\le\uc$ for all $h\in B(0,\delta)$. Consequently,
\begin{align}\label{eq1-BBM-ARMA}
&\lf|\int_\rn \psi(h)\,|h|^{mq}\,\varphi_j(h)\,dh
- \psi(0)\,\int_{|h|<\delta}|h|^{mq}\,\varphi_j(h)\,dh \r| \\
&\quad\le \int_{|h|<\delta} |\psi(h)-\psi(0)|\,|h|^{mq}\, \varphi_j(h)\,dh
+ \int_{\delta\le|h|\le R} |\psi(h)|\,|h|^{mq}\, \varphi_j(h)\,dh\notag\\
&\quad \le \uc \int_{|h|<\delta} |h|^{mq}\,\varphi_j(h)\,dh + \int_{\delta\le|h|\le R} |\psi(h)|\,|h|^{mq}\, \varphi_j(h)\,dh. \notag
\end{align}
For $\delta\le|h|<R$, we have
\begin{align*}
|h|^{mq}<R^{mq}= \lf(\frac{R}{\delta} \r)^{mq}\delta^{mq}
\le \lf(\frac{R}{\delta} \r)^{mq}
\min \lf\{|h|^{mq},\,1 \r\},
\end{align*}
which, together with $\VBBMmq$, yields
\begin{align*}
&\lim_{j\to\fz}\int_{\delta\le|h|\le R} |\psi(h)|\,|h|^{mq}\,\varphi_j(h)\,dh\\
&\quad\le \|\psi\|_{L^\fz(\rn)} \lf(\frac{R}{\delta} \r)^{mq} \lim_{j\to\fz}\int_{|h|\ge\delta} \varphi_j(h) \min\lf\{|h|^{mq},\,1 \r\} dh \\
&\quad=0.
\end{align*}
Moreover, $\NBBMmq$ and $\VBBMmq$ imply that
\begin{align*}
\lim_{j\to\infty}
\int_{|h|<\delta}|h|^{mq}\,\varphi_j(h)\,dh=1.
\end{align*}
Consequently, letting $j\to\infty$ on both sides of \eqref{eq1-BBM-ARMA}, we obtain
\begin{align*}
\limsup_{j\to\fz} \lf|\int_\rn \psi(h)\,|h|^{mq}\,\varphi_j(h)\,dh-\psi(0) \r|
\le\uc.
\end{align*}
Since $\uc\in(0,1)$ is arbitrary, this yields \eqref{eq-ARMA2-mq}.

\medskip

\textbf{Part 2: Proof that \eqref{eq-ARMA1-mq} and \eqref{eq-ARMA2-mq} imply $\BBMmq$.}
Choose $\psi_1,\psi_2 \in C_{\rm c}(\rn)$ such that
\begin{align}\label{eq-psi12}
0\le\psi_1,
\psi_2\le1,\quad
\supp\psi_1\subset B(0,1),\quad
\psi_1(0)=1,\quad \psi_2\equiv1\ \, \text{on}\ \, B(0,1).
\end{align}
Then,
\begin{align*}
\int_\rn \psi_1(h)\,|h|^{mq}\,\varphi_j(h)\,dh
\le \int_{|h|<1} |h|^{mq}\,\varphi_j(h)\,dh
\le \int_\rn \psi_2(h)\,|h|^{mq}\,\varphi_j(h)\,dh.
\end{align*}
By \eqref{eq-ARMA2-mq} and $\psi_1(0)=\psi_2(0)=1$, the first and third terms in the above inequality tend to $1$ as $j\to\fz$. Thus,
\begin{align*}
\lim_{j\to\infty} \int_{|h|<1}|h|^{mq}\,\varphi_j(h)\,dh=1,
\end{align*}
which establishes $\NBBMmq$.

Now, we prove $\VBBMmq$. For any $R>1$, we have
\begin{align*}
\int_{|h|\ge R} \varphi_j(h) \min\lf\{|h|^{mq},\,1 \r\} dh
&=\int_{|h|\ge R}\varphi_j(h)\,dh\\
&\le R^{-mq} \int_\rn \varphi_j(h) \min\lf\{|h|^{mq},\,R^{mq} \r\} dh.
\end{align*}
Consequently, by \eqref{eq-ARMA1-mq},
\begin{align}\label{eq-ARME-to-VBBM}
\limsup_{j\to\infty} \int_{|h|\ge R} \varphi_j(h) \min\lf\{|h|^{mq},\,1 \r\} dh
\le \frac{\mathcal{N}_{m,\,q}}{R^{mq}}.
\end{align}
Fix $\delta\in(0,1]$.
Choose $0\le \psi\in C_{\mathrm c}(\rn)$ such that
$\psi(0)=0$ and
$\psi(h)=1$ when $\delta\le|h|\le R$.
By \eqref{eq-ARMA2-mq}, we have
\begin{align*}
\int_{\delta\le|h|<R} \varphi_j(h) \min\lf\{|h|^{mq},\,1 \r\} dh
\le \int_\rn \psi(h)\,|h|^{mq}\,\varphi_j(h)\,dh,
\end{align*}
whose right-hand side tends to
$\psi(0)=0$
as $j\to\infty$. Combining this with \eqref{eq-ARME-to-VBBM}, we obtain
\begin{align*}
\limsup_{j\to\infty} \int_{|h|\ge\delta} \varphi_j(h) \min\lf\{|h|^{mq},\,1 \r\} dh
&= \limsup_{j\to\infty} \left(\int_{\delta\le|h|<R}+\int_{|h|\ge R} \right)\cdots \\
&\le \frac{\mathcal{N}_{m,\,q}}{R^{mq}}.
\end{align*}
Letting $R\to\infty$, we conclude that
\begin{align*}
\lim_{j\to\infty} \int_{|h|\ge\delta} \varphi_j(h) \min\lf\{|h|^{mq},\,1 \r\} dh=0.
\end{align*}
Thus $\VBBMmq$ holds, thereby completing the proof of (i).

\medskip

\textbf{Part 3: Proof that $\MSq$ implies
\eqref{eq-ARMA1-MS} and \eqref{eq-ARMA2-MS}.}
For every fixed $r>0$, we have
\begin{align*}
\int_{\rn}\varphi_j(h)
\min\left\{\frac{|h|^q}{r^q},\,1\right\}\,dh
\le \frac1{r^q}\int_{|h|<1}|h|^q\,\varphi_j(h)\,dh
+\int_{|h|\ge1}\varphi_j(h)\,dh.
\end{align*}
By $\NMSq$ and $\VMSq$, the right-hand side
converges to $1$. Thus $\mathcal{N}_q\le1$.

Let $\psi\in C_{\mathrm{c}}(\rn)$ such that $\operatorname{supp}\psi\subset B(0,R)$ for some $R\in(1,\fz)$. Given $\varepsilon\in(0,1)$, the continuity of $\psi$ at the origin yields some $M\in(R,\fz)$ such that for all $|h|\ge M$,
\begin{align*}
\left|\psi\left(\frac{h}{|h|^2} \right)-\psi(0) \right| \le \varepsilon.
\end{align*}
Since $\psi(h/|h|^2)=0$ whenever $0<|h|\le1/R$, it follows that
\begin{align}\label{eq-varphij}
&\left|\int_\rn \psi\left(\frac{h}{|h|^2} \right) \varphi_j(h)\,dh
- \psi(0) \int_{|h|\ge M} \varphi_j(h)\,dh \right| \\
&\quad\le \|\psi\|_{L^\fz(\rn)} \int_{1/R<|h|<M} \varphi_j(h)\,dh
+ \int_{|h|\ge M} \left|\psi\left(\frac{h}{|h|^2} \right)-\psi(0) \right| \varphi_j(h)\,dh \notag\\
&\quad\le R^q\,\|\psi\|_{L^\fz(\rn)} \int_{|h|<M} \varphi_j(h) \min\left\{|h|^q,\,1 \right\} dh
+ \varepsilon \int_{|h|\ge M} \varphi_j(h)\,dh, \notag
\end{align}
where in the last step we used the fact that $R\,\min\{|h|,1\}\ge 1$ for $|h|>1/R$.
By $\MSq$, we have
\begin{align}\label{eq-MS-tail-R}
\lim_{j\to\fz} \int_{|h|\ge M} \varphi_j(h)\,dh=1.
\end{align}
Taking the limit $j\to\infty$ on both sides of \eqref{eq-varphij}, together with $\VMSq$ and \eqref{eq-MS-tail-R}, we obtain
\begin{align*}
\limsup_{j\to\fz} \left|\int_\rn \psi\left(\frac{h}{|h|^2} \right)\varphi_j(h)\,dh-\psi(0) \right|\le\varepsilon.
\end{align*}
Letting $\varepsilon\to0$ yields \eqref{eq-ARMA2-MS}.

\medskip

\textbf{Part 4: Proof that
\eqref{eq-ARMA1-MS} and \eqref{eq-ARMA2-MS} imply $\MSq$.}
Let $\psi_1,\psi_2\in C_{\mathrm{c}}(\rn)$ satisfy \eqref{eq-psi12}.
Then
\begin{align*}
\int_\rn \psi_1\left(\frac{h}{|h|^2} \right) \varphi_j(h)\,dh
\le \int_{|h|\ge1} \varphi_j(h)\,dh
\le \int_\rn \psi_2\left(\frac{h}{|h|^2} \right)\varphi_j(h)\,dh.
\end{align*}
By \eqref{eq-ARMA2-MS} and $\psi_1(0)=\psi_2(0)=1$, both outer terms converge to $1$ as $j\to\infty$. This proves $\NMSq$.

To prove $\VMSq$, we fix $\delta\in[1,\fz)$ and $r\in(0,1)$. Since
\begin{align*}
\int_{|h|<r} \varphi_j(h) \min\left\{|h|^q,\,1 \right\}\, dh
=\int_{|h|<r}\varphi_j(h)\, |h|^q\,dh
\le r^q \int_\rn \varphi_j(h) \min\left\{\frac{|h|^q}{r^q},\,1 \right\} dh,
\end{align*}
it follows from \eqref{eq-ARMA1-MS} that
\begin{align}\label{eq-MS-small-ball}
\limsup_{j\to\fz} \int_{|h|<r} \varphi_j(h) \min\left\{|h|^q,\,1 \right\}\,dh
\le \mathcal{N}_q\, r^q.
\end{align}
Choose $0\le \psi\in C_{\mathrm{c}}(\rn)$ such that $\psi(0)=0$ and $\psi(h)=1$ when $\delta^{-1}\le|h|\le r^{-1}$. Then $\psi(h/|h|^2)=1$ when $r\le|h|<\delta$. Consequently,
\begin{align*}
\int_{r\le|h|<\delta} \varphi_j(h) \min\left\{|h|^q,\,1 \right\} dh
\le\int_\rn \psi\left(\frac{h}{|h|^2} \right) \varphi_j(h)\,dh
\to \psi(0)=0
\end{align*}
as $j\to\infty$ by \eqref{eq-ARMA2-MS}. Combining this with \eqref{eq-MS-small-ball} yields
\begin{align*}
\limsup_{j\to\fz} \int_{|h|<\delta} \varphi_j(h) \min\left\{|h|^q,\,1 \right\} dh
\le \mathcal{N}_q\, r^q
\qquad\text{for all } r\in(0,1).
\end{align*}
Letting $r\to0$ gives $\VMSq$. This completes the proof of (ii).
\end{proof}

\section{Standard $\Gamma$-convergence implies pointwise convergence}\label{sec-Gamma-to-point}

\subsection{An admissible condition}\label{sec-Perliminary3}

The following definition of admissible condition was introduced in \cite{Liu-Xiao-Zhang-2026-submit,Liu-Zhang-2026-submit},
which is indispensable for the study of non-trivial energy spaces.

\begin{definition}\label{def-adm}
Let $\alpha\in(0,\fz)$.
\begin{enumerate}[{\rm (i)}]
\item A nonnegative function $\varphi$ that is locally integrable on $\rn \setminus\{0\}$, with $\varphi(0)$ possibly equal to $+\infty$, is called \emph{$\alpha$-admissible} if
\begin{align*}
\int_{\rn} \varphi(h) \min\left\{|h|^{\alpha},\,1 \right\} dh<\fz.
\end{align*}

\item For any $j\in\nn$, let $0\le \varphi_j\in L_{\mathrm{loc}}^1(\rn \setminus\{0\})$, with $\varphi_j(0)$ possibly equal to $+\infty$.
The sequence $\{\varphi_j\}_{j\in\nn}$ is called \emph{uniformly $\alpha$-admissible} if
\begin{align*}
\limsup_{j\to\fz} \int_\rn \varphi_j(h) \min\lf\{|h|^{\alpha},\,1 \r\} dh<\fz.
\end{align*}
\end{enumerate}
\end{definition}

According to \cite[Proposition 2.10] {Liu-Zhang-2026-submit} (see also \cite[Remark 1.3]
{Liu-Xiao-Zhang-2026-submit} for the case $m=1$), when $\varphi$ is radial, every function in $C_{\rm c}^\fz(\rn)$
has finite $\mathscr{E}_{p,\,q}^{m,\,\varphi}$-energy if and only if $\varphi$ is $mq$-admissible.
Therefore, admissibility is essential for studying a non-trivial energy family.

\begin{proposition}\label{prop-BBM-MS=>adm}
Let $m\in\nn$ and $q\in(0,\fz)$. For any $j\in\nn$, assume that $0\le \varphi_j\in L_{\mathrm{loc}}^1(\rn \setminus\{0\})$,
with $\varphi_j(0)$ possibly equal to $+\infty$.
\begin{enumerate}[{\rm (i)}]
\item If $\{\varphi_j\}_{j\in\nn}$ satisfies the $\MSq$ condition, then it is uniformly $q$-admissible.

\item If $\{\varphi_j\}_{j\in\nn}$ satisfies the $\BBMmq$ condition, then
it is uniformly $mq$-admissible.
\end{enumerate}
\end{proposition}

\begin{proof}
The desired results are immediate consequences of the definitions of $\MSq$, $\BBMmq$, and uniform admissibility.
\end{proof}

In addition to Proposition \ref{prop-BBM-MS=>adm}, each statement in Theorems \ref{thm-MS} and \ref{thm-BBM} separately
guarantees the uniform admissibility of the full MS sequence or the indicated BBM subsequence. To prove this assertion,
we need the following lemma.

\begin{lemma}\label{Lem-Gamma-bdd}
Let $p\in[1,\fz)$, $q\in(0,\fz)$ and $m\in\nn$. For any $j\in\nn$, assume that $0\le \varphi_j\in L_{\mathrm{loc}}^1(\rn \setminus\{0\})$,
with $\varphi_j(0)$ possibly equal to $+\infty$. If $0\not\equiv f\in L^p(\rn)$ and $\{f_j\}_{j\in\nn}\subset L^p(\rn)$ satisfy
\begin{align*}
\lim_{j\to\infty}\|f_j-f\|_{L^p(\rn)}=0,
\end{align*}
then
\begin{align}\label{eq-Gamma=>adm}
\limsup_{j\to\infty}\int_\rn \varphi_j(h) \min\lf\{|h|^{mq},\,1 \r\} dh
\le C\,\limsup_{j\to\infty} \mathscr{E}_{p,\,q}^{m,\,\varphi_j}(f_j),
\end{align}
where $C$ is a positive constant depending on $m,n,p,q$, and $f$.
\end{lemma}

\begin{proof}
Since $0\not\equiv f\in L^p(\rn)$, there exists $R\in(1,\fz)$ such that $\|f{\bf 1}_{B(0,\,R)}\|_{L^p(\rn)}>0$ and
\begin{align*}
\|f{\bf 1}_{\rn\setminus B(0,\,R)}\|_{L^p(\rn)}
\le 2^{-(m+2)}\,\|f{\bf 1}_{B(0,\,R)}\|_{L^p(\rn)}.
\end{align*}
Fix such $R$ and write $g:=f{\bf 1}_{B(0,\,R)}$. By $f_j\to f$ in $L^p(\rn)$, there exists $J\in\nn$ such that, for every $j\ge J$,
\begin{align*}
\|f-f_j\|_{L^p(\rn)}
\le 2^{-(m+2)}\,\|g\|_{L^p(\rn)}.
\end{align*}
Combining the above two inequalities yields that for all $j\ge J$,
\begin{align}\label{eq-g-fj}
\|g-f_j\|_{L^p(\rn)}
\le 2^{-(m+1)}\,\|g\|_{L^p(\rn)}.
\end{align}

Now, fix $0\neq h\in\rn$ and set
$
N:=\lfloor{2R}/{|h|} \rfloor+1.
$
Then, we have
\begin{align}\label{eq-N<h}
N\min\{|h|,\,1\} \le \left(\frac{2R}{|h|}+2\right)\min\{|h|,\,1\}
\le2R+2.
\end{align}
On the one hand, since $\supp g\subset\overline{B(0,R)}$, it follows from Lemma \ref{Lem-LXlem2829}(iv) that
\begin{align*}
\|\Delta_{Nh}^mg\|_{L^p(\rn)}^p
= \sum_{j=0}^m \binom{m}{j}^p\,\|g\|_{L^p(\rn)}^p
\ge \|g\|_{L^p(\rn)}^p.
\end{align*}
On the other hand, applying \eqref{eq-g-fj} and Lemma \ref{Lem-LXlem2829}(i) yields
\begin{align*}
\|\Delta_{Nh}^m(g-f_j)\|_{L^p(\rn)} \le 2^m\,\|g-f_j\|_{L^p(\rn)} \le \frac{1}{2}\,\|g\|_{L^p(\rn)}.
\end{align*}
A combination of the last two estimates gives
\begin{align}\label{eq-Delta-Nh1}
\|\Delta_{Nh}^mf_j\|_{L^p(\rn)}
\ge \|\Delta_{Nh}^mg\|_{L^p(\rn)}-\|\Delta_{Nh}^m(g-f_j)\|_{L^p(\rn)}\ge \frac{1}{2}\,\|g\|_{L^p(\rn)}.
\end{align}
Meanwhile, it has been proved in \cite[(3.1)]{Triebel-ArchMath-2011} that
\begin{align*}
\|\Delta_{Nh}^mf_j\|_{L^p(\rn)}
\le N^m\,\|\Delta_h^mf_j\|_{L^p(\rn)}.
\end{align*}
Inserting this into \eqref{eq-Delta-Nh1} and using \eqref{eq-N<h}, we obtain
\begin{align*}
\|\Delta_h^mf_j\|_{L^p(\rn)}
\ge \frac{1}{2\,N^m}\,\|g\|_{L^p(\rn)}
\ge \frac{\|g\|_{L^p(\rn)}}{2^{m+1}\,(R+1)^m}\,\min\lf\{|h|^m,\,1 \r\}.
\end{align*}
This last inequality is trivially true for $h=0$. Thus,
\begin{align*}
\int_\rn \varphi_j(h) \min\lf\{|h|^{mq},\,1 \r\} dh
&\le \lf(\frac{2^{m+1}\,(R+1)^m}{\|g\|_{L^p(\rn)}} \r)^q\int_\rn \|\Delta_h^mf_j\|_{L^p(\rn)}^q\,\varphi_j(h)\,dh\\
&= \lf(\frac{2^{m+1}\,(R+1)^m}{\|g\|_{L^p(\rn)}} \r)^q \mathscr{E}_{p,\,q}^{m,\,\varphi_j}(f_j).
\end{align*}
Taking the limsup as $j\to\fz$ on both sides yields \eqref{eq-Gamma=>adm}. This completes the proof of Lemma \ref{Lem-Gamma-bdd}.
\end{proof}

\begin{proposition}\label{prop-Gamma-BBM-MS=>adm}
Let $p\in[1,\fz)$, $q\in(0,\fz)$ and $m\in\nn$. For any $j\in\nn$, assume that $0\le \varphi_j\in
L_{\rm loc}^1(\rn \setminus\{0\})$, with $\varphi_j(0)$ possibly equal to $+\infty$.
\begin{enumerate}
\item[{\rm (a)}] Any assertion (i)--(iv) of Theorem \ref{thm-MS} independently ensures that
$\{\varphi_j\}_{j\in\nn}$ is uniformly $q$-admissible.

\item[{\rm (b)}] Any assertion (i)--(iv) of Theorem \ref{thm-BBM} independently ensures that
$\{\varphi_{j_k}\}_{k\in\nn}$ is uniformly $mq$-admissible.
\end{enumerate}
In both items, the equivalences among the assertions (i)--(iv) of Theorem \ref{thm-MS} or Theorem \ref{thm-BBM} are not assumed a priori.
\end{proposition}

\begin{proof}
If assertion (i) of Theorem \ref{thm-MS} or Theorem \ref{thm-BBM} holds, then the desired uniform admissibility follows
from Proposition \ref{prop-BBM-MS=>adm}.

Assume that assertion (iii) of Theorem \ref{thm-MS} (resp. Theorem \ref{thm-BBM}) holds. Then, by the standard $\Gamma$-convergence condition $\mathbf{(\Gamma)_1}$, for any $0\not\equiv f\in C_{\rm c}^\infty(\rn)$ there exists a sequence $\{f_j\}_{j\in\nn}\subset L^p(\rn)$ converging to $f$ in $L^p(\rn)$ such that
\begin{align}\label{eq-limsup<Df}
\limsup_{j\to\infty} \mathscr{E}_{p,\,q}^{m,\,\varphi_j}(f_j)
\le \mathscr{D}(f),
\end{align}
where we take $m=1$ and $\mathscr{D}=\mathscr{D}_{p,\,q}^0$ when assuming Theorem \ref{thm-MS}(iii), and general $m\in\nn$ and $\mathscr{D}=\mathscr{D}_{p,\,q}^{m,\,\mu}$ when assuming Theorem \ref{thm-BBM}(iii). Combining this with Lemma \ref{Lem-Gamma-bdd}, we obtain
\begin{align*}
\limsup_{j\to\infty}\int_\rn \varphi_j(h) \min\left\{|h|^{mq},\,1 \right\} dh
\le C\,\mathscr{D}(f)
<\infty.
\end{align*}
Thus, the sequence $\{\varphi_j\}_{j\in\nn}$ is uniformly $q$-admissible (resp. $mq$-admissible).

If any of assertions (ii) and (iv) in Theorems \ref{thm-MS} and \ref{thm-BBM} holds, then by either pointwise convergence or strong $\Gamma$-convergence condition $\mathbf{(s\Gamma)_1}$, the limsup inequality \eqref{eq-limsup<Df} is satisfied by taking every $f_j\equiv f$.
Again, the desired uniform admissibility follows from Lemma \ref{Lem-Gamma-bdd}.
\end{proof}

\subsection{A technical estimate}\label{sec3.2}

The following proposition will be used in the proof of (iii) $\Rightarrow$ (ii) in both Theorem~\ref{thm-MS} and Theorem~\ref{thm-BBM}.

\begin{proposition}\label{prop-BBM-Gamma-to-cond}
Let $p\in[1,\fz)$, $q\in(0,\fz)$ and $m\in\nn$. For any $j\in\nn$, assume that $0\le \varphi_j\in L_{\rm loc}^1(\rn \setminus\{0\})$,
with $\varphi_j(0)$ possibly equal to $+\infty$. For any $f\in W^{m,\,p}(\rn)$, if $\{f_j\}_{j\in\nn}\subset L^p(\rn)$ converges to $f$ in $L^p(\rn)$ and
\begin{align}\label{eq-bdd-E(fj)}
\limsup_{j\to\fz} \mathscr{E}_{p,\,q}^{m,\,\varphi_j}(f_j)<\fz,
\end{align}
then
\begin{align}\label{eq-Gamma1-fj}
\limsup_{j\to\fz} \mathscr{E}_{p,\,q}^{m,\,\varphi_j}(f)
\le \limsup_{j\to\fz} \mathscr{E}_{p,\,q}^{m,\,\varphi_j}(f_j).
\end{align}
Moreover, if $p=1$ and $f\in\BVm(\rn)$, then \eqref{eq-Gamma1-fj} remains valid.
\end{proposition}

The key ingredient in the proof of Proposition~\ref{prop-BBM-Gamma-to-cond} is the following technical lemma.

\begin{lemma}\label{Lem-limit=0}
Let $p\in[1,\fz)$ and $m\in\nn$. If $f\in W^{m,\,p}(\rn)$ and $\{f_j\}_{j\in\nn}\subset L^p(\rn)$
converges to $f$ in $L^p(\rn)$, then
\begin{align}\label{eq-Ij-claim}
\lim_{j\to\fz} \lf(\sup_{h\in\rn\setminus\{0\}} \frac{\|\Delta_h^mf\|_{L^p(\rn)}
- \|\Delta_h^mf_j\|_{L^p(\rn)}}{\min\{|h|^m,\,1\}} \r)=0.
\end{align}
Moreover, if $p=1$ and $f\in\BVm(\rn)$, then \eqref{eq-Ij-claim} remains valid.
\end{lemma}

\begin{proof}
Let $p\in[1,\fz)$, $f\in L^p(\rn)$, and $\{f_j\}_{j\in\nn}\subset L^p(\rn)$ be a sequence converging to $f$ in
$L^p(\rn)$. For any $h_0\in\rn$,
applying Lemma~\ref{Lem-LXlem2829}(i) gives
\begin{align*}
\|\Delta_{h_0}^mf_j\|_{L^p(\rn)}
\le \|\Delta_{h_0}^mf\|_{L^p(\rn)}+\|\Delta_{h_0}^m(f_j-f)\|_{L^p(\rn)}
\le \|\Delta_{h_0}^mf\|_{L^p(\rn)}+2^m\,\|f_j-f\|_{L^p(\rn)},
\end{align*}
which implies
\begin{align*}
\sup_{h\in\rn\setminus\{0\}} \frac{\|\Delta_h^mf\|_{L^p(\rn)}-\|\Delta_h^mf_j\|_{L^p(\rn)}}{\min\{|h|^m,\,1\}}
&\ge \|\Delta_{h_0}^mf\|_{L^p(\rn)}-\|\Delta_{h_0}^mf_j\|_{L^p(\rn)}
\ge -2^m\,\|f-f_j\|_{L^p(\rn)}.
\end{align*}
Letting $j\to\infty$ yields
\begin{align}\label{eq-liminf>0}
\liminf_{j\to\fz} \lf(\sup_{h\in\rn\setminus\{0\}} \frac{\|\Delta_h^mf\|_{L^p(\rn)}
- \|\Delta_h^mf_j\|_{L^p(\rn)}}{\min\{|h|^m,\,1\}} \r)\ge 0.
\end{align}
Moreover, for any $\delta\in(0,1)$, applying Lemma \ref{Lem-LXlem2829}(i) implies
\begin{align*}
\lim_{j\to\fz} \lf(\sup_{|h|\ge\delta} \frac{\|\Delta_h^mf\|_{L^p(\rn)}
- \|\Delta_h^mf_j\|_{L^p(\rn)}}{\min\{|h|^m,\,1\}} \r)
&\le \lim_{j\to\fz} \lf(\sup_{|h|\ge\delta} \frac{\|\Delta_h^mf-\Delta_h^m f_j\|_{L^p(\rn)}}{\min\{|h|^m,\,1\}} \r) \\
&\le 2^m\,\delta^{-m} \lim_{j\to\fz} \|f-f_j\|_{L^p(\rn)} \notag\\
&=0 \notag
\end{align*}
and, hence,
\begin{align}\label{eq-h-large}
\lim_{\delta\to 0_+}\lim_{j\to\fz} \lf(\sup_{|h|\ge\delta} \frac{\|\Delta_h^mf\|_{L^p(\rn)}
- \|\Delta_h^mf_j\|_{L^p(\rn)}}{\min\{|h|^m,\,1\}} \r)=0.
\end{align}

Now, suppose that either $f\in W^{m,\,p}(\rn)$ with $p\in(1,\infty)$,
or $f\in \BVm(\rn)$ with $p=1$. Assume for the moment that we have obtained
\begin{align}\label{eq-h-small}
\lim_{\delta\to 0_+} \limsup_{j\to\fz} \lf(\sup_{0<|h|<\delta} \frac{\|\Delta_h^mf\|_{L^p(\rn)}
- \|\Delta_h^mf_j\|_{L^p(\rn)}}{\min\{|h|^m,\,1\}} \r)=0.
\end{align}
Then \eqref{eq-Ij-claim} follows directly from \eqref{eq-liminf>0}, \eqref{eq-h-large} and \eqref{eq-h-small}.

Note that the proof of \eqref{eq-liminf>0} also implies that the limit
in \eqref{eq-h-small} is nonnegative. Suppose, to the contrary, that
\eqref{eq-h-small} fails. Then there exists a constant $c\in(0,\infty)$
such that
\begin{align*}
\lim_{\delta\to 0_+}\limsup_{j\to\infty} \left(\sup_{0<|h|<\delta} \frac{\|\Delta_h^m f\|_{L^p(\mathbb{R}^n)}
- \|\Delta_h^m f_j\|_{L^p(\mathbb{R}^n)}}{\min\{|h|^m,1\}} \right)>c>0.
\end{align*}
Thus, for a sequence $\{\delta_j\}_{j\in\nn}$ that decreases to $0$, by picking a
sufficiently large index $j$ with the corresponding supremum exceeding
$c$, we obtain a sequence $\{h_j\}_{j\in\mathbb N}$
converging to $0$ such that, for all $j\in\mathbb N$,
\begin{align}\label{eq-fail}
\frac{\|\Delta_{h_j}^m f\|_{L^p(\mathbb{R}^n)}
- \|\Delta_{h_j}^m f_j\|_{L^p(\mathbb{R}^n)}}{|h_j|^m}\ge c>0.
\end{align}
Since $\msn$ is compact, passing to a subsequence we may assume that
$h_j/|h_j|\to\xi_0$ for some $\xi_0\in \msn$ as $j\to\infty$.  We will
show that if $f\in W^{m,\,p}(\rn)$ with $p\in[1,\infty)$, or
$f\in \BVm(\rn)$ with $p=1$, then
\begin{align}\label{eq-fail-lim}
\limsup_{j\to\infty} \frac{\|\Delta_{h_j}^mf\|_{L^p(\rn)}}{|h_j|^m}
\le \liminf_{j\to\infty} \frac{\|\Delta_{h_j}^m f_j\|_{L^p(\rn)}}{|h_j|^m}.
\end{align}
This contradicts \eqref{eq-fail}, which proves \eqref{eq-h-small} and
hence \eqref{eq-Ij-claim}.
The remainder of the proof is divided into two steps.

\medskip

{\bf Step 1: Proof of \eqref{eq-fail-lim} for $p\in (1,\fz)$ and $f\in W^{m,\,p}(\rn)$.}
By \eqref{eq-def-xi-cdot-nabla}, we have
\begin{align*}
\lf\|\lf(\frac{h_j}{|h_j|}\cdot\nabla \r)^mf-(\xi_0\cdot\nabla)^mf \r\|_{L^p(\rn)}
\le \sum_{\{\alpha\in\zz_+^n:\,|\alpha|=m\}} \frac{m!}{\alpha!} \lf|\lf(\frac{h_j}{|h_j|} \r)^\alpha-\xi_0^\alpha \r| \|\partial^\alpha f\|_{L^p(\rn)}.
\end{align*}
Since $h_j/|h_j|\to\xi_0$ as $j\to\infty$, it follows that
\begin{align}\label{eq-hj-xi0}
\lim_{j\to\infty}\lf\|\lf(\frac{h_j}{|h_j|}\cdot\nabla \r)^mf-(\xi_0\cdot\nabla)^mf \r\|_{L^p(\rn)}=0.
\end{align}
This, along with \eqref{eq-Delta-Lp-upper} in Lemma \ref{Lem-LXlem2829}(ii), implies
\begin{align}\label{eq-limsup-le}
\limsup_{j\to\fz}\,\frac{\|\Delta_{h_j}^mf\|_{L^p(\rn)}}{|h_j|^m}
\le \limsup_{j\to\fz} \lf\|\lf(\frac{h_j}{|h_j|}\cdot\nabla \r)^mf \r\|_{L^p(\rn)}
= \|(\xi_0\cdot\nabla)^mf\|_{L^p(\rn)}.
\end{align}

Denote by $p'$ the H\"older conjugate exponent of $p$. By duality and
integration by parts,
\begin{align}\label{eq-dual}
\|(\xi_0\cdot\nabla)^mf\|_{L^p(\rn)}
&= \sup_{\gfz{\psi\in C_{\rm c}^\fz(\rn)}{\|\psi\|_{L^{p'}(\rn)}\le1}}
\left|\int_\rn (\xi_0\cdot\nabla)^mf(x)\,\psi(x)\,dx\right|\\
&= \sup_{\gfz{\psi\in C_{\rm c}^\fz(\rn)}{\|\psi\|_{L^{p'}(\rn)}\le1}}
\left|(-1)^m\int_\rn f(x)\,(\xi_0\cdot\nabla)^m\psi(x)\,dx\right|. \notag
\end{align}
Next, we treat the right-hand side of \eqref{eq-dual} and show that it
is controlled by the liminf in \eqref{eq-fail-lim}.

Set $\tilde\xi_0:=-\xi_0$ and $\tilde h_j:=-h_j$ for all $j\in\nn$. For
any $\psi\in C_{\rm c}^\fz(\rn)$, by Lemma~\ref{lem-diffLpNorm} and the
convergence $|\tilde h_j|=|h_j|\to0$, we obtain
\begin{align*}
\lim_{j\to\infty}\left\|
\left(\frac{\tilde h_j}{|h_j|}\cdot\nabla\right)^m\psi-
\frac{\Delta_{\tilde h_j}^m\psi}{|h_j|^m} \right\|_{L^{p'}(\rn)}
&\le \lim_{j\to\infty}
\left(m |h_j| \sup_{\xi\in\msn} \|(\xi\cdot\nabla)^{m+1}\psi\|_{L^{p'}(\rn)}\right)\\
&=0.
\end{align*}
Meanwhile, the proof of \eqref{eq-hj-xi0} also implies
\begin{align*}
\lim_{j\to\infty}\left\|
\left(\frac{\tilde h_j}{|h_j|}\cdot\nabla\right)^m\psi-
(\tilde\xi_0\cdot\nabla)^m\psi \right\|_{L^{p'}(\rn)}=0.
\end{align*}
Combining these two limits with the Minkowski inequality yields
\begin{align}\label{eq-diff-LpNorm-psi-tilde}
\lim_{j\to\infty}\left\|
(\tilde\xi_0\cdot\nabla)^m\psi
-\frac{\Delta_{\tilde h_j}^m\psi}{|h_j|^m}\right\|_{L^{p'}(\rn)}=0.
\end{align}

For any $j\in\nn$, we write
\begin{align*}
\int_\rn \Delta_{h_j}^mf_j(z)\,\psi(z)\,dz
&= \sum_{i=0}^m (-1)^{m-i} \binom{m}{i}\,\int_\rn f_j(z+ih_j) \psi(z)\,dz \\
&= \sum_{i=0}^m (-1)^{m-i} \binom{m}{i}\,\int_\rn f_j(x) \psi(x-ih_j) \, dx \notag\\
&= \int_\rn f_j(x)\,\Delta_{\tilde h_j}^m\psi(x)\,dx.\notag
\end{align*}
By this identity and the H\"older inequality,
\begin{align*}
&\left|\int_\rn f(x)\,(\tilde\xi_0\cdot\nabla)^m\psi(x)\,dx-\int_\rn \frac{\Delta_{h_j}^mf_j(x)}{|h_j|^m}\, \psi(x)\,dx\right|\\
&\quad =\left|\int_\rn f(x)\,(\tilde\xi_0\cdot\nabla)^m\psi(x)\,dx-\int_\rn f_j(x)\,\frac{\Delta_{\tilde h_j}^m\psi(x)}{|h_j|^m}\,dx\right|\\
&\quad \le \left|\int_\rn \left[f(x)-f_j(x)\right]\,(\tilde\xi_0\cdot\nabla)^m\psi(x)\,dx\right|\\
&\qquad
+\left|\int_\rn f_j(x)\,\left[(\tilde\xi_0\cdot\nabla)^m\psi(x)-\frac{\Delta_{\tilde h_j}^m\psi(x)}{|h_j|^m}\right]
\,dx\right|\\
&\quad
\le  \|f-f_j\|_{L^p(\rn)} \left\|(\tilde\xi_0\cdot\nabla)^m\psi\right\|_{L^{p'}(\rn)}
+ \|f_j\|_{L^p(\rn)} \left\|
(\tilde\xi_0\cdot\nabla)^m\psi
-\frac{\Delta_{\tilde h_j}^m\psi}{|h_j|^m}\right\|_{L^{p'}(\rn)}.
\end{align*}
Since $f_j\to f$ in $L^p(\rn)$ and \eqref{eq-diff-LpNorm-psi-tilde}
holds, we have
\begin{align}\label{eq-Deltafj-f}
\int_\rn f(x)\,(\tilde\xi_0\cdot\nabla)^m\psi(x)\,dx
=\lim_{j\to\fz} \int_\rn \frac{\Delta_{h_j}^mf_j(x)}{|h_j|^m}\, \psi(x)\,dx.
\end{align}
From \eqref{eq-Deltafj-f} and  \eqref{eq-dual}, together with the relation
$(\tilde\xi_0\cdot\nabla)^m=(-1)^m(\xi_0\cdot\nabla)^m$, it follows that
\begin{align}\label{eq-liminf-ge}
\|(\xi_0\cdot\nabla)^mf\|_{L^p(\rn)}
&= \sup_{\gfz{\psi\in C_{\rm c}^\fz(\rn)}{\|\psi\|_{L^{p'}(\rn)}\le1}}
\left|\int_\rn f(x)\,(\tilde\xi_0\cdot\nabla)^m\psi(x)\,dx\right|\\
&=\sup_{\gfz{\psi\in C_{\rm c}^\fz(\rn)}{\|\psi\|_{L^{p'}(\rn)}\le1}}
\left|\liminf_{j\to\fz} \int_\rn \frac{\Delta_{h_j}^mf_j(x)}{|h_j|^m}\, \psi(x)\,dx\right|\notag\\
&\le \liminf_{j\to\fz} \frac{\|\Delta_{h_j}^mf_j\|_{L^p(\rn)}}{|h_j|^m}. \notag
\end{align}
Combining \eqref{eq-limsup-le} and \eqref{eq-liminf-ge} yields
\eqref{eq-fail-lim}  for $p\in (1,\fz)$ and $f\in W^{m,\,p}(\rn)$.

\medskip

{\bf Step 2: Proof of \eqref{eq-fail-lim} for $p=1$ and $f\in\BVm(\rn)$.}
Set $\nu_0:=(\xi_0\cdot\nabla)^mf$ and
$\nu_j:=\bigl(\frac{h_j}{|h_j|}\cdot\nabla\bigr)^mf$ for every $j\in\nn$.
In place of \eqref{eq-hj-xi0}, we now have
\[
\lim_{j\to\infty}\|\nu_j-\nu_0\|_{\mathcal{M}(\rn)}=0,
\]
which, together with \eqref{eq-Delta-BV-upper} in
Lemma~\ref{Lem-LXlem2829}(ii), yields
\begin{align}\label{eq-BVm-limsup-le}
\limsup_{j\to\fz}\,\frac{\|\Delta_{h_j}^mf\|_{L^1(\rn)}}{|h_j|^m}
&\le \limsup_{j\to\fz}
\left\|\left(\frac{h_j}{|h_j|}\cdot\nabla\right)^mf\right\|_{\mathcal{M}(\rn)}
= \limsup_{j\to\fz}\|\nu_j\|_{\mathcal{M}(\rn)}
= \|\nu_0\|_{\mathcal{M}(\rn)}.
\end{align}

In the $\BVm$ setting, the corresponding duality formula reads
\begin{align}\label{eq-dual-BV}
\|\nu_0\|_{\mathcal M(\rn)}
&= \sup_{\gfz{\psi\in C_{\rm c}^\fz(\rn)}{|\psi|\le1}}
\left|\int_\rn \psi\,d\nu_0\right|\\
&= \sup_{\gfz{\psi\in C_{\rm c}^\fz(\rn)}{|\psi|\le1}}
\left|(-1)^m\int_\rn f(x)\,(\xi_0\cdot\nabla)^m\psi(x)\,dx\right|. \notag
\end{align}
Moreover, for any $\psi\in C_{\rm c}^\fz(\rn)$, formula
\eqref{eq-diff-LpNorm-psi-tilde} remains valid with $p'=\infty$; hence,
\eqref{eq-Deltafj-f} holds. This gives
\begin{align*}
(-1)^m\int_\rn f(x)\,(\xi_0\cdot\nabla)^m\psi(x)\,dx
 =\lim_{j\to\fz}\int_\rn \frac{\Delta_{h_j}^mf_j(x)}{|h_j|^m}\,\psi(x)\,dx.
\end{align*}
Combining this with \eqref{eq-dual-BV}, we obtain
\begin{align}\label{eq-BVm-liminf-ge}
\|\nu_0\|_{\mathcal{M}(\rn)}
\le \liminf_{j\to\fz} \frac{\|\Delta_{h_j}^mf_j\|_{L^1(\rn)}}{|h_j|^m}.
\end{align}
Combining \eqref{eq-BVm-limsup-le} and \eqref{eq-BVm-liminf-ge} yields
\eqref{eq-fail-lim} for $p=1$ and $f\in\BVm(\rn)$.
\end{proof}

\begin{proof}[Proof of Proposition \ref{prop-BBM-Gamma-to-cond}]
Assume without loss of generality that $0\not\equiv f\in L^p(\rn)$; otherwise \eqref{eq-Gamma1-fj} holds trivially.
Then, by \eqref{eq-Gamma=>adm} and \eqref{eq-bdd-E(fj)}, we have
\begin{align}\label{eq-Gamma-bdd2}
M:= \limsup_{j\to\fz} \int_\rn \varphi_j(h) \min\lf\{|h|^{mq},\,1 \r\} dh<\fz.
\end{align}
For each $j\in\nn$, set
\begin{align*}
F_j:=\sup_{h\in\rn\setminus\{0\}}
\frac{\|\Delta_h^mf\|_{L^p(\rn)}-\|\Delta_h^mf_j\|_{L^p(\rn)}}{\min\{|h|^m,1\}}.
\end{align*}
Clearly, for every $h\in\rn$,
\begin{align*}
\|\Delta_h^mf\|_{L^p(\rn)}
\le\min\left\{|h|^m,\,1\right\}F_j+\|\Delta_h^mf_j\|_{L^p(\rn)}.
\end{align*}
Setting $r:=\max\{q,1\}$ and applying \eqref{eq-norm-minus-Euv}, we
obtain
\begin{align}\label{eq-limit0}
\lf(\mathscr{E}_{p,\,q}^{m,\,\varphi_j}(f) \r)^\frac{1}{r}
&\le  \lf( F_j^q\int_\rn \varphi_j(h) \min\lf\{|h|^{mq},\,1 \r\}\, dh \r)^\frac{1}{r}
+ \lf(\mathscr{E}_{p,\,q}^{m,\,\varphi_j}(f_j) \r)^\frac{1}{r}.
\end{align}
By \eqref{eq-Gamma-bdd2} and \eqref{eq-Ij-claim}, the first term on the
right-hand side of \eqref{eq-limit0} vanishes as $j\to\infty$, that is,
\begin{align*}
\lim_{j\to\infty}\left(F_j^q\int_\rn\varphi_j(h)\min\left\{|h|^{mq},\,1\right\}\,dh\right)^\frac{1}{r}=0.
\end{align*}
This proves \eqref{eq-Gamma1-fj} and hence completes the proof.
\end{proof}

\subsection{Proof of (iii) $\Rightarrow$ (ii) in Theorems \ref{thm-MS} and \ref{thm-BBM}}\label{ss-Gamma=>pw}

\begin{proof}[Proof of (iii) $\Rightarrow$ (ii) in Theorem \ref{thm-MS}]
Recall that $\mathscr{D}^0_{p,\,q}(f)= 2^\frac{q}{p}\,\|f\|_{L^p(\rn)}^q$. By (iii), we have
\begin{align*}
\Gamma\text{-}\lim_{j\to\fz} \mathscr{E}_{p,\,q}^{\varphi_j}
= \mathscr{D}^0_{p,\,q}
\quad\text{on}\,\ W^{1,\,p}(\rn).
\end{align*}
Let $f\in W^{1,\,p}(\rn)$. Applying $\mathbf{(\Gamma)_2}$ in Definition
\ref{def-Gamma} with the constant sequence $f_j=f$, we obtain
\begin{align*}
\liminf_{j\to\fz} \mathscr{E}_{p,\,q}^{\varphi_j}(f)
\ge \mathscr{D}^0_{p,\,q}(f).
\end{align*}
By $\mathbf{(\Gamma)_1}$, there exists a sequence $\{f_j\}_{j\in\nn}\subset L^p(\rn)$
such that $f_j\to f$ in $L^p(\rn)$ as $j\to\fz$ and
\begin{align*}
\limsup_{j\to\fz} \mathscr{E}_{p,\,q}^{\varphi_j}(f_j)
\le \mathscr{D}^0_{p,\,q}(f)=2^\frac{q}{p}\,\|f\|_{L^p(\rn)}^q<\fz.
\end{align*}
Combining this with Proposition \ref{prop-BBM-Gamma-to-cond} (with $m=1$ therein) gives
\begin{align*}
\limsup_{j\to\fz} \mathscr{E}_{p,\,q}^{\varphi_j}(f)
\le \limsup_{j\to\fz} \mathscr{E}_{p,\,q}^{\varphi_j}(f_j)
\le \mathscr{D}^0_{p,\,q}(f).
\end{align*}
This proves that $\{\mathscr{E}_{p,\,q}^{\varphi_j}(f)\}_{j\in\nn}$ converges to $\mathscr{D}^0_{p,\,q}(f)$
for every $f\in W^{1,\,p}(\rn)$. Thus, (ii) holds.
\end{proof}


\begin{proof}[Proof of (iii) $\Rightarrow$ (ii) in Theorem \ref{thm-BBM}]
Relabel the subsequence in (iii) by $j$. The proof is analogous to that of (iii) $\Rightarrow$ (ii) in Theorem
\ref{thm-MS} above. We only need to replace $W^{1,\,p}(\rn)$ by $W^{m,\,p}(\rn)$, $\mathscr{E}_{p,\,q}^{\varphi_j}$
by $\mathscr{E}_{p,\,q}^{m,\,\varphi_j}$ and $\mathscr{D}^0_{p,\,q}$ by $\mathscr{D}_{p,\,q}^{m,\,\mu}$, respectively;
the same argument then yields the desired conclusion.

The case $p=1$, with functions in $\BVm(\rn)$, can be treated similarly; the details are left to the reader.
\end{proof}

\section{Pointwise and strong $\Gamma$-convergence of the MS formula}\label{sec-MS}

The main aim of this section is to prove Theorem \ref{thm-MS}. By Remark \ref{rem-MS-BBM} and Section
\ref{ss-Gamma=>pw}, it suffices to show $\rm (ii) \ \Rightarrow\ (i)\ \Rightarrow\ (iv)$,
which are proved in Sections \ref{sec-MS1} and \ref{sec-MS2}, respectively.

\subsection{Proof of (ii) $\Rightarrow$ (i) in Theorem \ref{thm-MS}}\label{sec-MS1}

The following result yields the desired implication (ii) $\Rightarrow$ (i) in Theorem \ref{thm-MS}.

\begin{theorem}\label{Thm-mainresult-MS-varphi}
Let $p\in[1,\fz)$ and $q\in(0,\fz)$. For any $j\in\nn$, assume that $0\le \varphi_j\in
L_{\rm loc}^1(\rn \setminus\{0\})$, with $\varphi_j(0)$ possibly equal to $+\infty$.
Suppose that, for every $f\in C_{\rm c}^\fz(\rn)$,
\begin{align}\label{eq-MS-varphi}
\lim_{j\to\fz} \mathscr{E}_{p,\,q}^{\varphi_j}(f)=
 2^\frac{q}{p}\,\|f\|_{L^p(\rn)}^q.
\end{align}
Then, the sequence $\{\varphi_j\}_{j\in\nn}$ satisfies the $\MSq$ condition.
\end{theorem}

The proof of Theorem \ref{Thm-mainresult-MS-varphi} is based on the following two lemmas, which provide
the upper and lower estimates needed to verify the normalization condition $\NMSq$.

\begin{lemma}\label{Lem-MS-upper}
Let all hypotheses be as in Theorem \ref{Thm-mainresult-MS-varphi}. Then for any $\delta\in(0,\fz)$,
\begin{align}\label{eq-MS-upper}
\limsup_{j\to\fz} \int_{|h|\ge \delta} \varphi_j(h)\,dh\le1.
\end{align}
\end{lemma}

\begin{proof}
Let $0\not\equiv f\in C_{\rm c}^\fz(\rn)$ such that $\supp f\subset B(0,\delta/2)$. By Lemma
\ref{Lem-LXlem2829}(iv), when $|h|\ge\delta$, we have
\begin{align*}
\|\Delta_hf\|_{L^p(\rn)}=2^\frac{1}{p}\,\|f\|_{L^p(\rn)},
\end{align*}
which, together with \eqref{eq-MS-varphi}, implies
\begin{align*}
2^\frac{q}{p}\,\|f\|_{L^p(\rn)}^q
&\ge \limsup_{j\to\fz} \int_{|h|\ge\delta} \|\Delta_hf\|_{L^p(\rn)}^q\, \varphi_j(h)\,dh \\
&= \limsup_{j\to\fz} \int_{|h|\ge\delta} 2^\frac{q}{p}\,\|f\|_{L^p(\rn)}^q\, \varphi_j(h)\,dh.
\end{align*}
This gives \eqref{eq-MS-upper}.
\end{proof}

\begin{lemma}\label{Lem-MS-lower}
Let all hypotheses be as in Theorem \ref{Thm-mainresult-MS-varphi}. Then for any $\delta\in(0,\fz)$,
\begin{align}\label{eq-MS-lower}
\liminf_{j\to\fz} \int_{|h|\ge \delta} \varphi_j(h)\,dh\ge1.
\end{align}
\end{lemma}

\begin{proof}
Let $0\le f\in C_{\rm c}^\fz(\rn)$. For any $\lambda\in(0,\fz)$ and $x\in\rn$, set
\begin{align*}
f_\lambda(x):=f(\lambda x).
\end{align*}
Applying \eqref{eq-MS-varphi} to the function $f_\lambda$ gives
\begin{align*}
2^\frac{q}{p}\,\|f_\lambda\|_{L^p(\rn)}^q
= \liminf_{j\to\fz} \int_\rn \|\Delta_h f_\lambda\|_{L^p(\rn)}^q\, \varphi_j(h)\,dh.
\end{align*}
For any $h\in\rn$, observing that
\begin{align*}
\|f_\lambda\|_{L^p(\rn)}=\lambda^{-\frac{n}{p}}\,\|f\|_{L^p(\rn)}
\qquad\text{and}\qquad
\|\Delta_hf_\lambda\|_{L^p(\rn)}=\lambda^{-\frac{n}{p}}\,\|\Delta_{\lambda h}f\|_{L^p(\rn)},
\end{align*}
we therefore obtain
\begin{align}\label{eq-MS2}
2^\frac{q}{p}\,\|f\|_{L^p(\rn)}^q
&\le \liminf_{j\to\fz} \int_\rn \|\Delta_{\lambda h} f\|_{L^p(\rn)}^q\, \varphi_j(h)\,dh.
\end{align}

Fix $\delta\in(0,\fz)$. By $f\ge0$, we see that for any $\lambda\in(0,\infty)$ and $h\in\rn$,
\begin{align*}
\|\Delta_{\lambda h}f\|_{L^p(\rn)}^p
&= \int_\rn |f(x+{\lambda h})-f(x)|^p\,dx \\
&\le \int_\rn \max\lf\{|f(x+{\lambda h})|^p,\,|f(x)|^p \r\} dx \\
&\le \int_\rn \lf(|f(x+{\lambda h})|^p+|f(x)|^p \r) dx \\
&= 2\,\|f\|_{L^p(\rn)}^p,
\end{align*}
which in turn gives
\begin{align}\label{eq-MS-delta-fz}
\int_{|h|\ge\delta} \|\Delta_{\lambda h} f\|_{L^p(\rn)}^q\, \varphi_j(h)\,dh
\le 2^\frac{q}{p}\,\|f\|_{L^p(\rn)}^q \int_{|h|\ge\delta} \varphi_j(h)\,dh.
\end{align}

Since we have assumed \eqref{eq-MS-varphi}, it follows from Lemma \ref{Lem-Gamma-bdd}
(with every $f_j=f$ therein) that
\begin{align*}
M:=\limsup_{j\to\fz} \int_\rn \varphi_j(h) \min\lf\{|h|^q,\,1 \r\} dh<\fz.
\end{align*}
Note that Lemma \ref{Lem-LXlem2829}(ii) implies
\begin{align*}
\int_{|h|<\delta} \|\Delta_{\lambda h} f\|_{L^p(\rn)}^q\, \varphi_j(h)\,dh
\le \lambda^q\,\|f\|_{\mathring{W}^{1,p}(\rn)}^q \int_{|h|<\delta} |h|^q \varphi_j(h)\,dh.
\end{align*}
So, by taking the limsup as $j\to\infty$, we obtain
\begin{align*}
\limsup_{j\to\fz} \int_{|h|<\delta} \|\Delta_{\lambda h} f\|_{L^p(\rn)}^q\, \varphi_j(h)\,dh
&\le \lambda^q\,\|f\|_{\mathring{W}^{1,p}(\rn)}^q \limsup_{j\to\fz} \int_{|h|<\delta} |h|^q\,\varphi_j(h)\,dh \\
&\le M\max\{1,\,\delta^q\} \lambda^q\,\|f\|_{\mathring{W}^{1,p}(\rn)}^q.
\end{align*}
Consequently,
\begin{align}\label{eq-MS-0-delta}
\lim_{\lambda\to0_+}\limsup_{j\to\fz} \int_{|h|<\delta} \|\Delta_{\lambda h} f\|_{L^p(\rn)}^q\, \varphi_j(h)\,dh=0.
\end{align}

It is an elementary fact that for any two sequences $\{a_j\}_{j\in\nn}$ and $\{b_j\}_{j\in\nn}$ of positive numbers,
\begin{align*}
\inf_{j\ge n}\,(a_j+b_j)
\le \sup_{j\ge n} a_j+\inf_{j\ge n} b_j
\quad\text{for all}\ \, n\in\nn,
\end{align*}
which implies
\begin{align}\label{eq-limsupinf}
\liminf_{j\to\fz}\,(a_j+b_j)
\le \limsup_{j\to\fz} a_j+\liminf_{j\to\fz} b_j.
\end{align}
Applying this inequality together with \eqref{eq-MS2} and \eqref{eq-MS-delta-fz}, we obtain
\begin{align*}
2^\frac{q}{p}\,\|f\|_{L^p(\rn)}^q
&= \liminf_{j\to\fz} \lf(\int_{|h|< \delta} + \int_{|h|\ge \delta}\r) \|\Delta_{\lambda h} f\|_{L^p(\rn)}^q\, \varphi_j(h)\,dh\\
&\le \liminf_{j\to\fz} \lf( \int_{|h|< \delta} \|\Delta_{\lambda h} f\|_{L^p(\rn)}^q\, \varphi_j(h)\,dh
+ 2^\frac{q}{p}\,\|f\|_{L^p(\rn)}^q \int_{|h|\ge\delta} \varphi_j(h)\,dh\r)\\
&\le \limsup_{j\to\fz} \int_{|h|<\delta} \|\Delta_{\lambda h} f\|_{L^p(\rn)}^q\, \varphi_j(h)\,dh
+ 2^\frac{q}{p}\,\|f\|_{L^p(\rn)}^q \liminf_{j\to\fz} \int_{|h|\ge \delta} \varphi_j(h)\,dh.
\end{align*}
Further, letting $\lambda\to0_+$ and using \eqref{eq-MS-0-delta}, we arrive at
\begin{align*}
2^\frac{q}{p}\,\|f\|_{L^p(\rn)}^q
\le 2^\frac{q}{p}\,\|f\|_{L^p(\rn)}^q \liminf_{j\to\fz} \int_{|h|\ge\delta} \varphi_j(h)\,dh,
\end{align*}
which yields the desired estimate in \eqref{eq-MS-lower}.
\end{proof}

Applying Lemmas \ref{Lem-MS-upper} and \ref{Lem-MS-lower}, we now prove Theorem \ref{Thm-mainresult-MS-varphi}.

\begin{proof}[Proof of Theorem \ref{Thm-mainresult-MS-varphi}]
By \eqref{eq-MS-upper} and \eqref{eq-MS-lower}, we see that for any $\delta\in(0,\fz)$,
\begin{align}\label{eq-MS=1}
\lim_{j\to\fz} \int_{|h|\ge\delta} \varphi_j(h)\,dh=1.
\end{align}
In particular, $\NMSq$ follows immediately by letting $\delta=1$ in \eqref{eq-MS=1}.

It remains to show $\VMSq$. To this end, take $0\le f\in C_{\rm c}^\fz(\rn)$ with $\supp f\subset B(0,1/2)$.
By Lemma \ref{Lem-LXlem2829}(iii), there exist constants $\delta_0\in(0,1)$ and $c_0\in(0,\fz)$ such that
\begin{align*}
\|\Delta_h f\|_{L^p(\rn)}
\ge c_0\,|h|
\qquad\text{for all}\,\ |h|<\delta_0.
\end{align*}
On the one hand, applying Lemma \ref{Lem-LXlem2829}(iv) and the support condition on $f$ gives
\begin{align*}
\|\Delta_hf\|_{L^p(\rn)}=2^\frac{1}{p}\,\|f\|_{L^p(\rn)}
\quad \text{for all}\ \, |h|\ge 1,
\end{align*}
which, together with \eqref{eq-MS=1}, yields
\begin{align*}
\lim_{j\to\fz} \int_{|h|\ge 1}\|\Delta_hf\|_{L^p(\rn)}^q\,\varphi_j(h)\,dh
= 2^\frac{q}{p}\,\|f\|_{L^p(\rn)}^q \lim_{j\to\fz} \int_{|h|\ge 1} \varphi_j(h)\,dh
= 2^\frac{q}{p}\,\|f\|_{L^p(\rn)}^q.
\end{align*}
On the other hand, from \eqref{eq-MS-delta-fz} and \eqref{eq-MS=1}, it follows that
\begin{align*}
\limsup_{j\to\fz} \int_{|h|\ge \delta_0}\|\Delta_hf\|_{L^p(\rn)}^q\,\varphi_j(h)\,dh
&\le 2^\frac{q}{p}\,\|f\|_{L^p(\rn)}^q \lim_{j\to\fz} \int_{|h|\ge \delta_0} \varphi_j(h)\,dh\\
&= 2^\frac{q}{p}\,\|f\|_{L^p(\rn)}^q.
\end{align*}
Therefore, the last two estimates give
\begin{align*}
\lim_{j\to\fz} \int_{|h|\ge \delta_0}\|\Delta_hf\|_{L^p(\rn)}^q\,\varphi_j(h)\,dh
= 2^\frac{q}{p}\,\|f\|_{L^p(\rn)}^q.
\end{align*}
Combining the latter identity with \eqref{eq-MS-varphi}, we obtain
\begin{align*}
\lim_{j\to\fz} \int_{|h|< \delta_0}\|\Delta_hf\|_{L^p(\rn)}^q\,\varphi_j(h)\,dh =0.
\end{align*}
Furthermore, since
\begin{align*}
\int_{|h|<\delta_0} \|\Delta_{h} f\|_{L^p(\rn)}^q\, \varphi_j(h)\,dh
\ge \int_{|h|<\delta_0} (c_0|h|)^q\, \varphi_j(h)\,dh,
\end{align*}
it follows that
\begin{align}\label{eq-MS-claim}
\lim_{j\to\fz} \int_{|h|<\delta_0} |h|^q\, \varphi_j(h)\,dh=0.
\end{align}

Now, let $\delta\in(\delta_0,\fz)$. Using \eqref{eq-MS=1} again, we have
\begin{align*}
\lim_{j\to\fz} \int_{\delta_0\le |h|<\delta} \varphi_j(h) \min\lf\{|h|^q,\,1 \r\} dh
&\le \lim_{j\to\fz} \int_{\delta_0\le |h|<\delta} \varphi_j(h)\, dh\\
&= \lim_{j\to\fz} \int_{|h|\ge \delta_0} \varphi_j(h)\, dh
- \lim_{j\to\fz} \int_{|h|\ge \delta} \varphi_j(h)\, dh \\
&=0.
\end{align*}
Combining this formula with \eqref{eq-MS-claim}, we deduce that
\begin{align*}
&\lim_{j\to\fz} \int_{|h|<\delta} \varphi_j(h) \min\lf\{|h|^q,\,1 \r\} dh \\
&\quad = \lim_{j\to\fz} \lf(\int_{|h|<\delta_0} |h|^q\,\varphi_j(h)\,dh
+ \int_{\delta_0\le|h|<\delta} \varphi_j(h) \min\lf\{|h|^q,\,1 \r\} dh \r) \\
&\quad =0.
\end{align*}
This proves $\VMSq$ and thereby completes the proof of Theorem \ref{Thm-mainresult-MS-varphi}.
\end{proof}

\subsection{Proof of (i) $\Rightarrow$ (iv) in Theorem \ref{thm-MS}}\label{sec-MS2}

The implication $\rm (i)\Rightarrow(iv)$ in Theorem \ref{thm-MS} follows directly from the following result.

\begin{theorem}\label{Thm-MS-Gamma-converge}
Let $p\in[1,\fz)$, $q\in(0,\fz)$ and
$\{\varphi_j\}_{j\in\nn}$
satisfy $\MSq$.
Then, for any $f\in W^{1,\,p}(\rn)$,
\begin{align}\label{eq-MS-Gamma-conv1}
\lim_{j\to\fz} \int_\rn \lf\|\Delta_hf\r\|_{L^p(\rn)}^q \varphi_j(h)\,dh
= 2^\frac{q}{p}\,\|f\|_{L^p(\rn)}^q.
\end{align}
Moreover, if $f\in W^{1,\,p}(\rn)$ and $\{f_j\}_{j\in\nn}\subset L^p(\rn)$ such that $f_j\to f$ in
$L^p(\rn)$ as $j\to\fz$, then
\begin{align}\label{eq-MS-Gamma-conv2}
\liminf_{j\to\fz} \int_\rn \lf\|\Delta_hf_j \r\|_{L^p(\rn)}^q \varphi_j(h)\,dh
\ge 2^\frac{q}{p}\,\|f\|_{L^p(\rn)}^q.
\end{align}
\end{theorem}

\begin{proof}
The proof is divided into three steps.

\medskip

{\bf Step 1: preparation.}
Let $\delta\in(0,\fz)$ and $f\in W^{1,\,p}(\rn)$. By Lemma \ref{Lem-LXlem2829}(i) and (ii), the inequality
 \begin{align*}
\|\Delta_hf\|_{L^p(\rn)}
\le \min\lf\{|h|\,\|f\|_{\mathring{W}^{1,\,p}(\rn)},
\ 2\,\|f\|_{L^p(\rn)} \r\}
\end{align*}
holds for all $h\in\rn$,
which further implies
\begin{align*}
\int_{|h|<\delta} \lf\|\Delta_hf \r\|_{L^p(\rn)}^q \varphi_j(h)\,dh
&\le \int_{|h|<\delta} \lf(\min\lf\{|h|\,\|f\|_{\mathring{W}^{1,\,p}(\rn)},\
2\,\|f\|_{L^p(\rn)} \r\}\r)^q \varphi_j(h)\,dh.
\end{align*}
From this and $\VMSq$, it follows that for every $\delta\in(0,\fz)$,
\begin{align}\label{eq-MS-claim2-eq1}
\lim_{j\to\fz} \int_{|h|<\delta} \lf\|\Delta_hf \r\|_{L^p(\rn)}^q \varphi_j(h)\,dh=0.
\end{align}

\medskip

{\bf Step 2: proof of \eqref{eq-MS-Gamma-conv1}.}
Assume without loss of generality that $\|f\|_{L^p(\rn)}\neq 0$; otherwise both sides of \eqref{eq-MS-Gamma-conv1} are zero.
By Lemma \ref{Lem-LXlem2829}(v), for any $\uc\in(0,1)$, there exists $M\in(1,\fz)$ such that for all $|h|\ge M$,
\begin{align*}
\lf(2^{\frac q p}-\uc\r)\|f\|_{L^p(\rn)}^q
\le \|\Delta_hf\|_{L^p(\rn)}^q
\le \lf(2^{\frac q p}+\uc\r)\|f\|_{L^p(\rn)}^q.
\end{align*}
Moreover, by $\NMSq$ and $\VMSq$, we have
\begin{align*}
\lim_{j\to\fz}\int_{1\le |h|<M} \varphi_j(h)\,dh
= \lim_{j\to\fz} \int_{1\le |h|<M} \varphi_j(h) \min\lf\{|h|^q,\, 1 \r\} dh=0
\end{align*}
and, hence,
\begin{align*}
\lim_{j\to\fz}\int_{|h|\ge M} \varphi_j(h)\,dh
= \lim_{j\to\fz}\lf(\int_{|h|\ge 1} \varphi_j(h)\,dh
- \int_{1\le |h|<M} \varphi_j(h)\,dh\r)=1.
\end{align*}
Using these estimates, we obtain
\begin{align}\label{eq-MS-inf}
\liminf_{j\to\fz} \int_\rn \|\Delta_hf\|_{L^p(\rn)}^q\,\varphi_j(h)\,dh
&\ge \liminf_{j\to\fz} \int_{|h|\ge M} \|\Delta_hf\|_{L^p(\rn)}^q\, \varphi_j(h)\,dh \\
&\ge \lf(2^\frac{q}{p}-\uc\r) \|f\|_{L^p(\rn)}^q \liminf_{j\to\fz} \int_{|h|\ge M} \varphi_j(h)\,dh \notag \\
&= \lf(2^\frac{q}{p}-\uc \r) \|f\|_{L^p(\rn)}^q . \notag
\end{align}
Similarly, we have
\begin{align*}
\limsup_{j\to\fz} \int_{|h|\ge M} \|\Delta_hf\|_{L^p(\rn)}^q\, \varphi_j(h)\,dh
&\le \lf(2^\frac{q}{p}+\uc \r) \|f\|_{L^p(\rn)}^q \limsup_{j\to\fz}\int_{|h|\ge M} \varphi_j(h)\,dh \\
&= \lf(2^\frac{q}{p}+\uc \r) \|f\|_{L^p(\rn)}^q,
\end{align*}
which, combined with \eqref{eq-MS-claim2-eq1} (with $\delta=M$), implies
\begin{align}\label{eq-MS-sup}
&\limsup_{j\to\fz} \int_\rn \|\Delta_hf\|_{L^p(\rn)}^q\,\varphi_j(h)\,dh\\
&\quad\le \limsup_{j\to\fz} \int_{|h|< M} \|\Delta_hf\|_{L^p(\rn)}^q\, \varphi_j(h)\,dh
+ \limsup_{j\to\fz} \int_{|h|\ge M} \|\Delta_hf\|_{L^p(\rn)}^q\, \varphi_j(h)\,dh \notag\\
&\quad\le \lf(2^\frac{q}{p}+\uc \r) \|f\|_{L^p(\rn)}^q. \notag
\end{align}
Letting $\uc\to0_+$ in \eqref{eq-MS-inf} and \eqref{eq-MS-sup} then yields exactly \eqref{eq-MS-Gamma-conv1}.

\medskip

{\bf Step 3: proof of \eqref{eq-MS-Gamma-conv2}.}
For the right-hand side of \eqref{eq-MS-Gamma-conv2}, since $f\in W^{1,\,p}(\rn)$, it follows from \eqref{eq-MS-Gamma-conv1} and \eqref{eq-MS-claim2-eq1} that
\begin{align*}
2^\frac{q}{p}\,\|f\|_{L^p(\rn)}^q
= \lim_{j\to\fz} \int_\rn \|\Delta_h f\|_{L^p(\rn)}^q\,\varphi_j(h)\,dh
= \lim_{j\to\fz} \int_{|h|\ge 1} \|\Delta_h f\|_{L^p(\rn)}^q\,\varphi_j(h)\,dh.
\end{align*}
For the left-hand side of \eqref{eq-MS-Gamma-conv2}, one has
\begin{align*}
\liminf_{j\to\fz} \int_\rn \|\Delta_h f_j\|_{L^p(\rn)}^q\,\varphi_j(h)\,dh
\ge \liminf_{j\to\fz} \int_{|h|\ge 1} \|\Delta_h f_j\|_{L^p(\rn)}^q\, \varphi_j(h)\,dh.
\end{align*}
Therefore, to establish \eqref{eq-MS-Gamma-conv2}, it suffices to prove
\begin{align}\label{eq-MS-lim=}
\lim_{j\to\fz} \int_{|h|\ge1} \|\Delta_h f_j\|_{L^p(\rn)}^q\,\varphi_j(h)\,dh
= \lim_{j\to\fz} \int_{|h|\ge1} \|\Delta_h f\|_{L^p(\rn)}^q\, \varphi_j(h)\,dh.
\end{align}

We now prove \eqref{eq-MS-lim=}. By Lemma \ref{Lem-LXlem2829}(i), for all $j\in\nn$ and $h\in\rn$, we have
\begin{align*}
\|\Delta_hf_j-\Delta_h f\|_{L^p(\rn)}
\le 2\,\|f_j-f\|_{L^p(\rn)}.
\end{align*}
Since $f_j\to f$ in $L^p(\rn)$ as $j\to\fz$, for any $\uc\in(0,1)$, there exists $J\in\nn$ such that for all $j\ge J$ and $h\in\rn$,
\begin{align*}
\|\Delta_hf_j-\Delta_h f\|_{L^p(\rn)}<\uc.
\end{align*}
From this and \eqref{eq-norm-minus}, with $r=\max\{q,1\}$, it follows that
\begin{align*}
&\lf|\lf( \int_{|h|\ge 1} \|\Delta_hf_j\|_{L^p(\rn)}^q\, \varphi_j(h)\,dh\r)^\frac 1{r}
-\lf( \int_{|h|\ge 1} \|\Delta_hf\|_{L^p(\rn)}^q\, \varphi_j(h)\,dh\r)^\frac 1{r}\r|\\
& \quad\le \lf(\int_{|h|\ge 1} \|\Delta_hf_j-\Delta_h f\|_{L^p(\rn)}^q\, \varphi_j(h)\,dh \r)^\frac 1{r} \notag\\
&\quad \le \uc^\frac q{r} \lf( \int_{|h|\ge 1} \varphi_j(h)\,dh\r)^\frac{1}{r}.
\end{align*}
Then, letting $j\to\infty$ and using $\NMSq$ yields
\begin{align*}
&\limsup_{j\to\infty}\lf|\lf( \int_{|h|\ge 1} \|\Delta_hf_j\|_{L^p(\rn)}^q\,\varphi_j(h)\,dh \r)^\frac{1}{r}
-\lf( \int_{|h|\ge 1} \|\Delta_hf\|_{L^p(\rn)}^q\, \varphi_j(h)\,dh\r)^\frac{1}{r} \r| \\
&\quad \le \uc^\frac q{r} \limsup_{j\to\infty}\lf( \int_{|h|\ge 1} \varphi_j(h)\,dh\r)^\frac 1{r}\\
&\quad= \uc^\frac q{r}.
\end{align*}
Due to the arbitrariness of $\uc$, we conclude that
\begin{align*}
\limsup_{j\to\infty}\lf|\lf( \int_{|h|\ge 1} \|\Delta_hf_j\|_{L^p(\rn)}^q\,\varphi_j(h)\,dh\r)^\frac 1{r}
-\lf( \int_{|h|\ge 1} \|\Delta_hf\|_{L^p(\rn)}^q\, \varphi_j(h)\,dh\r)^\frac 1{r}\r|=0.
\end{align*}
This proves \eqref{eq-MS-lim=}. Consequently, we obtain \eqref{eq-MS-Gamma-conv2}.
\end{proof}

\section{Pointwise and strong $\Gamma$-convergence of the BBM formula}\label{sec-BBM}

This section is devoted to the proofs of Theorems \ref{thm-BBM} and \ref{thm-BBM-radial}. By Remark \ref{rem-MS-BBM}
and the discussions in Section \ref{ss-Gamma=>pw}, to prove Theorem \ref{thm-BBM}, it suffices to show
$\mathrm{(ii)}\ \Rightarrow\ \mathrm{(i)}\ \Rightarrow\ \mathrm{(iv)}.$ These two implications are proved in Sections
\ref{sec-BBM1} and \ref{sec-BBM2}, respectively. Further, Theorem \ref{thm-BBM-radial} is proved in Section \ref{sec-BBM4}.

\subsection{Proof of (ii) $\Rightarrow$ (i) in Theorem \ref{thm-BBM}}\label{sec-BBM1}

The following theorem establishes the desired implication.

\begin{theorem}\label{Thm-BBM-cond}
Let $p\in[1,\infty)$, $q\in(0,\infty)$ and $m\in\nn$. For any $j\in\nn$, let $0\le \varphi_j\in L_{\rm loc}^1(\rn\setminus\{0\})$,
with $\varphi_j$ possibly equal to  $+\infty$. Suppose that there exists a Borel probability measure $\mu$ on $\msn$ such that
\begin{align}\label{eq-BBM-varphi}
\lim_{j\to\infty}
\mathscr{E}_{p,\,q}^{m,\,\varphi_j}(f)
=\mathscr{D}_{p,\,q}^{m,\,\mu}(f)
\end{align}
for every $f\in C_{\mathrm c}^\infty(\rn)$. Then $\{\varphi_j\}_{j\in\nn}$ satisfies the $\BBMmq$ condition.
\end{theorem}

Since $f\in C_{\rm c}^\infty(\rn)$, when $p=1$, the right-hand side of \eqref{eq-BBM-varphi} coincides with
$\mathscr{D}_q^{m,\,\mu}(f)$. The proof of Theorem \ref{Thm-BBM-cond} is based on Lemmas \ref{Lem-DeltaNabla},
\ref{Lem-upper}, and \ref{Lem-lower}.

\begin{lemma}\label{Lem-DeltaNabla}
Let all the hypotheses be as in Theorem \ref{Thm-BBM-cond}. Then, for any $f\in C_{\rm c}^\fz(\rn)$,
\begin{align}\label{eq1-Delta-Nabla}
\lim_{\delta\to 0_+} \limsup_{j\to\infty} \int_{|h|<\delta} \|\Delta_h^mf-(h\cdot\nabla)^mf\|_{L^p(\rn)}^q\,\varphi_j(h)\,dh=0
\end{align}
and
\begin{align}\label{eq1-BBM-small-delta}
\lim_{\delta\to 0_+} \mathop{\operatorname{lim\,L}}_{j\to\infty} \int_{|h|<\delta} \|\Delta_h^mf\|_{L^p(\rn)}^q\,\varphi_j(h)\,dh
= \lim_{\delta\to 0_+} \mathop{\operatorname{lim\,L}}_{j\to\infty} \int_{|h|<\delta} \|(h\cdot\nabla)^mf\|_{L^p(\rn)}^q\,\varphi_j(h)\,dh.
\end{align}
Here and in what follows, the symbol $\operatorname{lim\,L}$ stands for either
$\liminf$ or $\limsup$, consistently within each formula.
\end{lemma}

\begin{proof}
Let $f\in C_{\rm c}^\infty(\rn)$ and $r:=\max\{q,\,1\}$. If we have proved \eqref{eq1-Delta-Nabla}, then by \eqref{eq-norm-minus}, we obtain
\begin{align*}
&\lim_{\delta\to 0_+} \mathop{\operatorname{lim\,L}}_{j\to\infty} \lf|\lf( \int_{|h|<\delta} \|\Delta_h^mf\|_{L^p(\rn)}^q\,\varphi_j(h)\,dh \r)^\frac 1{r}
- \lf( \int_{|h|<\delta} \|(h\cdot\nabla)^mf\|_{L^p(\rn)}^q\,\varphi_j(h)\,dh \r)^\frac{1}{r}\r| \\
& \quad\le \lim_{\delta\to 0_+} \limsup_{j\to\infty} \lf(\int_{|h|<\delta} \|\Delta_h^m f-(h\cdot\nabla)^mf\|_{L^p(\rn)}^q\,\varphi_j(h)\,dh \r)^\frac{1}{r} \\
& \quad=0,
\end{align*}
which establishes \eqref{eq1-BBM-small-delta}.

It remains to prove \eqref{eq1-Delta-Nabla}.
By $f\in C_{\rm c}^\infty(\rn)$ and Lemma \ref{lem-diffLpNorm}, there exists a constant
$C\in(0,\infty)$, depending only on $m$, $n$ and $f$, such that for all $h\in\rn$,
\begin{align*}
\|\Delta_h^mf-(h\cdot\nabla)^m f\|_{L^p(\rn)}
\le C\,|h|^{m+1}.
\end{align*}
Consequently, for any $\delta\in(0,1)$,
\begin{align*}
\int_{|h|<\delta} \|\Delta_h^m f-(h\cdot\nabla)^mf\|_{L^p(\rn)}^q\,\varphi_j(h)\,dh
&\le (C\delta)^q \int_{|h|<\delta} |h|^{mq}\,\varphi_j(h)\,dh\\
&= (C\delta)^q \int_{|h|<\delta} \varphi_j(h) \min\lf\{|h|^{mq},\,1 \r\} dh.
\end{align*}
Moreover, applying \eqref{eq-BBM-varphi} and Lemma \ref{Lem-Gamma-bdd} (with each $f_j$ therein chosen to be $f$) yields
\begin{align*}
M:=\limsup_{j\to\infty}\int_\rn \varphi_j(h) \min\lf\{|h|^{mq},\,1 \r\} dh<\infty.
\end{align*}
Thus,
\begin{align*}
\limsup_{j\to\infty} \int_{|h|<\delta} \|\Delta_h^m f-(h\cdot\nabla)^mf\|_{L^p(\rn)}^q\,\varphi_j(h)\,dh
\le (C\delta)^qM.
\end{align*}
Letting $\delta\to0_+$ yields the
desired limit in \eqref{eq1-Delta-Nabla}.
\end{proof}

\begin{lemma}\label{Lem-upper}
Let all the hypotheses be as in Theorem \ref{Thm-BBM-cond}.
Then, for any radial function $f\in C_{\rm c}^\infty(\rn)$,
\begin{align}\label{eq-BBM-smalldelta}
&\lim_{\delta\to 0_+} \mathop{\operatorname{lim\,L}}_{j\to\infty} \int_{|h|<\delta} \|\Delta_h^mf\|_{L^p(\rn)}^q\,\varphi_j(h)\,dh\\
&\quad = \mathscr{D}_{p,\,q}^{m,\,\mu}(f) \lf(\lim_{\delta\to 0_+} \mathop{\operatorname{lim\,L}}_{j\to\infty} \int_{|h|<\delta} |h|^{mq}\,\varphi_j(h)\,dh \r). \notag
\end{align}
Consequently,
\begin{align}\label{eq-upper}
\limsup_{\delta\to0_+} \lf(\limsup_{j\to\infty} \int_{|h|<\delta} |h|^{mq}\,\varphi_j(h)\,dh \r)\le1.
\end{align}
\end{lemma}

\begin{proof}
Let $f\in C_{\rm c}^\infty(\rn)$ be radial.
For any $\xi\in\msn$, choose an orthogonal
matrix $\mathcal{O}_\xi$ such that
$\mathcal{O}_\xi\mathbf e_1=\xi$, where $\mathbf e_1=(1,0,\ldots,0)\in\rn$.
By the chain rule, we have
\begin{align*}
\frac{\partial}{\partial x_1}
[f(\mathcal{O}_\xi x)]
=(\xi\cdot\nabla)f(\mathcal{O}_\xi x).
\end{align*}
Iterating this identity, we obtain
\begin{align*}
\frac{\partial^m[f(\mathcal{O}_\xi\cdot)]}{\partial x_1^m}(x)
=(\xi\cdot\nabla)^mf(\mathcal{O}_\xi x).
\end{align*}
By the radiality of $f$ and the invariance of Lebesgue measure under orthogonal transformations, we obtain
\begin{align*}
\|(\xi\cdot\nabla)^mf\|_{L^p(\rn)}
&=\|(\xi\cdot\nabla)^mf(\mathcal{O}_\xi\cdot)\|_{L^p(\rn)}\\
&=\left\|\frac{\partial^m[f(\mathcal{O}_\xi\cdot)]}
{\partial x_1^m}\right\|_{L^p(\rn)}=\left\|\frac{\partial^mf}{\partial x_1^m}\right\|_{L^p(\rn)}.
\end{align*}
In particular, the above quantity is independent of $\xi\in\msn$.
Since $\mu(\msn)=1$, it follows that
\begin{align}\label{eq-f-radial}
\mathscr{D}_{p,\,q}^{m,\,\mu}(f)
= \lf\|\frac{\partial^mf}{\partial x_1^m} \r\|_{L^p(\rn)}^q
= \|(\xi\cdot\nabla)^mf\|_{L^p(\rn)}^q.
\end{align}
Using \eqref{eq-f-radial}, for any $\delta\in(0,1)$, we have
\begin{align}\label{eq-sphe-coord}
\int_{|h|<\delta} \|(h\cdot\nabla)^mf\|_{L^p(\rn)}^q\,\varphi_j(h)\,dh
&= \int_{|h|<\delta} \lf\|\lf(\frac{h}{|h|}\cdot\nabla \r)^mf \r\|_{L^p(\rn)}^q |h|^{mq}\,\varphi_j(h)\,dh \\ \notag
&= \mathscr{D}_{p,\,q}^{m,\,\mu}(f)\,\lf(\int_{|h|<\delta} |h|^{mq}\,\varphi_j(h)\,dh \r).
\end{align}

Since $f\in C_{\rm c}^\fz(\rn)$, it follows that $f$ satisfies \eqref{eq1-BBM-small-delta}.
Substituting \eqref{eq-sphe-coord} into \eqref{eq1-BBM-small-delta} yields
\begin{align*}
&\lim_{\delta\to 0_+} \mathop{\operatorname{lim\,L}}_{j\to\infty} \lf( \int_{|h|<\delta} \|\Delta_h^mf\|_{L^p(\rn)}^q\,\varphi_j(h)\,dh\r)\\
&\quad= \lim_{\delta\to 0_+} \mathop{\operatorname{lim\,L}}_{j\to\infty} \lf( \int_{|h|<\delta} \|(h\cdot\nabla)^mf\|_{L^p(\rn)}^q\,\varphi_j(h)\,dh\r) \\
&\quad= \mathscr{D}_{p,\,q}^{m,\,\mu}(f)\,\lf(\lim_{\delta\to 0_+} \mathop{\operatorname{lim\,L}}_{j\to\infty} \int_{|h|<\delta} |h|^{mq}\,\varphi_j(h)\,dh \r).
\end{align*}
This proves \eqref{eq-BBM-smalldelta}.

For any radial $f\in C_{\rm c}^\infty(\rn)$ with $\mathscr{D}_{p,\,q}^{m,\,\mu}(f)\neq0$,
taking $\operatorname{lim\,L}=\limsup$ in
 \eqref{eq-BBM-smalldelta}
and using \eqref{eq-BBM-varphi}, we obtain
\begin{align*}
\mathscr{D}_{p,\,q}^{m,\,\mu}(f) \lf(\lim_{\delta\to 0_+} \limsup_{j\to\infty} \int_{|h|<\delta} |h|^{mq}\,\varphi_j(h)\,dh \r)
\le \limsup_{j\to\infty}\ \mathscr{E}_{p,\,q}^{m,\,\varphi_j}(f)
= \mathscr{D}_{p,\,q}^{m,\,\mu}(f).
\end{align*}
Dividing by the positive quantity $\mathscr{D}_{p,\,q}^{m,\,\mu}(f)$ gives \eqref{eq-upper}. This completes the proof.
\end{proof}

\begin{lemma}\label{Lem-lower}
Let all the hypotheses be as in Theorem \ref{Thm-BBM-cond}. Then, for any $\delta\in(0,\infty)$,
\begin{align}\label{eq-lower}
\liminf_{j\to\fz} \int_{|h|<\delta} |h|^{mq}\,\varphi_j(h)\,dh\ge1.
\end{align}
\end{lemma}

\begin{proof}
It suffices to prove \eqref{eq-lower} for $\delta\in(0,1)$. Choose a radial function $f\in C_{\rm c}^\fz(\rn)$ such that
$\supp f\subset B(0,1)$ and $\mathscr{D}_{p,\,q}^{m,\,\mu}(f)\neq0$. For any $\lambda\in(0,\fz)$ and $x\in\rn$, define
\begin{align*}
f_\lambda(x):=f(\lambda x).
\end{align*}
Applying \eqref{eq-BBM-varphi} to $f_\lambda$ gives
\begin{align}\label{eq-flz}
\mathscr{D}_{p,\,q}^{m,\,\mu}(f_\lambda)
= \lim_{j\to\fz} \mathscr{E}_{p,\,q}^{m,\,\varphi_j}(f_\lambda).
\end{align}
Observe that
\begin{align*}
\|(\xi\cdot\nabla)^mf_\lambda\|_{L^p(\rn)}
= \lambda^m\,\|(\xi\cdot\nabla)^mf(\lambda \,\cdot)\|_{L^p(\rn)}
= \lambda^{m-\frac{n}{p}}\,\|(\xi\cdot\nabla)^mf\|_{L^p(\rn)}
\end{align*}
and
\begin{align*}
\|\Delta_h^mf_\lambda\|_{L^p(\rn)}
= \|(\Delta_{\lambda h}^mf)(\lambda\,\cdot)\|_{L^p(\rn)}
= \lambda^{-\frac{n}{p}}\,\|\Delta_{\lambda h}^mf\|_{L^p(\rn)}.
\end{align*}
Combining these last two identities with \eqref{eq-flz}, we then obtain
\begin{align*}
\lambda^{(m-\frac{n}{p})q}\,\mathscr{D}_{p,\,q}^{m,\,\mu}(f)
= \lim_{j\to\fz} \mathscr{E}_{p,\,q}^{m,\,\varphi_j}(f_\lambda)
= \lambda^{-\frac{nq}{p}} \lim_{j\to\fz} \int_\rn \|\Delta_{\lambda h}^mf\|_{L^p(\rn)}^q\,\varphi_j(h)\,dh,
\end{align*}
or equivalently,
\begin{align*}
\mathscr{D}_{p,\,q}^{m,\,\mu}(f)
= \lambda^{-mq} \lim_{j\to\fz} \int_\rn \|\Delta_{\lambda h}^mf\|_{L^p(\rn)}^q\,\varphi_j(h)\,dh.
\end{align*}
Note that, for any nonnegative sequences $\{a_j\}_{j\in\nn}$ and $\{b_j\}_{j\in\nn}$,
there is the basic inequality (see \eqref{eq-limsupinf})
\begin{align*}
\liminf_{j\to\fz}\,(a_j+b_j)
\le \limsup_{j\to\fz} a_j+\liminf_{j\to\fz} b_j.
\end{align*}
This further implies that for any $\delta\in(0,1)$ and $\lambda\in(0,\infty)$,
\begin{align}\label{eq-flz1}
\mathscr{D}_{p,\,q}^m(f)
&\le \lambda^{-mq}\,\liminf_{j\to\fz} \int_{|h|<\delta} \|\Delta_{\lambda h}^mf\|_{L^p(\rn)}^q\,\varphi_j(h)\,dh \\
&\quad+ \lambda^{-mq}\,\limsup_{j\to\fz} \int_{|h|\ge\delta} \|\Delta_{\lambda h}^mf\|_{L^p(\rn)}^q\,\varphi_j(h)\,dh. \notag
\end{align}

Let us deal with the first term on the right-hand side of \eqref{eq-flz1}. By Lemma \ref{Lem-LXlem2829}(ii), we have
\begin{align*}
\|\Delta_{\lambda h}^mf\|_{L^p(\rn)}
\le \|(\lambda h\cdot\nabla)^mf\|_{L^p(\rn)}
= \lambda^m\,\|(h\cdot\nabla)^mf\|_{L^p(\rn)}.
\end{align*}
Using this estimate, together with \eqref{eq-sphe-coord}, we obtain
\begin{align*}
\lambda^{-mq} \int_{|h|<\delta} \|\Delta_{\lambda h}^mf\|_{L^p(\rn)}^q\,\varphi_j(h)\,dh
&\le \int_{|h|<\delta} \|(h\cdot\nabla)^mf\|_{L^p(\rn)}^q\,\varphi_j(h)\,dh \\
&= \mathscr{D}_{p,\,q}^{m,\,\mu}(f) \int_{|h|<\delta} |h|^{mq}\,\varphi_j(h)\,dh,
\end{align*}
which gives
\begin{align*}
\lambda^{-mq}\,\liminf_{j\to\fz} \int_{|h|<\delta} \|\Delta_{\lambda h}^mf\|_{L^p(\rn)}^q\,\varphi_j(h)\,dh
&\le \mathscr{D}_{p,\,q}^{m,\,\mu}(f) \lf(\liminf_{j\to\fz} \int_{|h|<\delta} |h|^{mq}\,\varphi_j(h)\,dh \r).
\end{align*}
Substituting this into \eqref{eq-flz1} gives that for all $\delta\in(0,1)$ and $\lambda\in(0,\infty)$,
\begin{align}\label{eq-flz2}
\mathscr{D}_{p,\,q}^{m,\,\mu}(f)
&\le \mathscr{D}_{p,\,q}^{m,\,\mu}(f) \lf(\liminf_{j\to\fz} \int_{|h|<\delta} |h|^{mq}\,\varphi_j(h)\,dh \r) \\
&\quad+ \lambda^{-mq}\,\limsup_{j\to\fz} \int_{|h|\ge\delta} \|\Delta_{\lambda h}^mf\|_{L^p(\rn)}^q\,\varphi_j(h)\,dh. \notag
\end{align}

In \eqref{eq-flz2}, we restrict to the case $\lambda:=\delta^{-2}$. Then, using the bound
$\|\Delta_{\delta^{-2} h}^mf \|_{L^p(\rn)}\le 2^m\,\|f\|_{L^p(\rn)}$, we have
\begin{align*}
\delta^{2mq}\,\limsup_{j\to\fz} \int_{|h|\ge\delta} \lf\|\Delta_{\delta^{-2}h}^mf \r\|_{L^p(\rn)}^q \varphi_j(h)\,dh
\le 2^{mq}\,\delta^{2mq}\,\|f\|_{L^p(\rn)}^q\,\limsup_{j\to\fz} \int_{|h|\ge\delta} \varphi_j(h)\,dh.
\end{align*}
Applying \eqref{eq-BBM-varphi} and Lemma \ref{Lem-Gamma-bdd} (with each $f_j$ therein chosen to be $f$), we see that
\begin{align*}
M:=\limsup_{j\to\fz} \int_\rn \varphi_j(h) \min\lf\{|h|^{mq},\,1 \r\} dh<\infty,
\end{align*}
which leads to
\begin{align*}
\limsup_{j\to\fz} \int_{|h|\ge\delta} \varphi_j(h)\,dh
\le \limsup_{j\to\fz} \lf(\delta^{-mq} \int_{|h|\ge\delta} \varphi_j(h) \min\lf\{|h|^{mq},\,1 \r\} dh \r)
\le M\delta^{-mq}.
\end{align*}
Therefore,
\begin{align}\label{eq-flz3}
&\lim_{\delta\to 0_+} \lf(\delta^{2mq}\,\limsup_{j\to\infty} \int_{|h|\ge\delta} \lf\|\Delta_{\delta^{-2} h}^mf \r\|_{L^p(\rn)}^q \varphi_j(h)\,dh \r)
\le M \lim_{\delta\to 0_+} 2^{mq}\,\delta^{mq}\,\|f\|_{L^p(\rn)}^q=0.
\end{align}

Finally, letting $\delta\to0_+$ on both sides of \eqref{eq-flz2} (with $\lambda=\delta^{-2}$ therein) and applying \eqref{eq-flz3}, we get
\begin{align*}
\mathscr{D}_{p,\,q}^{m,\,\mu}(f)
&\le \mathscr{D}_{p,\,q}^{m,\,\mu}(f) \lf(\lim_{\delta\to 0_+} \liminf_{j\to\fz} \int_{|h|<\delta} |h|^{mq}\,\varphi_j(h)\,dh \r).
\end{align*}
This implies \eqref{eq-lower} since we have chosen $f$ to satisfy $\mathscr{D}_{p,\,q}^{m,\,\mu}(f)\neq0$.
\end{proof}

Applying Lemmas \ref{Lem-upper} and \ref{Lem-lower}, we now prove Theorem \ref{Thm-BBM-cond}.

\begin{proof}[Proof of Theorem \ref{Thm-BBM-cond}]
From \eqref{eq-upper} and \eqref{eq-lower}, it follows directly that
\begin{align}\label{eq-limit-varphi}
\lim_{\delta\to0_+} \lim_{j\to\fz} \int_{|h|<\delta} |h|^{mq}\,\varphi_j(h)\,dh=1.
\end{align}
By \eqref{eq-limit-varphi} and \eqref{eq-BBM-smalldelta} (with $\operatorname{limL}$
taken to be $\liminf$), we obtain that, for any radial function
$f\in C_{\rm c}^\fz(\rn)$,
\begin{align*}
\lim_{\delta\to 0_+} \liminf_{j\to\fz} \int_{|h|<\delta} \|\Delta_h^mf\|_{L^p(\rn)}^q\,\varphi_j(h)\,dh
= \mathscr{D}_{p,\,q}^{m,\,\mu}(f).
\end{align*}
Noting that \eqref{eq-BBM-varphi} states that
\begin{align*}
\limsup_{j\to\fz} \int_\rn \|\Delta_h^mf\|_{L^p(\rn)}^q\,\varphi_j(h)\,dh
= \mathscr{D}_{p,\,q}^{m,\,\mu}(f),
\end{align*}
we then obtain
\begin{align}\label{eq-BBM-largedz=0}
&\lim_{\delta\to 0_+} \limsup_{j\to\fz} \int_{|h|\ge\delta} \|\Delta_h^mf\|_{L^p(\rn)}^q\,\varphi_j(h)\,dh\\
&\quad\le\lim_{\delta\to 0_+} \left(\limsup_{j\to\fz} \int_\rn - \liminf_{j\to\fz}\int_{|h|\ge\delta}
\right)\|\Delta_h^mf\|_{L^p(\rn)}^q\,\varphi_j(h)\,dh\notag\\
&\quad=0. \notag
\end{align}

Now, we fix some $\delta\in(0,\infty)$. Taking a radial function $0\not\equiv f\in C_{\rm c}^\fz(\rn)$ with
$\supp f\subset B(0, \delta/2)$, we deduce from Lemma \ref{Lem-LXlem2829}(iv) that
\begin{align*}
\|f\|_{L^p(B(0, \,\delta/2))}^q \int_{|h|\ge\delta} \varphi_j(h)\,dh
\le \int_{|h|\ge\delta} \|\Delta_h^mf\|_{L^p(\rn)}^q\,\varphi_j(h)\,dh.
\end{align*}
This, together with \eqref{eq-BBM-largedz=0}, implies
\begin{align*}
\lim_{\delta\to0_+} \limsup_{j\to\infty}\int_{|h|\ge\delta} \varphi_j(h)\,dh=0.
\end{align*}
Consequently, using the fact that $\min\{|h|^{mq},1\}\le 1$, we conclude that
\begin{align*}
\limsup_{j\to\fz} \int_{|h|\ge\delta} \varphi_j(h) \min\lf\{|h|^{mq},\,1 \r\} dh
\le \limsup_{j\to\infty} \int_{|h|\ge\delta} \varphi_j(h)\,dh =0.
\end{align*}
This proves $\VBBMmq$.

Furthermore, for any $\delta\in(0,1)$, by $\VBBMmq$, we have
\begin{align*}
\lim_{j\to\fz} \int_{\delta\le|h|<1} |h|^{mq}\,\varphi_j(h)\,dh
&= \lim_{j\to\fz} \int_{\delta\le|h|<1} \varphi_j(h) \min\lf\{|h|^{mq},\,1 \r\} dh \\
&\le \lim_{j\to\fz} \int_{|h|\ge\delta} \varphi_j(h) \min\lf\{|h|^{mq},\,1 \r\} dh \\
&=0.
\end{align*}
Combining this with \eqref{eq-limit-varphi} gives
\begin{align*}
\lim_{j\to\fz} \int_{|h|<1} |h|^{mq}\,\varphi_j(h)\,dh
= \lim_{\delta\to0_+} \lim_{j\to\fz} \lf(\int_{|h|<\delta}
+ \int_{\delta \le|h|<1} \r)|h|^{mq}\,\varphi_j(h)\,dh=1.
\end{align*}
Thus, we obtain $\NBBMmq$. This completes the proof of Theorem \ref{Thm-BBM-cond}.
\end{proof}

\subsection{Proof of (i) $\Rightarrow$ (iv) in Theorem \ref{thm-BBM}}\label{sec-BBM2}

We establish the implication (i) $\Rightarrow$ (iv) in Theorem \ref{thm-BBM} by proving the following theorem.

\begin{theorem}\label{Thm-Gamma-converge}
Let $p\in[1,\infty)$, $q\in(0,\infty)$, $m\in\nn$ and $\{\varphi_j\}_{j\in\nn}$
satisfy $\BBMmq$. Then, there exist a subsequence $\{\varphi_{j_k}\}_{k\in\nn}$
and a Borel probability measure $\mu$ on $\msn$ such that, for any $f\in W^{m,\,p}(\rn)$,
\begin{align}\label{eq-BBM-Gamma-conv1}
\limsup_{k\to\fz} \mathscr{E}_{p,\,q}^{m,\,\varphi_{j_k}}(f)
\le \mathscr{D}_{p,\,q}^{m,\,\mu}(f),
\end{align}
and, for any $f\in W^{m,\,p}(\rn)$ and any $\{f_k\}_{k\in\nn}\subset L^p(\rn)$
satisfying $f_k\to f$ in $L^p(\rn)$ as $k\to\fz$,
\begin{align}\label{eq-BBM-Gamma-conv2}
\liminf_{k\to\fz} \mathscr{E}_{p,\,q}^{m,\,\varphi_{j_k}}(f_k)
\ge \mathscr{D}_{p,\,q}^{m,\,\mu}(f).
\end{align}
Moreover, in the case $p=1$, the above results remain valid, with
$W^{m,\,p}(\rn)$ and $\mathscr{D}_{p,\,q}^{m,\,\mu}$ replaced by $\BVm(\rn)$
and $\mathscr{D}_{q}^{m,\,\mu}$, respectively.
\end{theorem}

The proofs of \eqref{eq-BBM-Gamma-conv1} and \eqref{eq-BBM-Gamma-conv2} rely on the following three lemmas:
Lemma \ref{Lem-lim-varphi=mu} is the key to \eqref{eq-BBM-Gamma-conv1}, while Lemmas \ref{Lem2-Gamma-converge}
and \ref{Lem-uc-to-0} are the key to \eqref{eq-BBM-Gamma-conv2}.

\begin{lemma}\label{Lem-lim-varphi=mu}
Let all the hypotheses be as in Theorem \ref{Thm-Gamma-converge}. Then there exist a subsequence $\{\varphi_{j_k}\}_{k\in\nn}$
and a Borel probability measure $\mu$ on $\msn$ such that, for any $f\in W^{m,\,p}(\rn)$,
\begin{align}\label{eq-lim-varphi=mu}
\lim_{k\to\fz} \int_{|h|<1} \|(h\cdot\nabla)^mf\|_{L^p(\rn)}^q\, \varphi_{j_k}(h)\,dh
= \mathscr{D}_{p,\,q}^{m,\,\mu}(f).
\end{align}
\end{lemma}

\begin{proof}
For any $j\in\nn$, define a Borel measure on the sphere $\msn$ by
\begin{align*}
d\mu_j(\xi):=\left(\int_0^1
t^{n+mq-1}\,\varphi_j(t\xi)\,dt\right)\,
d\sigma(\xi).
\end{align*}
Then, for any Borel set $E\subset\msn$, we have
\begin{align}\label{eq-def-muj}
\mu_j(E)=\int_E \int_0^1 t^{mq+n-1}\varphi_j(t\xi)\,dt\,d\sigma(\xi)
=\int_{|h|<1} {\bf 1}_E\lf(\frac{h}{|h|} \r) |h|^{mq}\,\varphi_j(h)\,dh.
\end{align}
Clearly, it follows from $\NBBMmq$ that
\begin{align*}
\lim_{j\to\fz} \mu_j(\msn)
= \lim_{j\to\fz} \int_{|h|<1} |h|^{mq}\,\varphi_j(h)\,dh=1.
\end{align*}
Consequently,
\begin{align*}
\sup_{j\in\nn} \mu_j(\msn)<\fz.
\end{align*}
Since $\msn$ is compact, we deduce from the Prokhorov theorem (see, for example, \cite[Theorem 8.6.2]{Bogachev-Book-2007})
that there exist a subsequence $\{\mu_{j_k}\}_{k\in\nn}$ and a Borel measure $\mu$ on $\msn$ such that for any $\psi\in C(\msn)$,
\begin{align*}
\lim_{k\to\fz} \int_\msn \psi(\xi)\,d\mu_{j_k}(\xi)
= \int_\msn \psi(\xi)\,d\mu(\xi).
\end{align*}
Taking $\psi\equiv 1$ yields $\mu(\msn)=1$.

Now, let $f\in W^{m,\,p}(\rn)$. Since the function $g(\xi)=\|(\xi\cdot\nabla)^mf\|_{L^p(\rn)}$ is
continuous on $\msn$, we have
\begin{align*}
\lim_{k\to\fz} \int_\msn \|(\xi\cdot\nabla)^mf\|_{L^p(\rn)}^q\, d\mu_{j_k}(\xi)
= \int_\msn \|(\xi\cdot\nabla)^mf\|_{L^p(\rn)}^q\,d\mu(\xi)
= \mathscr{D}_{p,\,q}^{m,\,\mu}(f).
\end{align*}
Thus, to obtain \eqref{eq-lim-varphi=mu}, it suffices to prove
\begin{align}\label{eq-claim-muj}
\int_{|h|<1} \lf\|\lf(\frac{h}{|h|}\cdot\nabla \r)^mf \r\|_{L^p(\rn)}^q |h|^{mq}\,\varphi_{j_k}(h)\,dh
= \int_\msn \|(\xi\cdot\nabla)^mf\|_{L^p(\rn)}^q\, d\mu_{j_k}(\xi).
\end{align}

Indeed, for any Borel set $E\subset\msn$, it follows from \eqref{eq-def-muj} that
\begin{align*}
\int_{|h|<1} {\bf 1}_E\lf(\frac{h}{|h|} \r) |h|^{mq}\,\varphi_{j_k}(h)\,dh
= \mu_{j_k}(E)
= \int_\msn {\bf 1}_E(\xi)\,d\mu_{j_k}(\xi).
\end{align*}
Since any nonnegative continuous function $g$ can be approximated from below by
nonnegative simple functions, it follows by linearity and the monotone
convergence theorem that
\begin{align*}
\int_{|h|<1} g\lf(\frac{h}{|h|} \r) |h|^{mq}\,\varphi_{j_k}(h)\,dh
= \int_\msn g(\xi)\,d\mu_{j_k}(\xi).
\end{align*}
Taking $g(\xi)=\|(\xi\cdot\nabla)^mf\|_{L^p(\rn)}^q$ in this identity yields
\eqref{eq-claim-muj}.
This completes the proof.
\end{proof}

Choose a radial function $0\le\eta\in C_{\rm c}^\fz(\rn)$ such that $\int_\rn \eta(x)\,dx=1$.
For any $\uc\in(0,\infty)$, define
\begin{align*}
\eta_\uc(x):= \uc^{-n}\,\eta(\uc^{-1}x)\qquad\text{for all}\ \, x\in\rn.
\end{align*}
Let $p\in[1,\infty)$. For each $f_k\in L^p(\rn)$ and $f\in L^p(\rn)$, we define
\begin{align}\label{eq-fjfvez}
\begin{cases}
 f_\uc:= f\ast\eta_\uc; \\
 f_{k,\uc}:= f_k\ast\eta_\uc \quad\text{for all }\, k\in\nn.
\end{cases}
\end{align}
Clearly, $f_\uc$ and $f_{k,\uc}$ belong to $C^\fz(\rn)$. The following lemma will be used in the proof of \eqref{eq-BBM-Gamma-conv2}.

\begin{lemma}\label{Lem2-Gamma-converge}
Let all the hypotheses be as in Theorem \ref{Thm-Gamma-converge}.
Then, for any $f\in W^{m,\,p}(\rn)$, any $\{f_k\}_{k\in\nn}\subset L^p(\rn)$ that converges to
 $f$ in $L^p(\rn)$, and any $\uc\in(0,\fz)$,
\begin{align}\label{eq-claim2-eq2}
\lim_{\delta\to0_+} \lim_{k\to\fz} \int_{|h|<\delta} \|\Delta_h^mf_{k,\uc}\|_{L^p(\rn)}^q\,\varphi_{j_k}(h)\,dh
= \mathscr{D}_{p,\,q}^{m,\,\mu}(f_\uc),
\end{align}
where $f_{k,\uc}$ and $f_\uc$ are defined as in \eqref{eq-fjfvez}, the subsequence $\{\varphi_{j_k}\}_{k\in\nn}$ and the measure $\mu$ are from Lemma \ref{Lem-lim-varphi=mu}.
\end{lemma}

\begin{proof}
Let $\delta\in(0,1)$. Then, it follows from $\VBBMmq$ that
\begin{align*}
\lim_{j\to\fz} \int_{\delta\le|h|<1} |h|^{mq}\,\varphi_j(h)\,dh
&= \lim_{j\to\fz} \int_{\delta\le|h|<1} \varphi_j(h) \min\lf\{|h|^{mq},\,1 \r\} dh \\
&\le \lim_{j\to\fz} \int_{|h|\ge\delta} \varphi_j(h) \min\lf\{|h|^{mq},\,1 \r\} dh \\
&= 0,
\end{align*}
which, together with $\NBBMmq$, gives
\begin{align}\label{eq-BBM-varphi-cond1-dz}
\lim_{j\to\fz} \int_{|h|<\delta} |h|^{mq}\,\varphi_j(h)\,dh
&= \lim_{j\to\fz} \lf(\int_{|h|<1}-\int_{\delta\le|h|<1}\r) |h|^{mq}\,\varphi_j(h)\,dh \\ \notag
&= \lim_{j\to\fz} \int_{|h|<1} |h|^{mq}\,\varphi_j(h)\,dh=1.
\end{align}

Let $f_{k,\uc}$ and $f_{\uc}$ be defined as in \eqref{eq-fjfvez}, and let $\{\varphi_{j_k}\}_{k\in\nn}$ and $\mu$
be given as in Lemma \ref{Lem-lim-varphi=mu}. Assume for the moment that we have proved
\begin{align}\label{eq-fjvez-lim=0}
\lim_{\delta\to0_+} \lim_{k\to\fz} \int_{|h|<\delta} \|\Delta_h^mf_{k,\uc}-(h\cdot\nabla)^mf_\uc\|_{L^p(\rn)}^q\,\varphi_{j_k}(h)\,dh=0.
\end{align}
Then, by \eqref{eq-fjvez-lim=0} and \eqref{eq-norm-minus}, we obtain
\begin{align*}
&\lim_{\delta\to 0_+} \lim_{k\to\fz} \int_{|h|<\delta} \|\Delta_h^mf_{k,\uc}\|_{L^p(\rn)}^q\,\varphi_{j_k}(h)\,dh \\ \notag
&\quad= \lim_{\delta\to0_+} \lim_{k\to\fz} \int_{|h|<\delta} \|(h\cdot\nabla)^m f_\uc\|_{L^p(\rn)}^q\,\varphi_{j_k}(h)\,dh.
\end{align*}
Moreover, we use \eqref{eq-lim-varphi=mu} in Lemma \ref{Lem-lim-varphi=mu} to deduce
\begin{align*}
\lim_{k\to\fz} \int_{|h|<\delta} \|(h\cdot\nabla)^m f_\uc\|_{L^p(\rn)}^q\,\varphi_{j_k}(h)\,dh
= \mathscr{D}_{p,\,q}^{m,\,\mu}(f_\uc).
\end{align*}
In view of the last two formulae and \eqref{eq-BBM-varphi-cond1-dz}, we have
\begin{align*}
\lim_{\delta\to 0_+} \lim_{k\to\fz} \int_{|h|<\delta} \|\Delta_h^mf_{k,\uc}\|_{L^p(\rn)}^q\,\varphi_{j_k}(h)\,dh
= \mathscr{D}_{p,\,q}^{m,\,\mu}(f_\uc),
\end{align*}
which is exactly \eqref{eq-claim2-eq2}.

It remains to prove \eqref{eq-fjvez-lim=0}; the argument is analogous to that for \eqref{eq1-Delta-Nabla}.
We give the details below. Since each $f_{k,\uc}\in C^\infty(\rn)$, it follows from Lemma \ref{lem-diffLpNorm} that for all $h\in\rn$,
\begin{align}\label{eq-diff-LpNorm-fjvez}
\|\Delta_h^mf_{k,\uc}-(h\cdot\nabla)^mf_{k,\uc}\|_{L^p(\rn)}
\le m|h|^{m+1} \sup_{\xi\in\msn} \|(\xi\cdot\nabla)^{m+1}f_{k,\uc}\|_{L^p(\rn)}.
\end{align}
For any $\xi\in\msn$, by \eqref{eq-fjfvez}, we have
\begin{align}\label{eq-normLp}
\|(\xi\cdot\nabla)^{m+1}f_{k,\uc}\|_{L^p(\rn)}
&\le \sum_{|\alpha|=m+1} \frac{(m+1)!}{\alpha!}\,\|\partial^\alpha f_{k,\uc}\|_{L^p(\rn)} \\
&= \sum_{|\alpha|=m+1} \frac{(m+1)!}{\alpha!}\,\|\partial^\alpha (f_k\ast\eta_\uc)\|_{L^p(\rn)}\notag\\
&= \sum_{|\alpha|=m+1} \frac{(m+1)!}{\alpha!}\,\|f_k\ast(\partial^\alpha\eta_\uc)\|_{L^p(\rn)}\notag\\
&\le \| f_k\|_{L^p(\rn)} \sum_{|\alpha|=m+1} \frac{(m+1)!}{\alpha!}\,\|\partial^\alpha\eta_\uc\|_{L^1(\rn)}. \notag
\end{align}
By the hypothesis of Theorem \ref{Thm-Gamma-converge}, the sequence $\{f_k\}_{k\in\nn}$
converges to $f$ in $L^p(\rn)$, which implies
\begin{align*}
M_f:=\sup_{k\in\nn} \|f_k\|_{L^p(\rn)}<\fz.
\end{align*}
Set
\begin{align*}
C_{\eta,\,\uc}:= \sum_{\ell=m}^{m+1} \sum_{|\alpha|=\ell} \frac{\ell!}{\alpha!}\,\|\partial^\alpha\eta_\uc\|_{L^1(\rn)}.
\end{align*}
Then, we deduce from \eqref{eq-diff-LpNorm-fjvez} and \eqref{eq-normLp} that
\begin{align*}
\|\Delta_h^mf_{k,\uc}-(h\cdot\nabla)^mf_{k,\uc}\|_{L^p(\rn)}
\le mM_f\,C_{\eta,\,\uc}\,|h|^{m+1}.
\end{align*}
Moreover, analogously to \eqref{eq-normLp}, we have
\begin{align*}
\|(h\cdot\nabla)^m f_{k,\uc}-(h\cdot\nabla)^mf_\uc\|_{L^p(\rn)}
&= \lf\|(h\cdot\nabla)^m \lf((f_k-f)\ast\eta_\uc \r)\r\|_{L^p(\rn)} \\
&\le |h|^m \sup_{\xi\in\msn} \lf\|(\xi\cdot\nabla)^m \lf((f_k-f)\ast\eta_\uc \r)\r\|_{L^p(\rn)} \notag \\
&\le |h|^m\,\|f_k-f\|_{L^p(\rn)} \sum_{|\alpha|=m} \frac{m!}{\alpha!}\,\|\partial^\alpha\eta_\uc\|_{L^1(\rn)} \notag \\
&\le C_{\eta,\,\uc}\,|h|^m\,\|f_k-f\|_{L^p(\rn)}. \notag
\end{align*}
By the last two estimates and the Minkowski inequality, we see that
\begin{align*}
&\|\Delta_h^mf_{k,\uc}-(h\cdot\nabla)^mf_{\uc}\|_{L^p(\rn)} \\
&\quad\le \|\Delta_h^mf_{k,\uc}-(h\cdot\nabla)^mf_{k,\uc}\|_{L^p(\rn)}
+ \|(h\cdot\nabla)^m f_{k,\uc}-(h\cdot\nabla)^mf_\uc\|_{L^p(\rn)} \notag \\
&\quad\le C\,|h|^m \lf(|h|+\|f_k-f\|_{L^p(\rn)}\r), \notag
\end{align*}
where $C$ is a positive constant depending on $m$, $n$, $M_f$, $\eta$ and $\uc$.
Consequently,
\begin{align*}
&\limsup_{k\to\fz} \int_{|h|<\delta} \lf\|\Delta_h^mf_{k,\uc}-(h\cdot\nabla)^mf_\uc \r\|_{L^p(\rn)}^q \varphi_{j_k}(h)\,dh \\
&\quad\le C^q \limsup_{k\to\fz} \int_{|h|<\delta} \lf(|h|+\|f_k-f\|_{L^p(\rn)} \r)^q |h|^{mq}\,\varphi_{j_k}(h)\,dh \\
&\quad\le C^q \limsup_{k\to\fz} \lf[\lf(\delta+\|f_k-f\|_{L^p(\rn)} \r)^q
\lf(\int_{|h|<\delta} |h|^{mq}\,\varphi_{j_k}(h)\,dh \r)\r] \\
&\quad\le (C\delta)^q \limsup_{k\to\fz} \int_{|h|<\delta} |h|^{mq}\,\varphi_{j_k}(h)\,dh \\
&\quad= (C\delta)^q,
\end{align*}
where the penultimate step follows from $\lim_{k\to\fz} \|f_k-f\|_{L^p(\rn)}=0$, and the last step
is due to \eqref{eq-BBM-varphi-cond1-dz}. Finally, letting $\delta\to0_+$ yields \eqref{eq-fjvez-lim=0}.
This finishes the proof of Lemma \ref{Lem2-Gamma-converge}.
\end{proof}

\begin{lemma}\label{Lem-uc-to-0}
Let $p\in[1,\fz)$ and $m\in\nn$.
Then, for any $f\in W^{m,\,p}(\rn)$,
\begin{align}\label{eq-uc-to-0-Wmp}
\lim_{\uc\to0_+} \|(\xi\cdot\nabla)^mf_\uc\|_{L^p(\rn)}
= \|(\xi\cdot\nabla)^mf\|_{L^p(\rn)}
\qquad\text{uniformly in}\,\ \xi\in\msn.
\end{align}
Moreover, for any $f\in\BVm(\rn)$,
\begin{align}\label{eq-uc-to-0-BVm}
\lim_{\uc\to0_+} \|(\xi\cdot\nabla)^mf_\uc\|_{\mathcal{M}(\rn)}
= \|(\xi\cdot\nabla)^mf\|_{\mathcal{M}(\rn)}
\qquad\text{uniformly in}\,\ \xi\in\msn.
\end{align}
Here, in \eqref{eq-uc-to-0-Wmp} and \eqref{eq-uc-to-0-BVm}, the family $\{f_\uc\}_{\uc\in(0,\infty)}$ is defined as in \eqref{eq-fjfvez}.
\end{lemma}

\begin{proof}
Let $f\in W^{m,\,p}(\rn)$. Since $p\in[1,\infty)$ and $\{\eta_\uc\}_{\uc\in(0,\fz)}$
is an approximation of the identity, it follows that for every $\alpha\in\zz_+^n$ with $|\alpha|=m$,
\begin{align*}
\lim_{\uc\to0_+} \|\partial^\alpha f_\uc-\partial^\alpha f\|_{L^p(\rn)}
= \lim_{\uc\to0_+} \|(\partial^\alpha f)\ast\eta_\uc-\partial^\alpha f\|_{L^p(\rn)}=0.
\end{align*}
By the multinomial expansion
and the fact that $|\xi^\alpha|\le1$ for $\xi\in\msn$, we have
\begin{align*}
\|(\xi\cdot\nabla)^mf_\uc-(\xi\cdot\nabla)^mf\|_{L^p(\rn)}
\le \sum_{|\alpha|=m} \frac{m!}{\alpha!}\|\partial^\alpha f_\uc-\partial^\alpha f\|_{L^p(\rn)},
\end{align*}
whose right-hand side is independent of $\xi\in\msn$. Combining the last two formulae shows that
\begin{align*}
\lim_{\uc\to0_+} \|(\xi\cdot\nabla)^mf_\uc-(\xi\cdot\nabla)^mf\|_{L^p(\rn)}=0
\qquad\text{uniformly in $\xi\in\msn$}.
\end{align*}
Thus, we obtain \eqref{eq-uc-to-0-Wmp}.

Next, let $f\in\BVm(\rn)$.
For any $\xi\in\msn$, set
\begin{align*}
\mu_\xi:=(\xi\cdot\nabla)^m f
=\sum_{|\alpha|=m}\frac{m!}{\alpha!}\,\xi^\alpha\,\partial^\alpha f\in\mathcal{M}(\mathbb{R}^n).
\end{align*}
Since $f_\varepsilon=f\ast\eta_\varepsilon$, we have
\begin{align*}
(\xi\cdot\nabla)^m f_\varepsilon
=\mu_\xi\ast\eta_\varepsilon.
\end{align*}
Thus, by the Young inequality for convolution with
finite measures,
\begin{align}\label{eq-Young-BV}
\|(\xi\cdot\nabla)^m f_\varepsilon\|_{\mathcal{M}(\mathbb{R}^n)}
=\|(\xi\cdot\nabla)^m f_\varepsilon\|_{L^1(\mathbb{R}^n)}
&=\|\mu_\xi\ast\eta_\varepsilon\|_{L^1(\mathbb{R}^n)} \\
&\le \|\mu_\xi\|_{\mathcal{M}(\mathbb{R}^n)}\,\|\eta_\varepsilon\|_{L^1(\mathbb{R}^n)}
= \|\mu_\xi\|_{\mathcal{M}(\mathbb{R}^n)}. \notag
\end{align}

Let $\psi\in C_{\rm c}(\rn)$ such that $|\psi|\le 1$ on $\rn$. By the uniform continuity of $\psi$, we obtain
$\eta_\uc\ast\psi\to\psi$ uniformly on $\rn$ as $\uc\to0_+$. Since $\mu_\xi$ is a finite measure and $\eta_\uc$
is radial, it follows that
\begin{align*}
\int_\rn \psi\,d\mu_\xi
&= \lim_{\uc\to0_+} \int_\rn (\eta_\uc\ast \psi)\,d\mu_\xi\\
&= \lim_{\uc\to0_+} \int_\rn (\eta_\uc\ast\mu_\xi)(x)\, \psi(x)\,dx\\
&= \lim_{\uc\to0_+} \int_\rn (\xi\cdot\nabla)^m f_\uc(x)\, \psi(x)\,dx\\
&\le \liminf_{\uc\to0_+} \|(\xi\cdot\nabla)^m f_\uc\|_{L^1(\rn)}\\
&= \liminf_{\uc\to0_+}\|(\xi\cdot\nabla)^m f_\uc\|_{\mathcal{M}(\rn)},
\end{align*}
where the second equality follows from the Fubini theorem and the symmetry of $\eta_\uc$. Taking the
supremum over all such $\psi$ and using the definition of $\|\cdot\|_{\mathcal{M}(\rn)}$, we get
\begin{align}\label{eq-reverse-BV}
\|\mu_\xi\|_{\mathcal{M}(\rn)}
\le \liminf_{\uc\to0}\|(\xi\cdot\nabla)^mf_\uc\|_{\mathcal{M}(\rn)}.
\end{align}
Combining \eqref{eq-reverse-BV} with \eqref{eq-Young-BV}, we conclude that
\begin{align*}
\lim_{\uc\to0_+}
\|(\xi\cdot\nabla)^m f_\uc\|_{\mathcal{M}(\rn)}
= \|\mu_\xi\|_{\mathcal{M}(\rn)}.
\end{align*}
This establishes the pointwise convergence in \eqref{eq-uc-to-0-BVm} for every fixed $\xi\in\msn$.

We now show that the convergence in \eqref{eq-uc-to-0-BVm} is uniform in $\xi$. For any two finite signed measures
$\mu$ and $\nu$ on $\rn$, we have
\begin{align*}
\lf|\|\mu\|_{\mathcal M(\rn)}
-\|\nu\|_{\mathcal M(\rn)}\r|
\le\|\mu-\nu\|_{\mathcal M(\rn)}.
\end{align*}
Combining this with \eqref{eq-def-xi-cdot-nabla}
and \eqref{eq-BV-norm}, for any $\uc>0$
and $\xi,\zeta\in\msn$, we obtain
\begin{align*}
\lf|\|(\xi\cdot\nabla)^mf_\uc\|_{\mathcal M(\rn)}-\|(\zeta\cdot\nabla)^mf_\uc\|_{\mathcal M(\rn)} \r|
&\le \|(\xi\cdot\nabla)^mf_\uc
-(\zeta\cdot\nabla)^mf_\uc\|_{\mathcal M(\rn)}\\
&\le \sum_{\{\alpha\in\zz_+^n:\,|\alpha|=m\}} \frac{m!}{\alpha!}\,|\xi^\alpha-\zeta^\alpha|\,\|\partial^\alpha f_\uc\|_{\mathcal M(\rn)}\\
&\le m|\xi-\zeta|\,\|f\|_{\mathringBVm(\rn)},
\end{align*}
where we have used the fact that for $|\alpha|=m$ and $\xi,\zeta\in\msn$,
\begin{align*}
|\xi^\alpha-\zeta^\alpha|\le m|\xi-\zeta|,
\end{align*}
together with the Young inequality for convolution with finite measures. Thus, the functions
\begin{align*}
\xi\mapsto\|(\xi\cdot\nabla)^mf_\uc\|_{\mathcal{M}(\rn)}
\end{align*}
are equi-Lipschitz on $\msn$. Since $\msn$ is compact, the pointwise convergence established above is uniform by the
Arzel\'a--Ascoli theorem. This completes the proof of the uniform convergence in \eqref{eq-uc-to-0-BVm}.
\end{proof}

Now, with Lemmas \ref{Lem-lim-varphi=mu}, \ref{Lem2-Gamma-converge} and \ref{Lem-uc-to-0} at hand,
we are ready to prove Theorem \ref{Thm-Gamma-converge}.

\begin{proof}[Proof of Theorem \ref{Thm-Gamma-converge}]
The proof is divided into three steps.

\medskip

\textbf{Step 1:\, Proof of \eqref{eq-BBM-Gamma-conv1} for functions in $W^{m,\,p}(\rn)$.}
Given any $f\in W^{m,\,p}(\rn)$, we write
\begin{align}\label{eq-BBMG-twopart}
\limsup_{j\to\fz} \mathscr{E}_{p,\,q}^{m,\,\varphi_j}(f)
&\le \limsup_{j\to\fz} \int_{|h|<1} \|\Delta_h^mf\|_{L^p(\rn)}^q\,\varphi_j(h)\,dh \\ \notag
&\quad+ \limsup_{j\to\fz} \int_{|h|\ge1} \|\Delta_h^mf\|_{L^p(\rn)}^q\,\varphi_j(h)\,dh.
\end{align}
Using $f\in W^{m,\,p}(\rn)$ and Lemma \ref{Lem-LXlem2829}(ii) yields that for any $h\in\rn$,
\begin{align}\label{eq1-remark-BV}
\|\Delta_h^mf \|_{L^p(\rn)} \le \|(h\cdot\nabla)^mf \|_{L^p(\rn)},
\end{align}
which implies
\begin{align*}
\int_{|h|<1} \|\Delta_h^mf\|_{L^p(\rn)}^q\,\varphi_j(h)\,dh
 &\le \int_{|h|<1} \|(h\cdot\nabla)^mf\|_{L^p(\rn)}^q\,\varphi_j(h)\,dh.
\end{align*}
By Lemma \ref{Lem-lim-varphi=mu}, we can find a subsequence $\{\varphi_{j_k}\}_{k\in\nn}$ and
a Borel probability measure $\mu$ on $\msn$ such that
\begin{align*}
\limsup_{k\to\fz} \int_{|h|<1} \|(h\cdot\nabla)^mf\|_{L^p(\rn)}^q\,\varphi_{j_k}(h)\,dh
= \mathscr{D}_{p,\,q}^{m,\,\mu}(f).
\end{align*}
From these last two inequalities, it follows that
\begin{align}\label{eq-BBMG-smallh}
\limsup_{k\to\fz} \int_{|h|<1} \|\Delta_h^mf\|_{L^p(\rn)}^q\,\varphi_{j_k}(h)\,dh
\le \limsup_{k\to\fz} \int_{|h|<1} \|(h\cdot\nabla)^mf\|_{L^p(\rn)}^q\,\varphi_{j_k}(h)\,dh
= \mathscr{D}_{p,\,q}^{m,\,\mu}(f).
\end{align}
Moreover, for any $h\in\rn$, by Lemma \ref{Lem-LXlem2829}(i), we have
\begin{align*}
\|\Delta_h^mf\|_{L^p(\rn)}
\le 2^m\,\|f\|_{L^p(\rn)}.
\end{align*}
From this and $\VBBMmq$, it follows that
\begin{align}\label{eq-BBMG-largeh}
\limsup_{k\to\fz} \int_{|h|\ge1} \|\Delta_h^mf\|_{L^p(\rn)}^q\,\varphi_{j_k}(h)\,dh
&\le 2^{mq}\,\|f\|_{L^p(\rn)}^q \limsup_{k\to\fz} \int_{|h|\ge1} \varphi_{j_k}(h)\,dh \\
&= 2^{mq}\,\|f\|_{L^p(\rn)}^q \lim_{k\to\fz} \int_{|h|\ge1} \varphi_{j_k}(h) \min\lf\{|h|^{mq},\,1 \r\} dh \notag \\
&= 0. \notag
\end{align}
Finally, substituting \eqref{eq-BBMG-smallh} and \eqref{eq-BBMG-largeh} into \eqref{eq-BBMG-twopart} applied along $j=j_k$,
we conclude that \eqref{eq-BBM-Gamma-conv1} holds.

\medskip

\textbf{Step 2:\, Proof of \eqref{eq-BBM-Gamma-conv2} for functions in $W^{m,\,p}(\rn)$.}
Fix $k\in\nn$ and $\uc\in(0,\infty)$. By the linearity of the convolution, for any $h\in\rn$, we have
\begin{align*}
\Delta_h^m f_{k,\uc} = \Delta_h^m(f_k\ast\eta_\uc) = (\Delta_h^mf_k)\ast\eta_\uc,
\end{align*}
which, together with the Young inequality and $\|\eta_\uc\|_{L^1(\rn)}=1$, implies
\begin{align*}
\|\Delta_h^mf_{k,\uc}\|_{L^p(\rn)}
= \|(\Delta_h^mf_k)\ast\eta_\uc\|_{L^p(\rn)}
\le \|\Delta_h^mf_k\|_{L^p(\rn)}.
\end{align*}
From this and \eqref{eq-claim2-eq2} in Lemma \ref{Lem2-Gamma-converge}, it follows that
\begin{align}\label{eq1-Gamma-conv}
\mathscr{D}_{p,\,q}^{m,\,\mu}(f_\uc)
= \lim_{\delta\to0_+} \liminf_{k\to\fz} \int_{|h|<\delta} \|\Delta_h^mf_{k,\uc}\|_{L^p(\rn)}^q\,\varphi_{j_k}(h)\,dh
\le \liminf_{k\to\fz} \mathscr{E}_{p,\,q}^{m,\,\varphi_{j_k}}(f_k).
\end{align}
For any $f\in W^{m,\,p}(\rn)$, the uniform convergence in \eqref{eq-uc-to-0-Wmp} implies
\begin{align*}
\lim_{\uc\to0_+} \int_\msn \|(\xi\cdot\nabla)^mf_\uc\|_{L^p(\rn)}^q\,d\mu(\xi)
= \int_\msn \|(\xi\cdot\nabla)^mf\|_{L^p(\rn)}^q\,d\mu(\xi),
\end{align*}
or equivalently,
\begin{align}\label{eq2-Gamma-conv}
\lim_{\uc\to0_+} \mathscr{D}_{p,\,q}^{m,\,\mu}(f_\uc)
= \mathscr{D}_{p,\,q}^{m,\,\mu}(f).
\end{align}
Now, combining \eqref{eq1-Gamma-conv} with \eqref{eq2-Gamma-conv} yields \eqref{eq-BBM-Gamma-conv2}.

\medskip

\textbf{Step 3:\, Proofs of \eqref{eq-BBM-Gamma-conv1} and \eqref{eq-BBM-Gamma-conv2} for functions in $\BVm(\rn)$.}
Let $f\in \BVm(\rn)$. In this case, the counterpart of \eqref{eq-lim-varphi=mu} reads
\begin{align}\label{eq-BV-lim-varphi=mu}
\lim_{k\to\fz} \int_{|h|<1} \|(h\cdot\nabla)^mf\|_{\mathcal M(\rn)}^q\, \varphi_{j_k}(h)\,dh
= \mathscr{D}_{q}^{m,\,\mu}(f).
\end{align}
This can be obtained by following the same lines as in the proof of Lemma \ref{Lem-lim-varphi=mu},
with every $\|\cdot\|_{L^1(\rn)}$-norm replaced by the total variation norm
$\|\cdot\|_{\mathcal M(\rn)}$.

Now we run the argument in \textbf{Step 1} using \eqref{eq-BV-lim-varphi=mu}.
The only other modification is to replace \eqref{eq1-remark-BV} by the following estimate
(see \eqref{eq-Delta-BV-upper}):
\begin{align*}
\|\Delta_{h}^mf\|_{L^1(\rn)}
\le \|(h\cdot\nabla)^mf\|_{\mathcal{M}(\rn)}.
\end{align*}
Then we obtain
\begin{align*}
\limsup_{k\to\fz} \mathscr{E}_{1,\,q}^{m,\,\varphi_{j_k}}(f)
\le \mathscr{D}_q^{m,\,\mu}(f).
\end{align*}
This is exactly the counterpart of \eqref{eq-BBM-Gamma-conv1}.

Even for $f\in \BVm(\rn)$, the functions $f_{k,\uc}$ and $f_\uc$ defined in \eqref{eq-fjfvez} belong to $W^{m,\,1}(\rn)$. Since $f_k\to f$ in $L^1(\rn)$, it follows that the limit in \eqref{eq-claim2-eq2} remains valid, with
$\mathscr{D}_{p,\,q}^{m,\,\mu}(f_\uc)=\mathscr{D}_{q}^{m,\,\mu}(f_\uc)$. Consequently, the current version of
\eqref{eq1-Gamma-conv} becomes
\begin{align*}
\mathscr{D}_q^{m,\,\mu}(f_\uc)
\le \liminf_{k\to\fz} \mathscr{E}_{1,\,q}^{m,\,\varphi_{j_k}}(f_k).
\end{align*}
Further, the uniform convergence in \eqref{eq-uc-to-0-BVm}
 gives
\begin{align*}
\lim_{\uc\to0_+} \mathscr{D}_q^{m,\,\mu}(f_\uc)
= \mathscr{D}_q^{m,\,\mu}(f).
\end{align*}
Therefore, for any $f\in \BVm(\rn)$ and any $\{f_k\}_{k\in\nn}\subset L^1(\rn)$
with $f_k\to f$ in $L^1(\rn)$ as $k\to\fz$, we obtain
\begin{align*}
\mathscr{D}_q^{m,\,\mu}(f)\le \liminf_{k\to\fz} \mathscr{E}_{1,\,q}^{m,\,\varphi_{j_k}}(f_k),
\end{align*}
which is the desired inequality \eqref{eq-BBM-Gamma-conv2}.
\end{proof}

\subsection{Proof of Theorem \ref{thm-BBM-radial}}\label{sec-BBM4}

\begin{proof}[Proof of Theorem \ref{thm-BBM-radial}]
Tracing the proof of Theorem \ref{thm-BBM}, we note that passing to a
subsequence is used only in the proof of (i) $\Rightarrow$ (iv), that is,
in showing that \(\BBMmq\) implies strong \(\Gamma\)-convergence.

When we revisit the proof of this part (that is, the arguments in Section
\ref{sec-BBM2}), we find that passing to a subsequence happens only in the
proof of Lemma \ref{Lem-lim-varphi=mu}, where we define
\begin{align*}
d\mu_{j}(\xi):= \lf(\int_0^1 t^{n+mq}\,\varphi_{j}(t\xi)\,\frac{dt}{t} \r) d\sigma(\xi),
\end{align*}
and then extract a subsequence $\{\mu_{j_k}\}_{k\in\nn}$ from
$\{\mu_j\}_{j\in\nn}$, verifying that it converges vaguely to a Borel
probability measure $\mu$ on $\msn$.

Now, assume that each $\varphi_j$ is radial. For simplicity, we write
$\varphi_j(t):=\varphi_j(h)$ whenever $|h|=t$. Then
\begin{align*}
d\mu_{j}(\xi)= \lf(\int_0^1 t^{n+mq}\,\varphi_{j}(t)\,\frac{dt}{t} \r) d\sigma(\xi).
\end{align*}
Using spherical coordinates, we have
\begin{align*}
\int_0^1 t^{n+mq}\,\varphi_{j}(t)\,\frac{dt}{t}
&=\fint_\msn \lf(\int_0^1 t^{n+mq}\,\varphi_{j}(t)\,\frac{dt}{t} \r) d\sigma(\xi) \\
&= \frac 1{\omega_{n-1}} \int_{|h|<1} |h|^{mq}\,\varphi_j(h)\,dh. \notag
\end{align*}
By $\NBBMmq$, we have
\begin{align*}
\lim_{j\to\infty} \int_0^1 t^{n+mq}\,\varphi_{j}(t)\,\frac{dt}{t}
= \lim_{j\to\infty}\left(\frac 1{\omega_{n-1}} \int_{|h|<1} |h|^{mq}\,\varphi_j(h)\,dh\right)
=\frac 1{\omega_{n-1}}.
\end{align*}
It follows that for any Borel set $E\subset\msn$,
\begin{align*}
\lim_{j\to\infty} \mu_j(E) = \frac1{\omega_{n-1}}\,\sigma(E).
\end{align*}
This proves that Lemma \ref{Lem-lim-varphi=mu} holds without passing to a subsequence, with $d\mu$ therein
taken to be $(\omega_{n-1})^{-1}\,d\sigma$. Consequently, all the equivalences in Theorem \ref{thm-BBM-radial}
hold along the full sequence.
\end{proof}

\subsection{A counterexample}\label{ss-eg}

The following counterexample shows that, for nonradial kernel sequences $\{\varphi_j\}_{j\in\nn}$,
passing to a subsequence is, in general, necessary for the implication (i) $\Rightarrow$ (iv)
in Theorem \ref{thm-BBM}.

\begin{example}\label{eg-BBM-counterexam}
\rm Let $p=q=2$, $m=1$ and $n=2$. For any $h=(h_1,h_2)\in\rr^2\setminus\{0\}$ and $j\in\nn$, let
\begin{align*}
\varphi_j(h):=
\begin{cases}
\dis \frac{3j^3}{\pi}\,\frac{h_1^2}{|h|^3}\,{\bf 1}_{\{|h|<1/j\}}(h)
&\quad\text{if}\,\ j\,\ \text{is even}; \\[3mm]
\dis \frac{3j^3}{\pi}\,h_2^2\,{\bf 1}_{\{|h|<1/j\}}(h)
&\quad\text{if}\,\ j\,\ \text{is odd}.
\end{cases}
\end{align*}
Each $\varphi_j$ is nonradial.
For any $j\in\nn$, by symmetry and polar coordinates, we have
\begin{align*}
\int_{|h|<1/j} |h|^2\,\frac{h_1^2}{|h|^3}\,dh
&= \int_{|h|<1/j} |h|^2\,\frac{h_2^2}{|h|^3}\,dh \\
&= \frac12 \int_{|h|<1/j} |h|^{-1}\,(h_1^2+h_2^2)\,dh
= \frac 12 \int_{|h|<1/j} |h|\,dh
= \frac{\pi}{3j^3}
\end{align*}
and, hence,
\begin{align*}
\int_{|h|<1} |h|^2\,\varphi_j(h)\,dh=1.
\end{align*}
Moreover, for any $\delta\in(0,1]$, if $j>2/\delta$, then
\begin{align*}
\int_{|h|\ge\delta} \varphi_j(h) \min\lf\{|h|^2,\,1 \r\} dh=0.
\end{align*}
Thus, the kernel sequence $\{\varphi_j\}_{j\in\nn}$ satisfies
the $\BBMmq$ condition with $m=1$ and $q=2$.

Now, for any $x=(x_1,x_2)\in\rr^2$, let
\begin{align*}
f(x):=e^{-(x_1^2+4x_2^2)}.
\end{align*}
Obviously, $f\in W^{1,\,2}(\rr^2)$.
For any $h=(h_1,h_2)\in\rr^2$, we have
\begin{align*}
\|\Delta_h f\|_{L^2(\rr^2)}^2
&=2\int_{\rr^2} f(x)^2\,dx
-2\int_{\rr^2} f(x)\,f(x+h)\,dx.
\end{align*}
A direct computation gives
\begin{align*}
\int_{\rr^2}f(x)^2\,dx
= \lf(\int_\rr e^{-2x_1^2}\,dx_1 \r) \lf(\int_\rr e^{-8x_2^2}\,dx_2 \r)
= \frac{\pi}{4}.
\end{align*}
Moreover, writing
\begin{align*}
f(x)\,f(x+h)
&= \exp\left(-2\left(x_1+\frac{h_1}{2} \right)^2
- 8\left(x_2+\frac{h_2}{2}\right)^2-\frac{h_1^2}{2}-2h_2^2 \right),
\end{align*}
we obtain
\begin{align*}
\int_{\rr^2} f(x)\,f(x+h)\,dx
= \frac{\pi}{4}\,e^{-(2^{-1}h_1^2+2h_2^2)}.
\end{align*}
Consequently,
\begin{align*}
\|\Delta_hf\|_{L^2(\rr^2)}^2
=\frac{\pi}{2}
\lf(1-e^{-(2^{-1}h_1^2+2h_2^2)}\r).
\end{align*}
For every $t\in(0,\fz)$, using the Taylor theorem yields
\begin{align*}
t-\frac{t^2}{2}<1-e^{-t}<t.
\end{align*}
Applying this with $t=2^{-1}h_1^2+2h_2^2$,
and using $(a+b)^2\le2(a^2+b^2)$ for all $a,b\in\rr$, we obtain
\begin{align}\label{eq-Deltahf-upper}
\|\Delta_hf\|_{L^2(\rr^2)}^2
\le\frac{\pi}{4}(h_1^2+4h_2^2)
\end{align}
and
\begin{align}\label{eq-Deltahf-lower}
\|\Delta_hf\|_{L^2(\rr^2)}^2
\ge\frac{\pi}{4}
\left(h_1^2+4h_2^2
-\left(\frac{h_1^2}{2}+2h_2^2\right)^2\right)
\ge
\frac{\pi}{4}
\left(h_1^2+4h_2^2
-\left(\frac{h_1^4}{2}+8h_2^4\right)\right).
\end{align}
Moreover, for $i=1,2$, we have
\begin{align}\label{eq-error}
\lim_{j\to\fz} j^3 \int_{|h|<1/j} \frac{h_i^2}{|h|^3} \lf(\frac{h_1^4}{2}+8h_2^4 \r) dh=0.
\end{align}

Combining \eqref{eq-Deltahf-upper}, \eqref{eq-Deltahf-lower} and \eqref{eq-error},
together with \eqref{eq-energy} and the choice of $\{\varphi_j\}_{j\in\nn}$, we obtain
\begin{align}\label{eq-Ef2k}
\lim_{k\to\fz} \mathscr{E}_{2,\,2}^{1,\,\varphi_{2k}}(f)
&= \frac{3}{4} \lim_{k\to\infty} (2k)^3 \int_{|h|<1/(2k)} \frac{h_1^2\,(h_1^2+4h_2^2)}{|h|^3}\,dh \\
&= \frac{3}{4} \lf(\frac{\pi}{4}+\frac{\pi}{3} \r)
= \frac{7\pi}{16},\notag
\end{align}
whereas
\begin{align}\label{eq-Ef2k+1}
\lim_{k\to\infty} \mathscr{E}_{2,\,2}^{1,\,\varphi_{2k+1}}(f)
&= \frac{3}{4} \lim_{k\to\fz} (2k+1)^3 \int_{|h|<1/(2k+1)} \frac{h_2^2\,(h_1^2+4h_2^2)}{|h|^3}\,dh \\
&= \frac{3}{4} \lf(\frac{\pi}{12}+\pi \r)
= \frac{13\pi}{16}, \notag
\end{align}
where we have used the following identities: for every $a\in(0,\fz)$,
\begin{align*}
\begin{cases}
\dis \int_{|h|<a} \frac{h_1^4}{|h|^3}\,dh
= \int_{|h|<a} \frac{h_2^4}{|h|^3}\,dh
= \frac{\pi a^3}{4}; \\[4mm]
\dis \int_{|h|<a} \frac{h_1^2\,h_2^2}{|h|^3}\,dh
= \frac{\pi a^3}{12};
\end{cases}
\end{align*}
these identities follow from polar coordinates and direct calculations.

In view of \eqref{eq-Ef2k}
and \eqref{eq-Ef2k+1}, the full sequence
$\{\mathscr{E}_{2,\,2}^{1,\,\varphi_j}(f)\}_{j\in\nn}$
does not converge.
Since either $\Gamma$-convergence
or strong $\Gamma$-convergence implies pointwise
convergence (see Section \ref{ss-Gamma=>pw} and Remark \ref{rem-strong-Gamma}), the full sequence of energy functionals
$\{\mathscr{E}_{2,\,2}^{1,\,\varphi_j}\}_{j\in\nn}$ converges in neither the $\Gamma$-sense nor the strong $\Gamma$-sense.
\end{example}

\section{Examples: Besov spaces and logarithmic Besov spaces}\label{sec-example}

The main aim of this section is to apply
Theorems~\ref{thm-MS} and
\ref{thm-BBM-radial}
to classical and logarithmic Besov energies by choosing
appropriate kernel sequences $\{\varphi_j\}_{j\in\nn}$.

\subsection{Pointwise and (strong) $\Gamma$-convergence for Besov spaces}\label{sec-example-Besov}

As an application of Theorem~\ref{thm-MS}, we obtain
the following MS formula for
classical Besov energies, together with the corresponding
$\Gamma$-convergence and strong $\Gamma$-convergence results.

\begin{theorem}[MS formula]\label{Thm-Besov-MS-Gamma}
Let $p\in[1,\fz)$ and $q\in(0,\fz)$. Then
\begin{align*}
\lim_{s\to0_+} \frac{sq}{\omega_{n-1}}\,\left(\|f\|_{\mathring{\Lambda}_{p,\,q}^s(\rn)}\right)^q
= 2^\frac{q}{p}\,\|f\|_{L^p(\rn)}^q
\quad\text{for every}\,\ f\in W^{1,\,p}(\rn).
\end{align*}
Moreover, the same convergence holds in the senses of
$\Gamma$-convergence and strong $\Gamma$-convergence.
\end{theorem}

\begin{proof}
For any $h\in\rn$, set
\begin{align}\label{eq-varphi-Besov-MS}
\varphi_s(h):= \frac{sq}{\omega_{n-1}}\,|h|^{-(n+sq)}.
\end{align}
A direct calculation shows that
\begin{align*}
\lim_{s\to0_+} \int_{|h|\ge1} \varphi_s(h)\,dh
= \lim_{s\to0_+} sq \int_1^\fz t^{-sq}\,\frac{dt}{t}=1
\end{align*}
and that, for any $\delta\in[1,\fz)$,
\begin{align*}
\lim_{s\to0_+} \int_{|h|<\delta} \varphi_s(h) \min\lf\{|h|^q,\,1 \r\} dh
&= \lim_{s\to0_+} sq \int_0^\delta t^{-sq-1}\min\lf\{t^q,\,1 \r\} dt\\
&= \lim_{s\to0_+} sq \lf(\int_0^1 t^{(1-s)q}\,\frac{dt}{t}
+ \int_1^\delta t^{-sq}\,\frac{dt}{t} \r) \\
&= \lim_{s\to0_+} \lf(\frac{s}{1-s}+1-\delta^{-sq} \r) \\
&=0.
\end{align*}
The last two formulae show that the family
$\{\varphi_s\}_{s\in(0,1)}$ satisfies the $\MSq$
condition as $s\to0_+$.
Consequently, Theorem~\ref{thm-MS} yields the
pointwise MS formula, together with its
$\Gamma$-convergence and strong $\Gamma$-convergence versions.
\end{proof}

\begin{remark}
Although Theorem~\ref{thm-MS} is formulated for
countable sequences of kernels, the family
$\{\varphi_s\}_{s\in(0,1)}$ in
\eqref{eq-varphi-Besov-MS} is continuously indexed.
This causes no difficulty, since the desired convergence
as $s\to0_+$ follows by applying Theorem~\ref{thm-MS}
to arbitrary sequences $s_j\to0_+$ and using the
sequential characterization of convergence.
The same argument applies to
Theorems~\ref{Thm-Besov-BBM-Gamma},
\ref{Thm-logBesov-MS-Gamma} and
\ref{Thm-logBesov-BBM-Gamma} below, without further comment.
\end{remark}

An application of Theorem \ref{thm-BBM-radial}, in turn, provides the following BBM formula
for Besov energies, as well as the corresponding
$\Gamma$-convergence and strong $\Gamma$-convergence.

\begin{theorem}[BBM formula]\label{Thm-Besov-BBM-Gamma}
Let $p\in[1,\infty)$, $q\in(0,\infty)$ and
$k\in\mathbb Z_+$. Then
\begin{align*}
\lim_{s\to1_-}
\frac{(1-s)q}{\omega_{n-1}}
\left(\|f\|_{\mathring{\Lambda}_{p,\,q}^{k+s}(\mathbb R^n)}\right)^q
=\mathscr D_{p,\,q}^{k+1}(f)
\qquad\text{for every }f\in W^{k+1,\,p}(\mathbb R^n).
\end{align*}
Moreover, the same convergence holds in the senses of
$\Gamma$-convergence and strong $\Gamma$-convergence.
If $p=1$, then the same conclusions hold
with $W^{k+1,\,1}(\mathbb R^n)$ and
$\mathscr D_{1,\,q}^{k+1}$ replaced by
$\mathrm{BV}^{k+1}(\mathbb R^n)$ and
$\mathscr D_q^{k+1}$, respectively.
\end{theorem}

\begin{proof}
For any $s\in(0,1)$ and $h\in\rn$, set
\begin{align*}
\psi_{s}(h):= \frac{(1-s)q}{\omega_{n-1}}\,|h|^{-n-(k+s)q}.
\end{align*}
By Theorem~\ref{thm-BBM-radial}, it suffices to show that for $m=k+1$
the family $\{\psi_s\}_{s\in(0,1)}$ satisfies the
$\BBMmq$ condition as $s\to1_-$.
Indeed, we have
\begin{align*}
\lim_{s\to1_-} \int_{|h|<1} |h|^{mq}\,\psi_{s}(h)\,dh
= \lim_{s\to1_-} (1-s)q \int_0^1 t^{(1-s)q}\,\frac{dt}{t}=1
\end{align*}
and, for any $\delta\in(0,1]$,
\begin{align*}
&\lim_{s\to1_-} \int_{|h|\ge\delta} \psi_{s}(h) \min\lf\{|h|^{mq},\,1 \r\} dh \\
&\quad= \lim_{s\to1_-} \frac{(1-s)q}{\omega_{n-1}} \int_{|h|\ge\delta} |h|^{-n-(k+s)q} \min\lf\{|h|^{(k+1)q},\,1 \r\} dh \\
&\quad= \lim_{s\to1_-} (1-s)q \lf(\int_\delta^1 t^{(1-s)q}\,\frac{dt}{t}+\int_1^\fz t^{-(k+s)q}\,\frac{dt}{t} \r) \\
&\quad= \lim_{s\to1_-} \lf(1-\delta^{(1-s)q} +\frac{1-s}{k+s} \r) \\
&\quad=0.
\end{align*}
Thus, the $\BBMmq$ condition holds for
$m=k+1$, and the desired conclusions follow
from Theorem~\ref{thm-BBM-radial}.
\end{proof}

\subsection{Pointwise and (strong) $\Gamma$-convergence for logarithmic Besov spaces}\label{sec-example-log-Besov}

Let $p\in[1,\fz)$ and $q\in(0,\fz)$.
An important class of examples in Definition
\ref{def-energy} arises from logarithmic Besov spaces:
\begin{enumerate}[{\rm (i)}]
\item Let $b>-1/q$. In \eqref{eq-energy}, taking $m=1$ and
\begin{align*}
\varphi(h):= |h|^{-n} \lf(\log\frac{e}{|h|\wedge1} \r)^{bq} {\bf 1}_{(0,\,e^{1/(bq+1)})}(|h|),
\end{align*}
we define the \emph{logarithmic Besov space $\mathring{\Lambda}_{p,\,q}^{0,\,b}(\rn)$} in terms of the following (quasi-)seminorm:
\begin{align}\label{eq-log0b}
\|f\|_{\mathring{\Lambda}_{p,\,q}^{0,\,b}(\rn)}
:= \lf(\mathscr{E}_{p,\,q}^\varphi(f) \r)^\frac{1}{q}
= \lf(\int_{|h|<e^{1/(bq+1)}} \|\Delta_hf\|_{L^p(\rn)}^q \lf(\log\frac{e}{|h|\wedge1} \r)^{bq} \frac{dh}{|h|^n} \r)^\frac{1}{q}.
\end{align}
This space describes functions with logarithmic smoothness at the zero-smoothness endpoint. The condition $bq>-1$ ensures that the above (quasi-)seminorm is finite for functions in $C_{\rm c}^\fz(\rn)$.

\item Let $b<-1/q$. In \eqref{eq-energy}, taking $m\in\nn$ and
\begin{align*}
\varphi(h):= |h|^{-(n+mq)} \lf(\log\frac{e}{|h|\wedge1} \r)^{bq},
\end{align*}
we then define the \emph{logarithmic Besov space $\mathring{\Lambda}_{p,\,q}^{m,\,b}(\rn)$}, given by
\begin{align}\label{eq-logmb}
\|f\|_{\mathring{\Lambda}_{p,\,q}^{m,\,b}(\rn)}
:= \lf(\mathscr{E}_{p,\,q}^{m,\,\varphi}(f) \r)^\frac{1}{q}
= \lf(\int_\rn |h|^{-mq}\,\|\Delta_h^mf\|_{L^p(\rn)}^q \lf(\log\frac{e}{|h|\wedge1} \r)^{bq}\,\frac{dh}{|h|^n} \r)^\frac{1}{q}.
\end{align}
This space describes smoothness slightly below
the integer order $m$, with a logarithmic deficit.
The condition $bq<-1$ ensures that
$C_{\rm c}^\fz(\rn)$ is contained in
$\mathring{\Lambda}_{p,\,q}^{m,\,b}(\rn)$.
\end{enumerate}
The study of function spaces with logarithmic smoothness
is closely related to the work of Br\'ezis and Wainger
\cite{Brezis-Wainger-CPDE-1980} and the subsequent
systematic investigations of Edmunds and Haroske
\cite{Haroske-ZAnalA-2000,
Edmunds-Haroske-DMath-1999,
Edmunds-Haroske-JAT-2000}.
These Besov spaces have attracted considerable attention and witnessed notable developments,
including embeddings \cite{Cobos-Dominguez-JMAA-2015,Dominguez-Tikhonov-AMS-2023,Dominguez-Haroske-Tikhonov-JMPA-2020,
Caetano-Gogatishvili-Opic-JAT-2008,Caetano-Gogatishvili-Opic-JAT-2011,Haroske-2007}, Fourier-analytical decompositions
\cite{Cobos-Dominguez-Triebel-JFA-2016,Cobos-Dominguez-JMAA-2015,Dominguez-Tikhonov-AMS-2023}, wavelet and
semigroup characterizations \cite{Cobos-Dominguez-Triebel-JFA-2016,Dominguez-Haroske-Tikhonov-JMPA-2020},
limiting interpolation \cite{Cobos-Dominguez-JAT-2015,Cobos-Dominguez-JFAA-2016,Cobos-Dominguez-JMAA-2015,
DeVore-Riemenschneider-Sharpley-JFA-1979,Liu-Xiao-Zhang-2026-submit} and applications to PDEs \cite{Brezis-Wainger-CPDE-1980,
Zuazua-CRMathAcad-2002,Brue-Nguyen-CMP-2021,Meyer-Seis-IndianaUniv-2024}.

In view of \eqref{eq-log0b} and \eqref{eq-logmb}, by choosing suitable approximation sequences and applying
Theorems \ref{thm-MS} and \ref{thm-BBM-radial}, we obtain the following two theorems.

\begin{theorem}[MS formula]\label{Thm-logBesov-MS-Gamma}
Let $p\in[1,\fz)$ and $q\in(0,\fz)$. Then
\begin{align*}
\lim_{b\to(-1/q)_+}
\frac{bq+1}{\omega_{n-1}}\,\left(\|f\|_{\mathring{\Lambda}_{p,\,q}^{0,\,b}(\rn)}\right)^q
= 2^\frac{q}{p}\,\|f\|_{L^p(\rn)}^q
\quad\text{for every}\,\ f\in W^{1,\,p}(\rn).
\end{align*}
Moreover, the same convergence holds in the senses of
$\Gamma$-convergence and strong $\Gamma$-convergence.
\end{theorem}

\begin{proof}
For any $b>-1/q$ and $h\in\rn$, set
\begin{align*}
\varphi_b(h):=\frac{bq+1}{\omega_{n-1}}\,|h|^{-n} \lf(\log\frac{e}{|h|\wedge1} \r)^{bq} {\bf 1}_{(0,\,e^{1/(bq+1)})}(|h|).
\end{align*}
Note that
\begin{align*}
\mathscr{E}_{p,\,q}^{\varphi_{b}}(f)
= \frac{bq+1}{\omega_{n-1}}\,\left(\|f\|_{\mathring{\Lambda}_{p,\,q}^{0,\,b}(\rn)}\right)^q.
\end{align*}
Moreover, observe that
\begin{align*}
\lim_{b\to(-1/q)_+} \int_{|h|\ge1} \varphi_{b}(h)\,dh
= \lim_{b\to(-1/q)_+} (bq+1) \int_1^{e^{1/(bq+1)}} \lf(\log\frac{e} {t\wedge1} \r)^{bq} \frac{dt}{t}=1
\end{align*}
and that, for any $\delta\in[1,\fz)$,
\begin{align*}
&\lim_{b\to(-1/q)_+} \int_{|h|<\delta} \varphi_{b}(h) \min\lf\{|h|^q,\,1 \r\} dh \\
&\quad= \lim_{b\to(-1/q)_+} (bq+1) \int_0^1 t^q \lf(\log\frac{e}{t} \r)^{bq} \frac{dt}{t}
+ \lim_{b\to(-1/q)_+} (bq+1) \int_1^{\min\{\delta,\,e^{1/(bq+1)}\}} \frac{dt}{t} \\
&\quad=0.
\end{align*}
Thus, the family $\{\varphi_{b}\}_{b>-1/q}$ satisfies the $\MSq$ condition as $b\to(-1/q)_+$.

Consequently, applying Theorem \ref{thm-MS} yields the desired pointwise convergence, $\Gamma$-convergence, and strong
$\Gamma$-convergence for the MS formula associated with logarithmic Besov spaces.
\end{proof}

\begin{theorem}[BBM formula]\label{Thm-logBesov-BBM-Gamma}
Let $p\in[1,\fz)$, $q\in(0,\fz)$ and $m\in\nn$. Then
\begin{align*}
\lim_{b\to(-1/q)_-} \frac{-bq-1}{\omega_{n-1}}\,\lf(\|f\|_{\mathring{\Lambda}_{p,\,q}^{m,\,b}(\rn)} \r)^q
= \mathscr{D}_{p,\,q}^m(f)
\quad\text{for every}\,\ f\in W^{m,\,p}(\rn).
\end{align*}
Moreover, the same convergence holds in the senses of
$\Gamma$-convergence and strong $\Gamma$-convergence.
\end{theorem}

\begin{proof}
For any $b<-1/q$ and $h\in\rn$, set
\begin{align*}
\varphi_{m,\,b}(h)
:= \frac{-bq-1}{\omega_{n-1}}\,|h|^{-(n+mq)} \lf(\log\frac{e}{|h|\wedge1} \r)^{bq}.
\end{align*}
Then,
\begin{align*}
\mathscr{E}_{p,\,q}^{m,\,\varphi_{m,\,b}}(f)
= \frac{-bq-1}{\omega_{n-1}}\,\lf(\|f\|_{\mathring{\Lambda}_{p,\,q}^{m,\,b}(\rn)} \r)^q.
\end{align*}
By Theorem~\ref{thm-BBM-radial}, it suffices to verify that the family $\{\varphi_{m,\,b}\}_{b<-1/q}$
satisfies the $\BBMmq$ condition as $b\to(-1/q)_-$.

Indeed, by polar coordinates, we have
\begin{align*}
\lim_{b\to(-1/q)_-} \int_{|h|<1} |h|^{mq}\,\varphi_{m,\,b}(h)\,dh
= \lim_{b\to(-1/q)_-} (-bq-1) \int_0^1 \lf(\log\frac{e}{t} \r)^{bq} \frac{dt}{t}=1.
\end{align*}
Moreover, for any $\delta\in(0,1]$,
\begin{align*}
&\lim_{b\to(-1/q)_-} \int_{|h|\ge\delta} \varphi_{m,\,b}(h) \min\lf\{|h|^{mq},\,1 \r\} dh \\
&\quad= \lim_{b\to(-1/q)_-} (-bq-1) \int_\delta^1 \lf(\log\frac{e}{t} \r)^{bq} \frac{dt}{t}
+ \lim_{b\to(-1/q)_-} (-bq-1) \int_1^\fz t^{-mq}\,\frac{dt}{t} \\
&\quad=0.
\end{align*}
Thus, the $\BBMmq$ condition holds, and the desired conclusions follow from Theorem~\ref{thm-BBM-radial}.
\end{proof}

\medskip

\noindent{\bf Statements and Declarations.}
The authors declare no competing financial interests. Data availability is not applicable to this paper.


\addcontentsline{toc}{section}{References}

\providecommand{\bysame}{\leavevmode\hbox to3em{\hrulefill}\thinspace}

\medskip

{\small\noindent Liguang Liu \&  Qishun Zhang

\smallskip

\noindent School of Mathematics,
Renmin University of China,
Beijing 100872, People's Republic of China

\smallskip

\noindent{\it E-mails:} \texttt{liuliguang@ruc.edu.cn}; \texttt{qszhang@ruc.edu.cn}
}

\end{document}